\documentclass[12pt]{article}

\usepackage[T1]{fontenc}
\usepackage{amsmath}
\usepackage{amssymb}
\usepackage{latexsym}
\usepackage{textcomp}
\usepackage{diagrams}
\diagramstyle{h=21pt,w=21pt,scriptlabels,midshaft,PostScript=dvips,nohug,shortfall=4pt}
\newarrow{Lig}{=}{=}{=}{=}{=}
\usepackage{graphicx}

\usepackage[noBBpl]{mathpazo}
\renewcommand{\texttt}[1]{{\fontfamily{pcr}\fontseries{m}\fontshape{n}%
\selectfont #1}}
\renewcommand{\tt}{\fontfamily{pcr}\fontseries{m}\fontshape{n}%
\selectfont }

\renewcommand{\d}{\partial}
\newcommand{\upperstar}{^{\raisebox{-0.25ex}[0ex][0ex]{\(\ast\)}}}
\newcommand{\lowerstar}{_{\raisebox{-0.33ex}[-0.5ex][0ex]{\(\ast\)}}}

\newcommand{\upperdot}{^\bullet}

\newcommand{\uppercirc}{^\circ}
\newcommand{\df}{\: {\raisebox{0.255ex}{\normalfont\scriptsize :\!\!}}=}

\newcommand{\isopil}{\stackrel{\raisebox{0.1ex}[0ex][0ex]{\(\sim\)}}%
			{\raisebox{-0.15ex}[0.28ex]{\(\rightarrow\)}}}
\newcommand{\Id}	{\operatorname{Id}}
\newcommand{\fat}[1]{\mathbf{{#1}}}
\providecommand{\qed}{\hspace*{\fill}$\Box$}

\newcommand{\shortsetminus}	{\,\raisebox{1pt}{\ensuremath{\mathbb r}\,}}
\newcommand{\ov}{\overline}

\newcommand{\ct}{\operatorname{ct}}

\usepackage{color}
\definecolor{lightgrey}{rgb}{0.9,0.9,0.9} 

\usepackage{listings}
\providecommand{\overskrift}[1]{\par\noindent\relax{\LARGE #1}\par\bigskip}

\newcommand{\hovedfont}{\normalfont\bfseries}

\usepackage{theorem}
\theorembodyfont{\normalfont\itshape}
\theoremheaderfont{\hovedfont}
	\theoremstyle{change}
\newtheorem{lemma}{Lemma.}[section]
\newtheorem{prop}[lemma]{Proposition.}
\newtheorem{satz}[lemma]{Theorem.}
\newtheorem{cor}[lemma]{Corollary.}
	\theorembodyfont{\normalfont}
\newtheorem{eks}[lemma]{Example.}

\newtheorem{taller}[lemma]{$\!\!$}
\newenvironment{blanko}[1]%
{\begin{taller}{\hovedfont #1}\normalfont}%
{\end{taller}}
{%
\begin{list}{\em Definition. }%
{\setlength{\labelsep}{0mm}\setlength{\leftmargin}{0mm}%
\setlength{\labelwidth}{0mm}\setlength{\listparindent}{\parindent}%
\setlength{\parsep}{\parskip}\setlength{\partopsep}{0mm}}%
\item%
}%
{%
\end{list}%
}
\newenvironment{dem}%
{%
\begin{list}{\em Proof. }%
{\setlength{\labelsep}{0mm}\setlength{\leftmargin}{0mm}%
\setlength{\labelwidth}{0mm}\setlength{\listparindent}{\parindent}%
\setlength{\parsep}{\parskip}\setlength{\partopsep}{0mm}}%
\item%
}%
{%
\qed\end{list}%
}
\newenvironment{dem*}[1]%
{%
\begin{list}{\em #1 }%
{\setlength{\labelsep}{0mm}\setlength{\leftmargin}{0mm}%
\setlength{\labelwidth}{0mm}\setlength{\listparindent}{\parindent}%
\setlength{\parsep}{\parskip}\setlength{\partopsep}{0mm}}%
\item%
}%
{%
\qed\end{list}%
}
{%
\begin{list}{\em Proof. }%
{\setlength{\labelsep}{0mm}\setlength{\leftmargin}{0mm}%
\setlength{\labelwidth}{0mm}\setlength{\listparindent}{\parindent}%
\setlength{\parsep}{\parskip}\setlength{\partopsep}{0mm}}%
\item%
}%
{%
\qed\end{list}%
}
\newenvironment{bevis*}[1]%
{%
\begin{list}{\em #1 }%
{\setlength{\labelsep}{0mm}\setlength{\leftmargin}{0mm}%
\setlength{\labelwidth}{0mm}\setlength{\listparindent}{\parindent}%
\setlength{\parsep}{\parskip}\setlength{\partopsep}{0mm}}%
\item%
}%
{%
\qed\end{list}%
}

\newenvironment{blanko*}[1]%
{%
\begin{list}{\bf {#1} }%
{\setlength{\labelsep}{0mm}\setlength{\leftmargin}{0mm}%
\setlength{\labelwidth}{0mm}\setlength{\listparindent}{\parindent}%
\setlength{\parsep}{\parskip}\setlength{\partopsep}{0mm}}%
\item%
}%
{%
\end{list}%
}

\newcounter{dummycounter}
\newenvironment{punkt-i}%
{%
	\begin{list}%
	{(\roman{dummycounter})}%
	{\usecounter{dummycounter}%
	\setlength{\itemsep}{0em}\setlength{\parsep}{0em}\setlength{\topsep}{0em}%
	\setlength{\itemindent}{0em}\setlength{\labelwidth}{1.8em}%
	\setlength{\labelsep}{0.6em}\setlength{\leftmargin}{2.4em}}%
}%
{\end{list}}

\providecommand{\kat}[1]{\text{\textbf{\textsl{#1}}}}
\newcommand{\Set}{\kat{Set}}

\newcommand{\tr}{\operatorname{tr}}
\newcommand{\const}{\operatorname{const}}

\providecommand{\N}{\mathbb{N}}

\usepackage{mathrsfs}

\newcommand{\EE}{\mathscr{E}}
\newcommand{\Cart}{\kat{Cart}}

\providecommand{\lastUpdate}[1]{#1}

\usepackage{texdraw}

\makeatletter
\renewcommand{\ps@headings}
	{\setlength{\headheight}{13pt}%
	 \setlength{\headsep}{12pt}%
	 \renewcommand{\@oddhead}{\parbox{\textwidth}{%
			\small
			\texttt{\jobname.tex \ \ \ 
			\lastUpdate{2010-02-20 00:59}
			\hfill [\thepage/\pageref{lastpage}]}
			\\ \rule[8pt]{\textwidth}{0.3pt}}%
	 }
	\renewcommand{\@oddfoot}{}
	\renewcommand{\@evenfoot}{}%
}
\makeatother

\begin{document}

\newcommand{\zoom}{\setlength{\unitlength}{1pt}\begin{picture}(30,10)(-5,2)
\put(2.5,5){\circle{5}}
\put(17.5,5){\circle*{3.5}}
\put(5,5){\line(1,0){12.5}}
\end{picture}}

\newcommand{\smallzoom}{\setlength{\unitlength}{0.8pt}\begin{picture}(30,10)(-5,2)
\put(2.5,5){\circle{5}}
\put(17.5,5){\circle*{3.5}}
\put(5,5){\line(1,0){12.5}}
\end{picture}}

\newcommand{\inlineTrivialOpetope}{%
\raisebox{-5pt}{
\begin{texdraw} \move (0 -10) \bsegment
    \move (0 0) \lvec (0 20)
    \move (0 10) \onedot \lcir r:5.5
  \esegment \end{texdraw} }\unskip}

  \newcommand{\inlineOnedotTree}{%
\raisebox{-4pt}{
\begin{texdraw} \linewd 0.5 \scriptsize
      \bsegment
	\move (0 0) \lvec (0 15) \move (0 7.5) \onedot
  \esegment \end{texdraw} } }

\newcommand{\inlineDotlessTree}{%
\raisebox{-4pt}{
\begin{texdraw} \linewd 0.5 \bsegment
    \move (0 0) \lvec (0 15) \move (1 0)
  \esegment \end{texdraw} } }

\newcommand{\polyFunct}[8]{
\begin{diagram}[w=2.5ex,h=3.7ex,tight]
&&#2&&\rTo^{#6}&&#3\\
&\ldTo^{#5}&&&&&&\rdTo^{#7}\\
#1&&&&#8&&&&#4
\end{diagram}
}

\newcommand{\bigleftbrace}[1]{\left\{\raisebox{0pt}[#1pt]{}\right.}
\newcommand{\bigrightbrace}[1]{\left.\raisebox{0pt}[#1pt]{}\right\}}

\everytexdraw{%
	\drawdim pt \linewd 0.5 \textref h:C v:C
	\setunitscale 1
}

\newcommand{\onedot}{
  \bsegment
    \move (0 0) \fcir f:0 r:2
  \esegment
}

\newcommand{\freeEllipsis}[3]{
    \writeps {
      gsave 
      #3 rotate 
    }
    \lellip rx:#1 ry:#2
    \writeps { 
      grestore
    }
}

\newcommand{\freeFillEllipsis}[4]{
    \writeps {
      gsave 
      #3 rotate 
    }
    \fellip f:#4 rx:#1 ry:#2
    \writeps { 
      grestore
    }
}

\newcommand{\Fam}{\kat{Fam}}
\newcommand{\Poly}{\kat{Poly}}

\pagestyle{headings}

\vspace*{24pt}

\begin{center}

\overskrift{Polynomial functors and opetopes}

\bigskip

\noindent
\textsc{Joachim Kock}, 
\textsc{Andr\'e Joyal}, 
\textsc{Michael Batanin}, and
\textsc{Jean-Fran\c cois Mascari}\\

\end{center}

\begin{abstract}
  We give an elementary and direct combinatorial definition of opetopes in terms
  of trees, well-suited for graphical manipulation and explicit computation.  To
  relate our definition to the classical definition, we recast the Baez-Dolan
  slice construction for operads in terms of polynomial monads: our opetopes
  appear naturally as types for polynomial monads obtained by iterating the
  Baez-Dolan construction, starting with the trivial monad.  We show that our
  notion of opetope agrees with Leinster's.  Next we observe a suspension
  operation for opetopes, and define a notion of stable opetopes.  Stable
  opetopes form a least fixpoint for the Baez-Dolan construction.  A final
  section is devoted to example computations, and indicates also how the
  calculus of opetopes is well-suited for machine implementation.
\end{abstract}

\section*{Introduction}

Among a dozen or so existing definitions of weak higher categories,
the opetopic approach is one of the most intriguing, since it is based
on a collection of `shapes' that had not previously been studied: the
opetopes.  
  Opetopes are combinatorial structures parametrising higher-dimensional
  many-in/one-out operations, and can be seen as higher-dimensional
  generalisations of trees.  They are important combinatorial structures
  on their own, `as pervasive in higher-dimensional algebra as simplices
  are in geometry', according to Leinster~\cite[p.216]{Leinster:0305049}.
Opetopes and opetopic higher categories were introduced by Baez and
Dolan in the seminal paper~\cite{Baez-Dolan:9702}, and the theory has
been developed further by
Hermida-Makkai-Power~\cite{Hermida-Makkai-Power:I},
Leinster~\cite{Leinster:0305049}, Cheng~\cite{Cheng:0304277},
\cite{Cheng:0304279}, \cite{Cheng:0304284}, \cite{Cheng:0304287}, and others.  It
is in a sense a theory from scratch, compared to several other
theories of higher categories which build on large bodies of
preexisting machinery and experience, e.g.~simplicial methods.  The
full potential of the opetopic approach may depend on a deeper
understanding of the combinatorics of opetopes.

At the conference on {\em $n$-categories: Foundations and
applications} at the IMA in Minneapolis, June 2004, much time was dedicated to
opetopes, but it became clear that a concise and direct definition of
opetopes was lacking, and that there was no practical way to represent
higher-dimensional opetopes on the blackboard.  In fact, there did not
seem to exist a general method to represent concrete opetopes in any
way, algebraic, graphical, or by machine.%
\footnote{ In fact a method {\em does} exist for algebraic/mechanical
representation: 
Hermida-Makkai-Power~\cite[final section]{Hermida-Makkai-Power:I} explain
how any opetope (there called multitope) in arbitrary dimension can be serialised
into a string of hash signs and stars, with two sorts of brackets. We shall not
go any further into that notation, but just to illustrate its flavour,
here is the representation of the $3$-opetope in \ref{ex-1-2-3}:
$$
\ulcorner\ulcorner\#\urcorner\ulcorner\#\urcorner\ulcorner\#\urcorner 
[\star]\urcorner
\ulcorner\ulcorner\#\urcorner\ulcorner\#\urcorner [\star]\urcorner
\ulcorner\ulcorner\#\urcorner [\star]\urcorner [\#]
[\#]
\ulcorner [\star]\urcorner
\ulcorner\ulcorner\#\urcorner [\star]\urcorner [\#] 
$$
We refer to \cite{Hermida-Makkai-Power:I} for instructions on how to parse this.}
The best definitions are very abstract and not very hands-on:
e.g.~Leinster's definition in terms of iterated free cartesian
monads~\cite{Leinster:0305049}, or the
Hermida-Makkai-Power~\cite{Hermida-Makkai-Power:I} definition of opetopic
sets (there called multitopic sets), followed by a theorem that this
category is a presheaf category, hence characterising a category of
opetopes (there called multitopes).

  As to graphical representations of opetopes in low dimensions, the
  current method is based on a polytope interpretation of opetopes (which
  is at the origin of the terminology: the word `opetope' comes from
  `operation' and `polytope').  Leinster~\cite[\S\,7.4]{Leinster:0305049}
  has constructed a geometric realisation functor which provides support
  for this interpretation, although the polytopes in general cannot be
  piece-wise linear objects in Euclidean space.  Moreover, geometrical
  objects in dimension higher than $3$ are inherently difficult to
  represent graphically, and currently one resorts to Lego-like drawings in
  which the individual faces of the polytopes are drawn separately, with
  small arrows as a recipe to indicate how they are supposed to fit
  together.

The goal of this paper is to come closer to the combinatorics.  Our initial
idea was to represent an opetope as a tree with some circles, which we now
call {\em constellations}.  This works in dimension $4$ (cf.~\ref{dim4}
below), but it does not seem to be sufficient to capture the possible
opetopes in dimension $5$ and higher.  Pursuing the idea, what we
eventually found was a representation in terms of a sequence of trees with
circles, and in fact it is basically the notion of metatree originally
proposed by Baez and Dolan.  
That notion was never really developed,
though: in the original paper~\cite{Baez-Dolan:9702} the claim that
metatrees could express opetopes was not really substantiated, and in the
subsequent literature there seems to be no mention of the metatree notion.
The presence of circles makes a conceptual
difference, and it also reveals a certain shortcoming in the original
notion of metatree, related to units (cf.~\ref{const v meta}).

We hasten to point out that our notion of opetope coincides with the notion due
to Leinster~\cite{Leinster:0305049} (cf.~the explicit comparison culminating in
Theorem~\ref{ours=leinsters}), not with the original Baez-Dolan definition: we
work consistently with non-planar trees, which means our opetopes are
`un-ordered' like abstract geometric objects, whereas the original Baez-Dolan
opetopes come equipped with an ordering of their faces.  In our version, the
planar aspect is only a particular feature of low dimensional opetopes.

While our opetopes agree with Leinster's, the description we provide is
completely elementary and does not even make reference to category theory.
We think that our description can serve as the famous
\mbox{`$5$-minute definition'} that was previously missing, and that it can
provide a convenient tool for communicating opetopical ideas.
We also indicate how our approach is well-suited for machine 
manipulation.

Opetopes were introduced to parametrise higher-dimensional substitution
operations.  Surprisingly, opetopes arise also in another way, namely from
computads and higher-dimensional pasting theory, and we wish to mention
that a very different combinatorial approach has been developed in this
setting by Palm~\cite{Palm}.  A {\em computad} is a strict
$\omega$-category which is dimension-wise free.  This notion was devised by
Street~\cite{Street:computads} as a tool for describing higher-dimensional
compositions in strict $n$-categories.  In the works of
Johnson~\cite{Johnson:n-pasting} and Power~\cite{Power:2-pasting},
\cite{Power:n-pasting}, different combinatorial and topological
representations of computads (called pasting schemes) were given, starting
from B\'enabou's pasting diagrams for $2$-categories and the dual
graphical language of string diagrams.  The subtleties encountered are
related with the fact that the category of computads is not a presheaf
category.  A computad is called {\em many-to-one} if the codomain of every
indeterminate in dimension $k+1$ is itself an indeterminate (in dimension
$k$).  Harnik, Makkai and Zawadowski~\cite{Harnik-Makkai-Zawadowski}
established an equivalence of categories between many-to-one computads and
multitopic sets.  In particular, the category of many-to-one computads {\em
is} a presheaf category.  Palm~\cite{Palm} has given a purely combinatorial
description of this presheaf category.  He introduces a notion of {\em
dendrotopes}, certain decorated Hasse diagrams, and shows that dendrotopic
sets (their presheaves) are equivalent to many-to-one computads.  Hence, by
the theorems of Harnik-Makkai-Zawadowski and Hermida-Makkai-Power,
dendrotopes should correspond to opetopes.  However, a direct combinatorial
comparison has not been given at this time.

\bigskip

Let us briefly outline the organisation of the exposition.  In the
first section we give the definition of opetopes in a direct
combinatorial way, without reference to category theory.  The crucial
ingredient is the correspondence between {\em non-planar} trees and
nestings of circles: an opetope is merely a sequence of such
correspondences, with an initial condition. 
We give the definition in two steps: first the elementary `5-minute
definition' with examples, then we develop the involved notions of
trees and constellations more formally and compare with Baez-Dolan
metatrees.  It is possible
to jump directly from the `5-minute definition' 
to Section~\ref{Sec:calculus}, where the same
elementary and hands-on approach is pursued to describe in detail how
to compute sources and targets of opetopes, and how to compose them.
However, such a reading would ignore the theoretical justification for
the definitions and constructions.

In Section~\ref{Sec:poly} we review some basic facts about polynomial functors,
notably their graphical interpretation which is the key point to relate the
formal constructions with explicit combinatorics.

Section~\ref{Sec:BD} forms the theoretical heart of this work: we give an easy
account of the Baez-Dolan slice construction in the setting of polynomial
monads.  From the graphical description of polynomial functors we see that the
Baez-Dolan construction is about certain decorated trees.  The double Baez-Dolan
construction gives trees decorated with trees, subject to complicated
compatibility conditions.  We show that these compatibility conditions are
completely encoded by drawing circles in trees.  Iterating the Baez-Dolan
construction involves the correspondence between trees and nestings, and it
readily follows (Theorem~\ref{op=op}) that the opetopes defined in
Section~\ref{Sec:opetopes} arise precisely as types for the polynomial
monads produced by
iterating the Baez-Dolan construction, starting from the trivial monad.
We compare the polynomial Baez-Dolan construction with
Leinster's version of the Baez-Dolan construction, and conclude 
(Theorem~\ref{ours=leinsters}) that our
notion of opetope agrees with Leinster's~\cite{Leinster:0305049}.

In the short Section~\ref{Sec:stable}, we observe a suspension
operation for opetopes, and define a notion of stable opetopes.  The
stable opetopes also form a polynomial monad, and we show this is the
least fixpoint for the Baez-Dolan construction (for pointed monads).

In Section~\ref{Sec:calculus}, we show by way of examples how the 
calculus of opetopes works in practice: we are concerned with 
computing sources and target of opetopes, and with composing them.
In the Appendix we briefly describe a machine implementation of the
`calculus of opetopes' based on XML, including a mechanism for automated
graphical output.


\begin{blanko*}{Acknowledgements.}
  We are grateful to John Baez and Peter May for organising the very inspiring
  workshop in Minneapolis, and to Eugenia Cheng and Michael Makkai for patiently
  telling us about opetopes at that occasion.  We are grateful to editors and
  referees for their attentiveness.  We thank our respective financing
  institutions: the research of A.~J.~was supported by the NSERC; M.~B.~was
  supported by the Australian Research Council; J.-F.~M.~was supported by a
  grant from the CNR in the framework of an IMA-CNR collaboration; J.~K.,
  currently supported by grants MTM2006-11391 and MTM2007-63277 of the Spanish
  Ministry of Education and Science, was previously supported by a grant from
  the CIRGET at the UQAM and insisted on including this opetope drawing, as
  expression of his gratitude and admiration:

\begin{center}\begin{texdraw}
    \setunitscale 0.8

  \linewd 1.2
\move (0 12) \lvec (35 12) 
\lvec (30 25) \lvec (60 45) 
\lvec (50 52) \lvec (60 72)  
\lvec (40 65) \lvec (35 75) 
\lvec (15 55) \lvec (25 95) 
\lvec (10 90) \lvec (0 110)

\move (0 12)   \lvec (-35 12) 
\lvec (-30 25) \lvec (-60 45) 
\lvec (-50 52) \lvec (-60 72)  
\lvec (-40 65) \lvec (-35 75) 
\lvec (-15 55) \lvec (-25 95) 
\lvec (-10 90) \lvec (0 110)

\move (0 -16) \lvec (0 18.5) \onedot

\linewd 0.2
\move (30 25) \lvec (-30 25)
\lvec (-15 55) \lvec (15 55) \lvec (30 25)
\move (30 25) \lvec (50 52) \lvec (40 65) \lvec (15 55)
\lvec (10 90) \lvec (-10 90) \lvec (-15 55)
\lvec (-40 65) \lvec (-50 52)
\lvec (-30 25)

\move (0 18.5)
\lvec (0 38) \onedot \lvec (0 75) \onedot
\lvec (0 96) \onedot \lvec (10 105)
\move (0 96) \lvec (-10 105)

\move (0 18.5) \lvec (-40 18.5)
\move (0 18.5) \lvec (40 18.5)
\move (0 38) \lvec (35 50) \onedot
\move (0 38) \lvec (-35 50) \onedot

\move (0 75) \lvec (17 85) \onedot
\move (0 75) \lvec (-17 85) \onedot
\move (17 85) \lvec (17 98) \move (17 85) \lvec (27 80)
\move (-17 85) \lvec (-17 98) \move (-17 85) \lvec (-27 80)

\move (35 50) \lvec (32 66) \onedot
\lvec (25 73) \move (32 66) \lvec (43 73)
\move (35 50) \lvec (49 62) \onedot \lvec ( 48 72)
\move (49 62) \lvec (58 58)
\move (35 50) \lvec (50 44) \onedot \lvec (59 53)
\move (50 44) \lvec (53 35)

\move (-35 50) \lvec (-32 66) \onedot
\lvec (-25 73) \move (-32 66) \lvec (-43 73)
\move (-35 50) \lvec (-49 62) \onedot \lvec ( -48 72)
\move (-49 62) \lvec (-58 58)
\move (-35 50) \lvec (-50 44) \onedot \lvec (-59 53)
\move (-50 44) \lvec (-53 35)

\end{texdraw}

\end{center}
\end{blanko*}

\section{Opetopes}

\label{Sec:opetopes}

We first give the quick definition of opetope, through the
notions of tree, constellation, and zoom.
Afterwards we develop these notions more carefully.

\subsection{The `5-minute definition' of opetope}

\begin{blanko}{Trees.}
  The fundamental concept is that of a tree.  Our trees are non-planar
  finite rooted trees with boundary: they have any number of input edges
  (called leaves), and have precisely one output edge (called the root
  edge) always drawn at the bottom.  There is a partial order in which the
  root is the maximal element and the leaves are minimal elements. 
  The following drawings
  should suffice to exemplify trees, but beware that the planar aspect
  inherent in a drawing should be disregarded:
\begin{center}\begin{texdraw}
  \linewd 0.5 \footnotesize
  \move (-50 0)
  \bsegment
    \move (0 0) \lvec (0 30)
  \esegment
  
  \move (0 0)
  \bsegment
    \move (0 0) \lvec (0 18) \onedot
  \esegment
  
  \move (50 0)
  \bsegment
    \move (0 0) \lvec (0 36)
    \move (0 18) \onedot
  \esegment
  
  \move (105 0)
  \bsegment
    \move (0 0) \lvec (0 15) \onedot
    \lvec (-5 33) \onedot
    \move (0 15) \lvec (-12 28) \onedot
    \move (0 15) \lvec (4 43)
    \move (0 15) \lvec (12 40)
  \esegment
  
  \move ( 170 0)
  \bsegment
    \move (0 0) \lvec (0 18) \onedot
    \lvec (-6 32) \onedot
    \lvec (-12 57)
    \move (0 18) \lvec (4 40) \onedot
    \lvec (20 50) \onedot
    \lvec (15 65)
    \move (20 50) \lvec (25 65)
    \move (4 40) \lvec (9 54) \onedot
    \move (4 40) \lvec (-4 61)
  \esegment
\end{texdraw}\end{center}
A formal definition of tree is given in \ref{tree-formal}.
An alternative formalism is developed in \cite{Kock:0807}.
\end{blanko}

\begin{blanko}{Nestings.}\label{nest}
  Another graphical representation of the same structure is given in terms
  of nested circles in the plane.  We prefer to talk about nested spheres
  in space to avoid any idea of planarity when in a moment we combine the
  notion with trees.  A {\em nesting} is a finite collection of
  non-intersecting spheres and dots, which either consists of a single dot
  (and no spheres) or has one outer sphere, containing all the other
  spheres and dots.

  The dots of a nesting correspond to the leaves of the tree.  The
  outer sphere corresponds to the root edge of the tree, and the special
  case of a nesting which consists solely of one dot corresponds to
  the dotless tree.  The partial order is simply inclusion.

  The following drawings of nestings correspond exactly to the
  five trees drawn above.
  
\begin{center}\begin{texdraw}
  \linewd 0.5 \footnotesize

  \move (-45 0)
  \bsegment
    \move (0 15) \onedot
  \esegment
    
  \move (0 0)    
  \bsegment
    \move (0 17) \lcir r:12
  \esegment

  \move (50 0)
  \bsegment
    \move (0 17) \onedot \lcir r:12
  \esegment

  \move (110 0)
  \bsegment
    \move (0 21) \lcir r:20
    \move (-5 13) \lcir r:6
    \move (5 26) \lcir r:6
    \move (-6 30) \onedot
    \move (10 13) \onedot
  \esegment
    
  \move ( 190 0)
  \bsegment
    \move (0 30) \lcir r:30
    \move (-15 15) \lcir r:7
    \move (8 32) \lcir r:19
    \move (16 38) \lcir r:5
    \move (7 24) \lcir r:9
    \move (-16 14) \onedot
    \move (4 25) \onedot
    \move (10 23) \onedot
    \move (3 40) \onedot
  \esegment
\end{texdraw}\end{center}
\end{blanko}

\begin{blanko}{Correspondences.}\label{corr}
  A {\em correspondence} between a nesting $S$ and a tree $T$ 
  consists of specified bijections
  $$\begin{array}{rrr}
    \text{dots}(S) &\leftrightarrow  & \text{leaves}(T) \\
    \text{spheres}(S) & \leftrightarrow & \text{dots}(T) 
  \end{array}$$
  respecting the partial orders.
Here is a typical picture:
\begin{center}
\begin{texdraw}
  \bsegment \linewd 0.3 \footnotesize
    \move (0 40) \lellip rx:35 ry:40
    \htext (24 18) {$a$}
    \move (-12 23) \lcir r:13
    \htext (-18 28) {$b$}
    \move (10 52) \lcir r:17 \htext (-2 56) {$c$}
    \move (14 56) \lcir r:7 \htext (14 57) {$d$}
    \move (-18 52) \onedot \htext (-23 52) {$e$}
    \move (-7 20) \onedot \htext (-11 20) {$f$}
    \move (15 43) \onedot \htext (10 43) {$g$}
  \esegment

  \move (145 0)
  \bsegment \linewd 0.5 \footnotesize
    \move (0 0) \lvec (0 15) \onedot \lvec (-5 70) 
    \move (0 15) \lvec (-15 30) \onedot \lvec (-20 65)
    \move (0 15) \lvec (20 38) \onedot
    \lvec (12 50) \onedot
    \move (20 38) \lvec (30 65)
    \htext (-5 12) {$a$}
    \htext (-5 75) {$e$}
    \htext (-20 30) {$b$}
    \htext (-20 71) {$f$}
    \htext (25 38) {$c$}
    \htext (12 57) {$d$}
    \htext (32 70) {$g$}
  \esegment
  \end{texdraw}\end{center}
 The bijections are indicated by the labels $a,b,c,d,e,f,g$.
\end{blanko}

\begin{blanko}{Constellations.}\label{const}
  A {\em constellation} is a superposition of a tree with a nesting
  with common set of dots, and such that
  each sphere cuts a subtree. Here is an example:
  \begin{center}
    \begin{texdraw}
    \bsegment 
    \linewd 0.5 
    \move (0 -5) \lvec (0 20) \onedot
    \lvec (-15 30) \onedot
    \move (0 20) 
    \lvec (20 38) \onedot
    \lvec (28 52) \onedot
    \lvec (28 90)

    \move (20 38) \lvec (-20 60)
    \onedot \lvec (-30 85)
    \move (-20 60) \lvec (-10 96)
    
	  \linewd 0.3
		  \move (-9 25) \freeEllipsis{12}{20}{54}
		  \move (24 45) \lellip rx:12 ry:18
		  \move ( -5 52) \lcir r:6
		  \move (0 46) \lcir r:44
    \esegment

    \end{texdraw}
  \end{center}
  More precisely, it is a configuration $C$ of edges, dots, and
  spheres, such that 
  \begin{punkt-i}
    \item edges and dots form a tree
  (called the {\em underlying tree} of $C$),
    
  \item dots and
  spheres form a nesting (the {\em underlying nesting}
  of $C$),
  
  \item for each sphere, the edges and dots contained
  in it form a tree again.
  
  \end{punkt-i}
  
  \medskip
  
A purely combinatorial definition of constellation is given in 
\ref{const-formal}.
  
  \bigskip

  Let us briefly take a look at some degenerate examples.
  In a constellation without a sphere, the underlying
  nesting is necessarily a single dot.  Hence the possibilities in
  this case are exhausted by the set of trees with only one dot:
\begin{center}\begin{texdraw}
  \linewd 0.5 \footnotesize
  
  \move (0 0)
  \bsegment
    \move (0 0) \lvec (0 15) \onedot
  \esegment
  
  \move (40 0)
  \bsegment
    \move (0 0) \lvec (0 27)
    \move (0 15) \onedot
  \esegment
  
  \move (85 0)
  \bsegment
    \move (0 0) \lvec (0 15) \onedot
    \lvec (-6 27)
    \move (0 15) \lvec (6 27)
  \esegment
  
  \move (135 0)
  \bsegment
    \move (0 0) \lvec (0 15) \onedot
    \lvec (-8 27)
    \move (0 15) \lvec (0 30)
    \move (0 15) \lvec (8 27)
  \esegment
  
  \htext (185 15) {etc.}
\end{texdraw}\end{center}
  
In a constellation without dots, the underlying tree must be a single edge.
There must be an outer sphere, so such constellations may look like
these examples:
\begin{center}\begin{texdraw}
  \linewd 0.5 \footnotesize
  
  \move (0 0)
  \bsegment
    \move (0 0) \lvec (0 34) 
    \move (0 17) \lcir r:10
  \esegment
  
  \move (60 0)
  \bsegment
    \move (0 0) \lvec (0 44)
    \move (0 16) \lcir r:4
    \move (0 28) \lcir r:4
    \move (0 22) \lcir r:16
  \esegment
\end{texdraw}\end{center}
Note that every sphere must contain a segment of a line, since there 
is no such thing as the empty tree.

Finally, we draw a few examples of constellations without leaves:
\begin{center}\begin{texdraw}
  \linewd 0.5 \footnotesize
  
  \move (0 0)
  \bsegment
    \move (0 0) \lvec (0 15) \onedot
  \esegment
  
  \move (40 0)
  \bsegment
    \move (0 0) \lvec (0 25) \onedot
     \move (0 15) \lcir r:4
     \move (0 19) \lcir r:14
  \esegment
  
  \move (95 0)
  \bsegment
    \move (0 0) \lvec (0 12) \onedot 
    \lvec (-8 21) \onedot \lcir r:6
    \move (0 12) \lvec (4 27) \onedot
    \move (0 12) \lvec (10 22) \onedot
 \move (0 22) \lcir r:17
 \esegment
\end{texdraw}\end{center}

In \ref{Thm:slicetwice} it is shown that constellations represent, in a
precise sense, trees of trees, which is the reason for their importance.
We want to iterate the idea of trees of trees by repeating the step of
drawing spheres.  To do this, we shift the nesting to a tree and
iterate.  In our terminology, we zoom:
\end{blanko}

\begin{blanko}{Zooms.}\label{zoom}
  A {\em zoom} from constellation $A$ to constellation $B$, written
  $$
  A \zoom B ,
  $$
  is a correspondence between the underlying nesting of $A$ and
  the underlying tree of $B$.  In other words, there are specified two
  bijections: $$\begin{array}{rrr} 
    \text{dots}(A) &\leftrightarrow  & \text{leaves}(B) \\
  \text{spheres}(A) & \leftrightarrow & \text{dots}(B)
  \end{array}$$
  respecting the partial orders.  

  Here is an example:
  \begin{center}\begin{texdraw}
  \bsegment \linewd 0.5 \scriptsize
    \move (0 -10) \lvec (0 20) \onedot
    \lvec ( -80 40)
    \move (0 20) \lvec (-20 70) \onedot
    \lvec (-50 80) \onedot
    \lvec (-83 110)
    \move (-20 70) \lvec (-40 170)
    \move (-20 70) \lvec (0 110) \onedot
    \move (0 20) \lvec (15 70) \onedot
    \move (0 20) \lvec (40 60) \onedot
    \lvec (25 130) \onedot
    \lvec (0 175)
    \move (25 130) \lvec (70 140) 
    \move (40 60) \lvec (83 100)
  \linewd 0.3
    \move (0 80) \lellip rx:70 ry:80
    \move (-28 115) \lcir r:7 \htext (-25 115) {$7$}
    \move (-31 132) \lcir r:7 \htext (-29 132) {$6$}
   
    \move (-27 123) \freeEllipsis{16}{21}{16} \htext (-15 126) {$5$}

    \move (-36 73) \freeEllipsis{18}{25}{70} \htext (-40 61) {$4$}

    \move (22 38) \freeEllipsis{36}{14}{45}
    \move (18 45) \freeEllipsis{42}{32}{52} \htext (-1 65) {$2$}

    \htext (61 110) {$1$}
    \htext (31 35) {$3$}
    \htext (7 20) {$8$}
    \htext (18 130) {$11$}
    \htext (41 54) {$9$}
    \htext (22 70) {$10$}
    \htext (-25 65) {$12$}
    \htext (-52 74) {$13$}
    \htext (7 110) {$14$}
  \esegment

  \htext (90 15) {\large $A$}
  \htext (150 15) {\large $B$}
  \htext (120 15) {$\zoom$}

  \move (250 5)
  \bsegment \linewd 0.5 \scriptsize
    \move (0 -20) \lvec (0 20) \onedot
    \lvec ( -35 30) \onedot
    \lvec (-60 70) \onedot
    \lvec (-90 90)
    \move (-60 70) \lvec (-75 125)
    \move (-35 30) \lvec (-55 145)    
    \move (0 20) \lvec (-30 152) 
    \move (0 20) \lvec (0 155)
    \move (0 20) \lvec (15 70) \onedot
    \lvec (50 140)
    \move (15 70) \lvec ( 80 120)
    \move ( 0 20) \lvec (50 40) \onedot
    \lvec (40 68) \onedot
    \move (50 40) \lvec (62 64) \onedot
  \linewd 0.3
    \move (-32 33) \lcir r:9
    \move (35 110) \lcir r:8
    \move (-10 64) \lcir r:5
    \move (-8 46) \freeEllipsis{42}{35}{30}
    \move (58 52) \freeEllipsis{26}{12}{61}
    \move (0 65) \lellip rx:80 ry:70
    \htext (-5 16) {$1$}
    \htext (-30 35) {$2$}
    \htext (-65 67) {$3$}
    \htext (20 66) {$4$}
    \htext (56 40) {$5$}
    \htext (68 64) {$7$}
    \htext (40 75) {$6$}
    \htext (-94 91) {$9$}
    \htext (-77 129) {$8$}
    \htext (-56 151) {$10$}
    \htext (-30 158) {$14$}
    \htext (0 161) {$11$}
    \htext (86 123) {$13$}
    \htext (52 145) {$12$}
  \esegment
\end{texdraw}\end{center}
The bijections are indicated with numbers.

We also wish to exhibit the two most degenerate zooms:
\begin{center}\begin{texdraw}
  \bsegment \linewd 0.5 \scriptsize
    \move (0 5) \lvec (0 20) \onedot
    \htext (-7 20) {$x$}
  \esegment

  \htext (40 15) {$\zoom$}
  \move (90 0)
  \bsegment \linewd 0.5 \scriptsize
    \move (0 0) \lvec (0 30) 
  \linewd 0.3
    \move (0 15) \lcir r:10.5
    \htext (0 35) {$x$}
  \esegment

  \move (200 0)
  \bsegment
    \bsegment \linewd 0.5 \scriptsize
      \move (0 0) \lvec (0 30) 
    \linewd 0.3
      \move (0 15) \lcir r:11
      \htext (7.5 15) {$x$}
    \esegment
    \htext (45 15) {$\zoom$}
    \move (85 0)
    \bsegment \linewd 0.5 \scriptsize
      \move (0 5) \lvec (0 20) \onedot
      \htext (7 20) {$x$}
    \esegment
  \esegment
\end{texdraw}\end{center}  
\end{blanko}

\begin{blanko}{Zoom complexes.}\label{zoom-complex}
  A {\em zoom complex} of degree $n\geq 0$ is a sequence of zooms
  $$
  X_0 \zoom X_1 \zoom X_2 \zoom X_3 \zoom \dots \zoom X_n .
  $$
\end{blanko}

\begin{blanko}{Opetopes.}\label{def}
  An {\em opetope} of dimension $n\geq 0$ is defined to be a zoom complex
  $X$ of degree $n$ starting like this:
\begin{equation}\label{fig:def}\begin{texdraw}
  \linewd 0.5 

  \bsegment
    \move (0 0) \lvec (0 40)
    \move (0 20) \onedot \linewd 0.3 \lcir r:10
    \htext (0 -20) {\normalsize $X_0$}
  \esegment
  \htext (40 -20) {$\zoom$}

  \move (80 0)
  \bsegment
  \linewd 0.5
    \move (0 0) \lvec (0 40)
    \move (0 20) \onedot \linewd 0.3 \lcir r:10
    \htext (0 -20) {\normalsize $X_1$}
  \esegment
  \htext (125 -20) {$\zoom$}

  \move (170 0)
  \bsegment
    \move (0 -5)
    \lvec (0 51) 
    \move (0 23) \onedot
    \linewd 0.3
    \lcir r:5
    \lcir r:10
    \lcir r:15
    \lcir r:20
    \htext (0 -20) {$X_2$}
    \esegment
\end{texdraw}\end{equation}
   Here, $X_0$ and $X_1$ are exactly as drawn, while $X_2$ is
   described verbally as having one dot and one leaf (necessary in
   order to be in zoom relation with $X_1$), and having any finite
   number of linearly nested spheres.  (We consider two opetopes the
   same if they only differ by the names of the involved elements.)
\end{blanko}

\begin{blanko}{Remark.}
  This definition of opetope should be attributed to Baez and
  Dolan~\cite{Baez-Dolan:9702} who introduced the notion of opetope in
  terms of a slice construction for symmetric operads (a polynomial
  analogue of which we shall call the {\em Baez-Dolan construction}
  (Section~\ref{Sec:BD})), and offered an alternative description in terms
  of sequences of trees called metatrees.  Definition~\ref{def} features
  important adjustments to the Baez-Dolan notion of metatree, as we shall
  explain in \ref{const v meta}
\end{blanko}

\begin{blanko}{Examples.}\label{ex-1-2-3}
  A $0$-opetope is the zoom complex
  \inlineTrivialOpetope\ (there is only one such), and a $1$-opetope is
  the zoom complex \inlineTrivialOpetope $\smallzoom\!\!$
  \inlineTrivialOpetope\ (again there is only one such).  The
  $2$-opetopes are in bijection with the natural numbers, counting the
  linearly nested spheres in $X_2$.  
  
  For $n\geq 3$, there are no
  restrictions on the constellations $X_n$, except to be in zoom
  relation with $X_{n-1}$.  
  For example, if there are $n$ spheres in
  $X_2$, then the zoom condition forces $X_3$ to be a straight line
  with $n$ dots on (and the bijection between spheres and dots is
  uniquely determined since the linear nesting of the spheres in $X_2$
  must correspond to the linear arrangement of the dots in $X_3$), and
  any nesting can be drawn on top of that. Here is an example:

    \begin{center}
  \begin{texdraw}
    
    \move (-165 -25)
      \bsegment
    \move (0 0) \lvec (0 40)
    \move (0 20) \onedot \lcir r:10
    \htext (0 -30) {\normalsize $X_0$}
  \esegment
  \htext (-125 -5) {$\zoom$}

  \move (-80 -25)
  \bsegment
    \move (0 0) \lvec (0 40)
    \move (0 20) \onedot \lcir r:10
    \htext (0 -30) {\normalsize $X_1$}
  \esegment
  \htext (-40 -5) {$\zoom$}

    \move (5 -5)
     \bsegment
  
  \move (0 -30) \lvec (0 30) 
  \move (0 0) \onedot
  \linewd 0.3
  \lcir r:5
  \lcir r:10
  \lcir r:15
    \htext (0 -50) {\normalsize $X_2$}
  \esegment

 \htext (52 -5) {$\zoom$}

 \move (105 0)

   \bsegment
    \linewd 0.5

  \move (0 -35)
  \lvec (0 35) 
  \move (0 -20) \onedot
  \move (0 8) \onedot
  \move (0 20) \onedot
    
  \linewd 0.3
  
  \move (0 15) \lellip rx:10 ry:14
  \move (0 2) \lellip rx:18 ry:30
  \move (0 -20)
    \lcir r:5
  \move (0 -7)
    \lcir r:4
  \move (0 20)
    \lcir r:5
    \htext (0 -55) {\normalsize $X_3$}

  \esegment
\move (0 -70)
\end{texdraw}\end{center}

Clearly the information encoded in $X_0$, $X_1$ and $X_2$ is redundant, and
a $3$-opetope is completely specified by a $X_3$ of this form: a line with
dots and `spheres'.  This is equivalent to specifying a {\em
planar} tree.  The planarity comes about because there is a line organising
the dots in $X_3$, which in turn is a consequence of the linear nesting of
the spheres in $X_2$.
Here is the planar tree corresponding to the $3$-opetope above:
\begin{center}\begin{texdraw}
	\bsegment
  \linewd 0.5 \scriptsize
  
  \move (0 5)
  \lvec (0 20) \onedot
  \lvec (-12 32) \onedot
  \lvec (-15 45) \onedot
  \lvec (-15 55)
  \move (-12 32)
  \lvec (0 55)
  \move (0 20) \lvec (0 32) \onedot
  \move (0 20) \lvec (10 38) \onedot
  \lvec (10 55) \move (0 60)
  \esegment
\end{texdraw}\end{center}


\noindent and here is how this $3$-opetope would be represented in the
polytope style, as in Leinster's book~\cite{Leinster:0305049} and in the
work of Cheng:
      \begin{center}\begin{texdraw}
  \linewd 0.4 \setunitscale 1.4
	\bsegment
  \move (0 0)
  \lvec (40 0)
  \lvec (20 15)
  \lvec (0 0)
  \lvec (0 20)
  \lvec (20 15)
  \clvec (12 0)(28 0)(20 15)
  \clvec (30 15)(40 10)(40 0)
  \move (0 0)
  \clvec (-5 5)(-5 15)(0 20)
  \esegment

  \move (80 0)
  	\bsegment
  \move (0 0)
  \lvec (40 0)
  \lvec (22 17)
  \lvec (0 20)
  \lvec (0 0)
  \esegment

  \move (57 10)
  \bsegment
  \linewd 0.2
  \move (0 0) \lvec (6 0)
  \lvec (3 3) \move (6 0)
  \lvec (3 -3)
  \move (0 1.5) \lvec (4.5 1.5)
  \move (0 -1.5) \lvec (4.5 -1.5)
  \esegment
\end{texdraw}\end{center}

\end{blanko}

\begin{blanko}{Remark.}
  The two-step initial condition in the definition of opetope may look strange,
  and in any case the first two constellations are redundant in terms of
  information.  (As we just saw, for $n\geq 3$ also $X_2$ is redundant, since
  the configuration of dots in $X_3$ completely determines $X_2$.)  The
  justifications for including $X_0$ and $X_1$ are first of all to cover also
  dimension $0$ and $1$ in an uniform way, and make the opetope dimension match
  the degree of the complex.  Second, those leading
  \raisebox{-2pt}{\inlineTrivialOpetope} will play a key role in the notion of
  stable opetopes in \ref{suspension}.  From the theoretical viewpoint, which we
  take up in the next section, the point is that $X_0$ and $X_1$ represent the
  trivial polynomial functor (the identity functor on $\Set$), from which
  iterated application of the Baez-Dolan construction (\ref{BD}) will generate all
  the opetopes in higher dimension, cf.~Theorem~\ref{op=op}.
  The extra condition imposed on $X_2$ (the
  linear nesting of the spheres) is also explained by that construction.  The
  fact that there are no extra conditions on $X_n$ for $n\geq 3$ expresses a
  remarkable feature of the double Baez-Dolan construction, at the heart of this
  paper, namely that the double Baez-Dolan construction generates
  constellations, cf.~Theorem~\ref{Thm:slicetwice}.
\end{blanko}

\begin{eks}
  \label{dim4}
  A $4$-opetope is a zoom complex of degree $4$ like this example:
  \begin{center}
  \begin{texdraw}
    
    \move (-165 -25)
      \bsegment
    \move (0 0) \lvec (0 40)
    \move (0 20) \onedot \lcir r:10
    \htext (0 -35) {\normalsize $Y_0$}
  \esegment
  \htext (-125 -5) {$\zoom$}

  \move (-80 -25)
  \bsegment
    \move (0 0) \lvec (0 40)
    \move (0 20) \onedot \lcir r:10
    \htext (0 -35) {\normalsize $Y_1$}
  \esegment
  \htext (-40 -5) {$\zoom$}

    \move (5 -5)
    \bsegment
 \scriptsize
 
 \move (0 -25) \lvec (0 25) 
 \move (0 0) \onedot
 \lcir r:6
 \lcir r:13
 \htext (9 0) {$b$}
 \htext (14 -10) {$a$}
   \htext (0 -55) {\normalsize $Y_2$}
 \esegment
 
 \htext (52 -5) {$\zoom$}

 \move (105 0)

 \bsegment
 \scriptsize
    \move (0 -45) \lvec (0 45) 
   \move (0 -13) \onedot \htext (6 -13) {$a$}
   \move (0 0)      \lcir r:5 \htext (9 0) {$r$}
   \move (0 -25)      \lcir r:5 \htext (9 -25) {$q$}

   \move (0 20) \onedot \lcir r:11 \htext (5 20) {$b$}
   \htext (14 20) {$s$}
   \htext (24 -20) {$p$}
   \move (0 0) \lellip rx:22 ry:40
     \htext (0 -60) {\normalsize $Y_3$}
 \esegment

 \htext (165 -5){\zoom}

\move (240 0)
\bsegment 
\scriptsize
  \move (0 0) \lcir r:32
  \move (-5 10)  \onedot
  \move (-15 -9) \lcir r:6
  \htext (-10 10) {$s$}   
  \move (0 -40) \lvec (0 -20) \onedot
  \lvec (-40 10)
  \move (0 -20) \lvec (-10 40) \htext (-11 45) {$b$}
  \move (0 -20) \lvec (12 10) \onedot \lcir r:10 \htext (17 10) {$q$}
  \move (0 -20) \lvec (14 -15) \onedot \htext (19 -13) {$r$}
  \htext (-5 -25) {$p$} \htext (-45 12) {$a$}

  \htext (0 -60) {\normalsize $Y_4$}
  \move (0 -70)
\esegment
\end{texdraw}\end{center}
As discussed, it would be enough to indicate $Y_3 \zoom Y_4$,
and if we furthermore
  take advantage of the linear order in $Y_3$ and 
  make the convention that $Y_4$ should be a {\em planar} tree, where the 
  clockwise planar order expresses the (downwards) linear order in $Y_3$,
  then also $Y_3$ is redundant, and we can represent the $4$-opetope
  by the single constellation:
  \begin{center}\begin{texdraw}
  \bsegment \linewd 0.3 \scriptsize
    \move (0 0) \lcir r:32
    \move (6 15) \lcir r:8 
    \move (-15 -9) \onedot  \htext (-18 -13) {$s$}
    \htext (-7 -23) {$p$}   
    \htext (-44 15) {$b$}   
  \linewd 0.5
    \move (0 -40) \lvec (0 -20) \onedot
    \lvec (-40 10)
    \move (0 -20) \lvec (10 40) \htext (11 45) {$a$}
    \move (0 -20) \lvec (17 -8) \onedot  \move (18 -8) \linewd 0.3 \lcir r:8 \htext (22 -8) {$q$}
    \linewd 0.5 \move (0 -20) \lvec (-14 15) \onedot \htext (-19 13) {$r$}

    \htext (0 -60) {\normalsize $Y_4$}
  \esegment
  \end{texdraw}\end{center}
  (While such economy can sometimes be practical,
  conceptually it is rather an 
  obfuscation.)
\end{eks}

\begin{eks}\label{opEx}
  We finish with an example of a $5$-opetope, just to point out
  that there is no longer any natural planar structure on the
  underlying trees in degree $d\geq 5$.
  Arguing as above, to specify
  a $5$-opetope it is enough to specify a single zoom $Z_4 
  \zoom Z_5$, provided we 
  understand that the tree in $Z_4$ is planar (and 
  hence allows us to reconstruct the previous constellation).  Here is 
  an example of a $5$-opetope represented in this economical manner:
  \begin{center}
    \begin{texdraw}
  \bsegment \linewd 0.3 \scriptsize
    \move (0 40) \lcir r:35 \htext (25 21) {$p$}
    \move (-2 15) \lcir r:6 \htext (-4 15) {$x$}
    \move (0 43) \freeEllipsis{16}{28}{90} \htext (0 53) {$y$}
  \linewd 0.5
    \move (0 0) \lvec (0 37) \onedot
    \lvec (-14 46) \onedot \lvec (-34 70)
    \move (0 37) \lvec (14 46) \onedot
    \lvec (7 82)
    \move (14 46)
    \lvec (21 80)
    \htext (-5 34) {$b$}
    \htext (17 41) {$c$}
    \htext (-16 40) {$a$}
    
    \htext (0 -20) {\normalsize $Z_4$}
  \esegment
  
 \htext (85 30){\zoom}

  \move (165 0)
  \bsegment \linewd 0.3 \scriptsize
    \move (0 40) \lcir r:35
    \move (4 58) \lcir r:6
  \linewd 0.5
    \move (0 0) \lvec (0 20) \onedot
    \lvec (-15 30) \onedot
    \move (0 20) \lvec (15 42) \onedot
    \lvec (-10 80) \htext (-11 85) {$a$}
    \move (15 42) \lvec (18 77) \htext (18 83) {$b$}
    \move (15 42) \lvec ( 35 63) \htext (38 66) {$c$}
    \htext (6 18) {$p$}
    \htext (19 37) {$y$}
    \htext (-20 30) {$x$}

    \htext (0 -20) {\normalsize $Z_5$}
  \esegment
  \end{texdraw}\end{center}
\end{eks}

\subsection{Formal definitions: trees and constellations}

While the presented definition of opetopes is appealing in its simplicity,
scrutiny of the definition raises some questions: what exactly is meant by tree?
Is it a combinatorial notion?  In that case, what does it mean to draw circles
on a tree?  And when we say `tree', `constellation', or `opetope', do we refer
to concrete specific sets with structure or do we refer to isomorphism classes
of such?  In this subsection we give the definitions a more formal treatment.
We show in particular that the notion of constellation is purely combinatorial
and does not depend on geometric realisation.  Secondly, the analysis will
clarify the relation to Baez-Dolan metatrees (and uncover the shortcoming with
these).  Thirdly, the insight provided by the formal viewpoint will be helpful
for understanding the constructions in Section~\ref{Sec:BD} and the calculations
in Section~\ref{Sec:calculus}.

The question of explicit-sets-with-structure versus their isomorphism classes
deserves a remark before the definitions.  We want to define the various notions
(trees, constellations, zoom complexes, opetopes) in terms of finite sets with
some structure, in order to classify as combinatorial notions.  As such these
objects form a proper class.  On the other hand, naturally we are mostly
interested in these structures up to isomorphism.  Our choice will be to stick
with the explicit finite-sets-with-structure as long as the objects may possess
non-trivial automorphisms (which is the case for trees, constellations, and zoom
complexes), but consider isomorphism classes for rigid objects like $P$-trees
(trees decorated by a polynomial endofunctor $P$, as introduced in
\ref{P-trees}) and opetopes.  Hence an opetope will be defined as a set of
isomorphism classes of certain (rigid) objects (this was implicit in \ref{def},
and in particular there will be only a small set of them.  This is in accordance
with previous definitions of opetopes in the literature --- in fact this issue
had not previously come up since there was no combinatorial description 
available.

\begin{blanko}{Graphs.}
  By a {\em graph} we understand a pair $(T_0,T_1)$, where $T_0$ is a set,
  and $T_1$ is a set of subsets of $T_0$ of cardinality $2$.  The elements
  in $T_0$ are called {\em vertices}, and the elements in $T_1$ {\em
  edges}.  An edge $\{x,y\}$ is said to be {\em incident} to a vertex
  $v$ if $v\in \{x,y\}$.  We say a vertex is of valence $n$ if the set
  of incident edges is of cardinality $n$.
  
  The geometric realisation of a graph is the CW-complex with a $0$-cell
  for each vertex, and for each edge a $1$-cell attached at the points
  corresponding to its two incident vertices.
\end{blanko}

\begin{blanko}{Trees.}\label{tree-formal}
  By a {\em finite rooted tree with boundary} we mean a finite graph
  $T=(T_0,T_1)$, connected and simply connected, equipped with a pointed
  subset $\d T$ of vertices of valence $1$, called the {\em boundary}.
  We will not need other kinds
  of trees than finite rooted trees with boundary, and we will simply call 
  them {\em trees}.  (An alternative tree formalism is developed in 
  \cite{Kock:0807}.)
  
  The basepoint $t_0 \in \d T$ is called the {\em output vertex}, and
  the remaining vertices in $\d T$ are called {\em input vertices}.  Most
  of the time we shall not refer to the boundary vertices at all, and
  graphically a boundary vertex is just represented as a loose end of the
  incident edge.  Edges incident to input vertices are called {\em leaves}
  or {\em input edges} of the tree, while the unique edge incident to the
  output vertex is called the {\em root edge} or the {\em output edge}
  of the tree.
  The vertices in $T_0 \shortsetminus \d T$ are called {\em
  nodes} or {\em dots}; we draw them as dots.  A tree may have zero
  dots, in which case it is just a single edge (together with two boundary
  vertices, which we suppress);  we call such a tree a {\em unit tree}.
  Not every vertex of valence $1$ needs to be a
  boundary vertex: those which are not are called {\em null-dots}.
  
    \begin{center}\begin{texdraw}
  \linewd 0.5 \footnotesize
  \move (-50 5)
  \bsegment
    \move (0 0) \lvec (0 25)
    \move (0 -10) \htext {\footnotesize unit tree}
  \esegment
  
  \move (56 0)
  \bsegment
    \move (0 0) \lvec (0 15) \onedot
    \move (0 15) \lvec (-8 28) \onedot
    \move (0 15) \lvec (12 40)
        \move (-36 29) \htext {\footnotesize null-dot $\to$}
    \move (27 33) \htext {\footnotesize $\leftarrow$ leaf}
    \move (-30 7) \htext {\footnotesize root edge $\to$}

  \esegment
  
\end{texdraw}\end{center}

    The standard graphical representation of trees is justified by geometric
  realisation.  Note that leaves and root are realised by half-open
  intervals, and we keep track of which are which by always drawing the
  root at the bottom. 
  
  An {\em isomorphism of trees} is an isomorphism of the underlying graphs
  preserving root and leaves.  A tree can be recovered up to isomorphism by its
  geometric realisation.  We shall frequently be interested only in the
  isomorphism classes.  This was implicit in the `5-minute' definition.
  
  If $T=(T_0,T_1,\d T, t_0)$ is a tree, the set $T_0$ has a natural poset
  structure $a\leq b$, in which the input vertices and null-dots are
  minimal elements and the output vertex is the maximal element. 
  We say $a$ is a {\em
  child} of $b$ if $a\leq b$ and $\{a,b\}$ is an edge.  Each dot has one
  output edge, and the remaining incident edges are called input edges of
  the dot.

\end{blanko}

\begin{blanko}{Nestings and correspondences.}
  Nestings (as in \ref{nest}) are just another graphical representation of
  an abstract tree $(T_0, T_1, \d T, t_0)$.  Graphically, a {\em nesting}
  is a collection of non-intersecting spheres and dots, which either
  consists of a single dot (and no spheres) or has one outer sphere,
  containing all the other spheres and dots.  We identify two nestings if
  there is an isotopy between them.  We shall need some more terminology
  about nestings, expanding the dictionary between trees and nestings.
%
  A sphere that does not contain any other spheres or dots is called
  a {\em null-sphere}.  These correspond exactly to the null-dots of a
  tree.  The region bounded on the outside by a sphere $S$ and on the
  inside by the dots and spheres contained in $S$ is called a {\em layer}.
  The layers of a nesting correspond to the nodes of the tree.
  An inner sphere mediates between
  two layers just like an inner edge in a tree sits between two nodes.  We
  will often confuse a layer with its outside bounding sphere. 

\end{blanko}

\begin{blanko}{Towards a combinatorial definition of constellations.}%
\label{towards-const-formal}
  In \ref{const} we defined a constellation as a tree with a sphere
  nesting on top, more precisely as a
  configuration $C$ of edges, dots, and
  spheres (in $3$-space), such that: (i) edges and dots form a tree,
  (ii) dots and
  spheres form a nesting, and (iii) for each sphere, the edges and dots contained
  in it form a tree again.  
   This definition has a clear intuitive content, and plays an important
   role as convenient tool for manipulating constellations and opetopes,
   just like we usually manipulate trees in terms of their geometrical
   aspect, not in terms of abstract graphs.  However, the definition
   depends on geometric realisation, and it is not clear at this point of
   our exposition that it is a rigorous notion at all.  It is likely
   that the definition can be formalised geometrically by talking about
   isotopy classes of such configurations of (progressive) line segments,
   dots, and spheres in Euclidean space.  We shall not go further into
   this.  We wish instead to stress that the notion can be given in purely
   combinatorial terms.
%
%
%
  The idea is to capture the structure by specifying some bijections between
  the underlying tree and the tree corresponding to the nesting.
  For this to work it is necessary to mark 
  the position of the null-spheres by temporarily turning them into dots.
  This is formalised through the notion of subdivision of trees:
\end{blanko}

\begin{blanko}{Subdivision and kernels.}
  A {\em linear tree} is a tree in which every dot has exactly one input edge.
  The unit tree is an example of a linear tree.  (Up to isomorphism) 
  there is one linear tree for each natural number.
  A {\em subdivision} of a tree $T$ is a tree $T'$ obtained by 
  replacing each edge by a linear tree.  
  We draw the new dots as white dots.  Here is a picture of a tree and
  a subdivision:
   \begin{equation}\label{subdiv-fig}
    \begin{texdraw}
    \bsegment 
    \linewd 0.5 
    \move (0 -3) \lvec (0 20) \onedot
    \lvec (-15 40) \onedot 
    \lvec (-30 54)
    \move (-15 40) \lvec (-12 64)
    \move (0 20) \lvec (5 70)
    \move (0 20) \lvec (18 32) \onedot
	    \move (0 -20) \htext {$T$}
    \esegment
    \move (100 0)
    \bsegment 
    \move (0 -3) \lvec (0 20) \onedot
    \lvec (-15 40) \onedot 
    \lvec (-30 54)
    \move (-15 40) \lvec (-12 64)
    \move (0 20) \lvec (5 70)
    \move (0 20) \lvec (18 32) \onedot
	  \linewd 0.3
	    \move (0 9) \lcir r:2
	    \move (-7.5 30) \lcir r:2
	    \move (2 40) \lcir r:2
	    \move (3.6 56) \lcir r:2

	    \move (0 -20) \htext {$T'$}
    \esegment
    \end{texdraw}
  \end{equation}

  When we speak about dots of a subdivided tree we mean the union of
  old and new dots:
  $$
  \operatorname{dots}(T') = \operatorname{blackdots}(T') +
  \operatorname{whitedots}(T')
  $$
  (note that $\operatorname{blackdots}(T')=\operatorname{dots}(T)$).
  \medskip

%
  
   If $T$ is a tree, every subset $K \subset \operatorname{dots}(T)$ spans a
   full subgraph $K^\dagger$, where an edge of $K^\dagger$ is an edge of $T$
   connecting two nodes of $K$.  We call $K$ a {\em kernel} if the graph
   $K^\dagger$ is non-empty and connected.  A kernel $K$ spans a tree with
   boundary $K^\ddagger$, whose edges are those of $T$ adjacent to an element of
   $K$; the dots of $K^\ddagger$ are the elements of $K$ and the boundary
   vertices of $K^\ddagger$ are those vertices of $K^\dagger$ not in $K$.  In
   other words, a sphere containing exactly the dots of a kernel cuts a tree, as
   in condition (iii) of \ref{const}.  In the following picture, $K=\{r,u,v\}$
   is an example of a kernel, $K^\dagger$ is indicated with fat edges, and the
   tree $K^\ddagger$ is what's inside the sphere:
   \begin{center}
    \begin{texdraw}
        \bsegment 
    \linewd 0.5 
    \move (0 0) \lvec (0 20) \onedot
    \lvec (-27 65) 
    \move (-18 50) \onedot 
    \move (0 20) \lvec (0 40) \onedot
    \lvec (-6 56) \onedot
    \move (0 40) \lvec (6 56) \onedot
    \move (0 20) \lvec (24 65) 
    \move (0 20) \lvec (15 30) \onedot

    \htext (6 18) {\footnotesize $r$}
    \htext (5 40) {\footnotesize $u$}
    \htext (20 32) {\footnotesize $v$}

    \move (5 28) \freeEllipsis{23}{21}{45}
    
    \linewd 1.8
    \move (0 40) \lvec (0 20) \lvec ( 15 30)
    \esegment

    \end{texdraw}
  \end{center}
%
  (When we speak of kernels of a subdivided tree we refer to
  all dots, black and white.) 

%
\end{blanko}

%

\begin{blanko}{Combinatorial definition of constellation.}\label{const-formal}
  A {\em constellation} $C: T \to N$ between two trees $T$ and $N$
  is a triple $(T',\sigma_{\bullet}, \sigma_{\circ})$,
  where $T'$
  is a subdivision 
  of $T$,  and $\sigma_\bullet$ and $\sigma_\circ$ are bijections
  \begin{eqnarray*}
    \sigma_\bullet : \operatorname{blackdots}(T') & \isopil & 
    \operatorname{leaves}(N)  \\
    \sigma_\circ : \operatorname{whitedots}(T') & \isopil & 
    \operatorname{nulldots}(N)
  \end{eqnarray*}
  such that the sum map
  $ \sigma \df \sigma_\bullet +\sigma_\circ $
  satisfies the {\em kernel rule}:
\begin{equation}\label{kernel-rule}
  \text{
   for each $x \in \operatorname{dots}(N)$, the set 
   $
   \{t \in \operatorname{dots}(T') \mid  \sigma(t) \leq x \}
   $
   is a kernel in $T'$.
   }
   \end{equation}
   An {\em  isomorphism of constellations} consists of isomorphisms
   of the underlying (subdivided) trees compatible with the structural 
   bijections.
\end{blanko}

Here is a picture of a constellation in this sense:
   \begin{equation}\label{const-formal-ex}
    \begin{texdraw}
    \bsegment 
    \linewd 0.5 
    \move (0 -5) \lvec (0 20) \onedot
    \lvec (-15 45) \onedot 
    \lvec (-35 70)
    \move (-15 45) \lvec (-10 80)
        \move (-18 40) \htext{\footnotesize $a$}
    \move (-4 16) \htext{\footnotesize $b$}

    \move (0 20) \lvec (30 80)
	  \linewd 0.3
	    \move (6.5 33) \lcir r:2
	    \move (11 43) \lcir r:2
	    \move (19 58) \lcir r:2
    \move (13 30) \htext {\footnotesize $x$}
    \move (18 40) \htext {\footnotesize $y$}
    \move (26 56) \htext {\footnotesize $z$}

	    \move (0 -20) \htext {$T$}
    \esegment
    \htext (50 -17) {$\stackrel{C}{\longrightarrow}$}
    \move (100 0)
        \bsegment 
    \linewd 0.5 
    \move (0 0) \lvec (0 20) \onedot
    \lvec (-27 65) 
    \move (-9 35) \onedot \htext (-15 32) {\footnotesize $p$}
    \move (-18 50) \onedot \htext (-24 47) {\footnotesize $q$} 
    \move (-28 70) \htext {\footnotesize $a$}
    \move (0 20) \lvec (0 40) \onedot
    \lvec (-6 56) \onedot
    \move (0 40) \lvec (6 56) \onedot
    \move (0 20) \lvec (24 65) \move (25 71) \htext {\footnotesize $b$}
    \move (0 20) \lvec (15 30) \onedot
    
    \move (-7 63.5) \htext {\footnotesize $x$}
    \move (7 63.5) \htext {\footnotesize $y$}
    \move (22 33) \htext {\footnotesize $z$}
    	    \move (0 -20) \htext {$N$}

    \esegment

    \end{texdraw}
  \end{equation}
  (The white dots are not a part of $T$; they represent the subdivision
  of $T$ which is a part of the data constituting $C$.)
  
  \medskip
  
Let us compare the definition of constellation given in
\ref{const-formal} with the drawings of \ref{const}, justifying that the
latter constitute a faithful graphical representation of the former.
Given a constellation according to definition~\ref{const-formal}, as in 
Figure~(\ref{const-formal-ex}),  for each dot $x$ in $N$ that is 
not a null-dot, draw a sphere
  in $T'$ around those dots in $T'$ corresponding to the descendant leaves
  and null-dots of $x$ in $N$, as in the kernel rule~(\ref{kernel-rule}).
  The kernel rule tells us 
  that this sphere cuts a tree (as in \ref{const}).  The sphere must
  be drawn inside the sphere corresponding to the parent node of $x$ (if any);
  this ensures that the spheres are non-intersecting and that the resulting 
  nesting corresponds to the tree $N$.
  (Name the spheres and white dots in $T'$ by the corresponding dots in $N$.)
  To finish the construction, replace the white dots in $T'$ by null-spheres. 
  \begin{equation}\label{whitedots}
    \begin{texdraw}
    \bsegment 
    \linewd 0.5 
    \move (0 -5) \lvec (0 20) \onedot
    \lvec (-15 45) \onedot 
    \lvec (-35 70)
    \move (-15 45) \lvec (-10 80)
    \move (-16 44) \linewd 0.3  \lcir r:9 \lcir r:13
        \move (-18 40) \htext{\footnotesize $a$}
    \move (-4 16) \htext{\footnotesize $b$}

    \linewd 0.5 
    \move (0 20) \lvec (30 80)
	  \linewd 0.3
	    \move (11 36) \freeEllipsis{14}{12}{65}
	    \move (6.5 33) \lcir r:4
	    \move (11 43) \lcir r:4
	    \move (19 58) \lcir r:4
	    \move (0 40) \lcir r:35
    \move (13 30) \htext {\footnotesize $x$}
    \move (18 40) \htext {\footnotesize $y$}
    \move (26 56) \htext {\footnotesize $z$}

	    \move (-50 35) \htext {$C$}
	    \move (0 -20) \htext {$T$}
    \esegment

    \move (100 0)
        \bsegment 
    \linewd 0.5 
    \move (0 0) \lvec (0 20) \onedot
        \lvec (-27 65) 
    \move (-9 35) \onedot \htext (-15 32) {\footnotesize $p$}
    \move (-18 50) \onedot \htext (-24 47) {\footnotesize $q$}

    \htext (-28 70) {\footnotesize $a$}
    \move (0 20) \lvec (0 40) \onedot
    \lvec (-6 56) \onedot
    \move (0 40) \lvec (6 56) \onedot
    \move (0 20) \lvec (24 65) \move (25 71) \htext {\footnotesize $b$}
    \move (0 20) \lvec (15 30) \onedot
    	    \move (0 -20) \htext {$N$}
    \move (-7 63.5) \htext {\footnotesize $x$}
    \move (7 63.5) \htext {\footnotesize $y$}
    \move (22 33) \htext {\footnotesize $z$}

    \esegment

    \end{texdraw}
  \end{equation}
   It is now
  clear that the left-hand side of the picture is a constellation in the
  sense of \ref{const}.  (Note that if $N$ has no dots, then in particular it 
  has no null-dots, so $T=T'$. Furthermore in this case $N$ must have precisely 
  one leaf, so $T$ has just one dot and is therefore a 
  constellation even without any spheres drawn.)
  
  Conversely, given an constellation $C$ in the sense of \ref{const}, with
  underlying tree $T$,
    \begin{center}
    \begin{texdraw}
    \bsegment 
    \move (-50 40) \htext{$C$}
    \linewd 0.5 
    \move (0 -5) \lvec (0 20) \onedot
    \lvec (-15 45) \onedot 
    \lvec (-35 70)
    \move (-15 45) \lvec (-10 80)

    \move (0 20) \lvec (30 80)
	  \linewd 0.3
    \move (-15 45) \lcir r:9 \lcir r:13
	    \move (9 38) \freeEllipsis{14}{9}{65}
	    \move (6.5 33) \lcir r:4
	    \move (11 43) \lcir r:4
	    \move (19 58) \lcir r:4
	    \move (0 40) \lcir r:35
	    \move (0 -20) \htext {$T$}
    \esegment
    \end{texdraw}
  \end{center}
  the preceding arguments can be reversed to
  construct a constellation according to the combinatorial
  definition~\ref{const-formal}: first draw the tree $N$ corresponding to
  the underlying nesting of $C$ (using the spheres as names for the dots in $N$)
  (this gives Figure~(\ref{whitedots})), then erase
  all the spheres in $C$ except the null-spheres, and draw the null-spheres
  so small that they look like (white) dots --- they constitute now a
  subdivision of $T$.  At this point we have a constellation in the
  sense of \ref{const-formal}: the bijections $\sigma_\bullet$ and
  $\sigma_\circ$ are already part of the correspondence between the
  underlying nesting of $C$ and the tree $N$, and each dot $x\in
  \operatorname{dots}(N)$ corresponds to a sphere in $C$, so the
  kernel rule~(\ref{kernel-rule})
  is just a reformulation of the condition that each
  sphere cuts a tree.
  
  It is clear that a constellation in the sense of \ref{const-formal} can be 
  recovered uniquely from its \ref{const}-interpretation.

\begin{blanko}{Zooms and zoom complexes, revisited.}\label{zoom-complex-rev}
  Now that the notion of constellation has been formalised, the definitions of
  zoom (\ref{zoom}) and zoom complex (\ref{zoom-complex}) are already formal.
  Let us unravel these notions by plugging in the combinatorial definition 
  of constellation (\ref{const-formal}).
  Given a zoom 
     \begin{center}
    \begin{texdraw}
    \bsegment 
    \linewd 0.5 
    \move (0 -5) \lvec (0 20) \onedot
    \lvec (-15 45) \onedot 
    \lvec (-35 70)
    \move (-15 45) \lvec (-10 80)
        \move (-18 40) \htext{\tiny $a$}
    \move (-4 16) \htext{\tiny $b$}

    \move (0 20) \lvec (30 80)
	  \linewd 0.3
    \move (-16 44) \lcir r:9
	    \move (11 36) \freeEllipsis{14}{12}{65}
	    \move (6.5 33) \lcir r:4
	    \move (11 43) \lcir r:4
	    \move (19 58) \lcir r:4
	    \move (0 40) \lcir r:35
    \move (13 30) \htext {\tiny $x$}
    \move (18 40) \htext {\tiny $y$}
    \move (26 56) \htext {\tiny $z$}

	    \move (0 -20) \htext {$C_1$}
    \esegment

    \move (70 30) \htext {$\zoom$}
    \move (160 0)
        \bsegment 
    \linewd 0.5 
    \move (0 -5) \lvec (0 20) \onedot
    \lvec (-13 33) \onedot 
    \lvec (-50 70) 
\linewd 0.3    	    \move (-32 52) \lcir r:4
    \move (-55 75) \htext {\tiny $a$}
    \linewd 0.5 
    \move (0 20) \lvec (0 40) \onedot
    \lvec (-6 56) \onedot
    \move (0 40) \lvec (6 56) \onedot
    \move (0 20) \lvec (35 75) \move (37 80) \htext {\tiny $b$}
    \move (0 20) \lvec (15 30) \onedot
    
    \move (-7 63.5) \htext {\tiny $x$}
    \move (7 63.5) \htext {\tiny $y$}
    \move (22 33) \htext {\tiny $z$}
    \linewd 0.3
        	    \move (21 52) \lcir r:4

    	    \move (0 40) \freeEllipsis{50}{40}{0}
    	    \move (0 50) \freeEllipsis{12}{14}{0}
    	    \move (6 41) \freeEllipsis{27}{30}{0}

    	    \move (0 -20) \htext {$C_2$}

    \esegment

    \end{texdraw}
  \end{center}
  the formal definition of constellation (\ref{const-formal}) leads to this 
  drawing:
\begin{center}
    \begin{texdraw}
      \setunitscale 0.9
      \tiny
      \bsegment
      \move (-50 -40)
      \lvec (160 -40) \lvec (160 100) \lvec (-50 100) \lvec (-50 -40)
      
      \move (0 0)
    \bsegment 
    \linewd 0.5 
    \move (0 -5) \lvec (0 20) \onedot
    \lvec (-15 45) \onedot 
    \lvec (-35 70)
    \move (-15 45) \lvec (-10 80)
        \move (-18 40) \htext{$a$}
    \move (-4 16) \htext{$b$}

    \move (0 20) \lvec (30 80)
	  \linewd 0.3
	    \move (6.5 33) \lcir r:2
	    \move (11 43) \lcir r:2
	    \move (19 58) \lcir r:2
    \move (13 30) \htext {$x$}
    \move (18 40) \htext {$y$}
    \move (26 56) \htext {$z$}

    \esegment
    \htext (50 17) {\normalsize $\stackrel{C_1}{\longrightarrow}$}

    \move (100 0)
        \bsegment 
    \linewd 0.5 
    \move (0 -5) \lvec (0 20) \onedot
    \lvec (-15 45) \onedot 
    \lvec (-27 65) \move (-28 70) \htext {$a$}
    \move (0 20) \lvec (0 40) \onedot
    \lvec (-6 56) \onedot
    \move (0 40) \lvec (6 56) \onedot
    \move (0 20) \lvec (24 65) \move (25 71) \htext {$b$}
    \move (0 20) \lvec (15 30) \onedot
    
    \move (-7 63.5) \htext {$x$}
    \move (7 63.5) \htext {$y$}
    \move (22 33) \htext {$z$}

    \esegment
\esegment

\move (195 30) \htext {$\zoom$}

\move (280 0)

      \bsegment
      \move (-50 -40)
      \lvec (160 -40) \lvec (160 100) \lvec (-50 100) \lvec (-50 -40)
      
      \move (0 0)
        \bsegment 
    \linewd 0.5 
    \move (0 -5) \lvec (0 20) \onedot
    \lvec (-15 45) \onedot 
    \lvec (-27 65) 
    \move (-28 70) \htext {$a$}
    \move (-20 42) \htext {$w$}

    \move (0 20) \lvec (0 40) \onedot
    \lvec (-6 56) \onedot
        \move (-5 38) \htext {$v$}
        \move (-5 17) \htext {$u$}

    \move (0 40) \lvec (6 56) \onedot
    \move (0 20) \lvec (24 65) \move (25 71) \htext {$b$}
    \move (0 20) \lvec (15 30) \onedot
    
    \move (-7 63.5) \htext {$x$}
    \move (7 63.5) \htext {$y$}
    \move (22 33) \htext {$z$}
    \linewd 0.3
    	    \move (-23 58) \lcir r:2
	    \move (-27 53) \htext {$p$}
	    \move (16 49.5) \lcir r:2
	    \move (20 45) \htext {$q$}

    \esegment
    \htext (50 17) {\normalsize $\stackrel{C_2}{\longrightarrow}$}

    \move (100 0)
        \bsegment 
    \linewd 0.5 
    \move (0 -5) \lvec (0 20) \onedot
    \move (0 20) \lvec (15 30) \onedot
        	    \move (21 32) \htext {$p$}

    \move (0 20) \lvec (0 35) \onedot
    \lvec (-18 40) \onedot
    	    \move (-24 42) \htext {$q$}

    \move (0 35) \lvec (-7 50) \onedot
    \lvec (-26 70) 
    \move (-7 50) \lvec (-17 78)
    \move (-7 50) \lvec (-5 80)
      \move (-28 75) \htext {$x$}
      \move (-18 82) \htext {$y$}
      \move (-5 85) \htext {$v$}

    \move (0 35) \lvec (7 80)
    \move (0 35) \lvec (20 75)
    \move (0 35) \lvec (27 65)
      \move (9 85) \htext {$u$}
      \move (22 80) \htext {$w$}
      \move (30 69) \htext {$z$}

    \esegment
\esegment

    \end{texdraw}
  \end{center}
The defining property of zoom means the two trees in
the middle coincide (modulo the subdivision, which is rather a
part of the structure of $C_2$), so we can overlay the two constellations:
   \begin{center}
    \begin{texdraw}
      \setunitscale 0.9
      \tiny
      \bsegment
      \move (-50 -20)
      \lvec (155 -20) \lvec (155 100) \lvec (-50 100) \lvec (-50 -20)
      
      \move (0 0)
    \bsegment 
    \linewd 0.5 
    \move (0 -5) \lvec (0 20) \onedot
    \lvec (-15 45) \onedot 
    \lvec (-35 70)
    \move (-15 45) \lvec (-10 80)
        \move (-18 40) \htext{$a$}
    \move (-4 16) \htext{$b$}

    \move (0 20) \lvec (30 80)
	  \linewd 0.3
	    \move (6.5 33) \lcir r:2
	    \move (11 43) \lcir r:2
	    \move (19 58) \lcir r:2
    \move (13 30) \htext {$x$}
    \move (18 40) \htext {$y$}
    \move (26 56) \htext {$z$}

    \esegment

\esegment

    \move (100 5)


      \bsegment
      \move (-50 -20)
      \lvec (155 -20) \lvec (155 100) \lvec (-50 100) \lvec (-50 -20)
      
      \move (0 0)
        \bsegment 
    \linewd 0.5 
    \move (0 -5) \lvec (0 20) \onedot
    \lvec (-15 45) \onedot 
    \lvec (-27 65) 
    \move (-28 70) \htext {$a$}
    \move (-20 42) \htext {$w$}

    \move (0 20) \lvec (0 40) \onedot
    \lvec (-6 56) \onedot
        \move (-5 38) \htext {$v$}
        \move (-5 17) \htext {$u$}

    \move (0 40) \lvec (6 56) \onedot
    \move (0 20) \lvec (24 65) \move (25 71) \htext {$b$}
    \move (0 20) \lvec (15 30) \onedot
    
    \move (-7 63.5) \htext {$x$}
    \move (7 63.5) \htext {$y$}
    \move (22 33) \htext {$z$}
    
    	    \move (-23 58) \lcir r:2
	    \move (-27 53) \htext {$p$}
	    \move (16 49.5) \lcir r:2
	    \move (20 45) \htext {$q$}

    \esegment

    \move (105 -8)
        \bsegment 
    \linewd 0.5 
    \move (0 0) \lvec (0 20) \onedot
    \move (0 20) \lvec (15 30) \onedot
        	    \move (21 32) \htext {$p$}

    \move (0 20) \lvec (0 35) \onedot
    \lvec (-18 40) \onedot
    	    \move (-24 42) \htext {$q$}

    \move (0 35) \lvec (-7 50) \onedot
    \lvec (-26 70) 
    \move (-7 50) \lvec (-17 78)
    \move (-7 50) \lvec (-5 80)
      \move (-28 75) \htext {$x$}
      \move (-18 82) \htext {$y$}
      \move (-5 85) \htext {$v$}

    \move (0 35) \lvec (7 80)
    \move (0 35) \lvec (20 75)
    \move (0 35) \lvec (27 65)
      \move (9 85) \htext {$u$}
      \move (22 80) \htext {$w$}
      \move (30 69) \htext {$z$}

    \esegment
\esegment

    \end{texdraw}
  \end{center}
  In conclusion, a zoom is a sequence of three trees connected by 
  constellations: 
  $$
  T_0 \stackrel{C_1}{\longrightarrow} 
  T_1 \stackrel{C_2}{\longrightarrow} 
  T_2 .
  $$
  
  Similarly,  a zoom complex
is a sequence of trees and constellations
  \begin{equation}\label{zoomT}
  T_0 \stackrel{C_1}{\longrightarrow} T_1 \stackrel{C_2}{\longrightarrow} 
  T_2 \ \ \ \cdots \ \ \ T_{n-1} \stackrel{C_n}{\longrightarrow} T_n .
  \end{equation}
  
  An {\em isomorphism of zoom complexes} is a sequence of isomorphisms of
  constellations, compatible with the zoom bijections.  In the viewpoint of
  (\ref{zoomT}) it is a sequence of isomorphisms of subdivided trees compatible 
  with the structural bijections of \ref{const-formal}.  Note that a zoom 
  complex of any degree may allow non-trivial automorphisms.  For example,
  the following zoom complex has a non-trivial involution:
\begin{center}\begin{texdraw}
  \linewd 0.5 \tiny
  
  \move (0 0)
  \bsegment
    \move (0 0) \lvec (0 15) \onedot 
    \lvec (0 45)
    \move (0 15)\lvec (-15 40) 
    \move (0 15)\lvec (15 40) 
    
    \move (-7 26) \lcir r:4
    \move (7 26) \lcir r:4
    \htext (-10 20) {$u$}
    \htext (10 20) {$v$}
    
 \move (0 22) \lcir r:17
 \esegment
 \htext (40 10){\zoom}
   \move (80 0)
  \bsegment
    \move (0 0) \lvec (0 13) \onedot 
    \lvec (0 45) 
    \move (0 13)\lvec (-8 23) \onedot
    \move (0 13)\lvec (8 23) \onedot
    \htext (-10 28) {$u$}
    \htext (10 28) {$v$}
 \move (0 22) \lcir r:17
 \esegment
 \htext (120 10){\zoom}
   \move (160 0)
  \bsegment
    \move (0 0) \lvec (0 15) \onedot 
    \lvec (0 45)
    \move (0 15)\lvec (-15 40) 
    \move (0 15)\lvec (15 40) 
    \htext (-18 43) {$u$}
    \htext (18 43) {$v$}
    
    \move (-7 26) \lcir r:4
    \move (7 26) \lcir r:4
 \move (0 22) \lcir r:17
 \esegment
\end{texdraw}\end{center}
\end{blanko}

\begin{blanko}{Opetopes, revisited.}\label{opetope-metatree}
  We defined the $k$-opetopes to be the isomorphism classes of zoom complexes
  of degree $k$ subject to an initial condition (\ref{def}).  Observe that
  such zoom complexes are rigid objects (i.e.~have no non-trivial automorphisms).  
  Indeed,
  any non-trivial automorphism of a zoom complex $C$ induces a non-trivial
  automorphism already on the underlying tree of $C_0$ because of the structural
  bijections in the definition of zoom.  Clearly the initial condition
  precludes non-trivial automorphisms in $C_0$.

  We saw in \ref{ex-1-2-3} that an opetope of dimension $2$ can be re\-presented
  by a linear tree, and an opetope of dimension $3$ by a planar tree, which is
  the same thing as a nesting on a linear tree.  In other words, an opetope of
  dimension $3$ can be represented as a constellation $T_2\to T_3$, where $T_2$
  is a linear tree.  In general, an opetope of dimension $n\geq 3$ can be
  represented by a sequence of trees and constellations
\begin{equation}\label{opetopeT}
  T_2 \stackrel{C_3}\longrightarrow T_3 \stackrel{C_4}\longrightarrow
  \cdots\stackrel{C_n}\longrightarrow T_n ,
\end{equation}
  with $T_2$ a linear tree, or equivalently, as
  \begin{equation}\label{opetopeC}
  C_3\zoom C_4 \zoom\cdots 
  \zoom  C_n ,
  \end{equation}
  where $C_3$ is the constellation associated to a planar tree as in
  \ref{ex-1-2-3}.  The sequence (\ref{opetopeC}) is
  graphically redundant compared to the sequence (\ref{opetopeT}),
  but drawing the redundant spheres is very practical as they explicitly witness
  the validity of the kernel rule~(\ref{kernel-rule}).
%
\end{blanko}

\begin{blanko}{Relation with Baez-Dolan metatrees.}\label{const v meta}
  The viewpoint on zoom complexes given in \ref{zoom-complex-rev} provides
  an explicit comparison with the notion of metatree
  introduced by Baez and Dolan~\cite{Baez-Dolan:9702}.  There are two 
  important differences.
  
  A {\em metatree} (cf.~\cite{Baez-Dolan:9702}, pp.~176--177) is
  essentially a sequence of trees $T_0,\ldots,T_n$ not allowed to have
  null-dots, with specified bijections $\sigma_\bullet{}_i :
  \operatorname{dots}(T_{i-1}) \isopil \operatorname{leaves}(T_i)$
  satisfying the kernel rule~(\ref{kernel-rule}).
  In other words, it is the special case of a
  zoom complex where the trees have no null-dots, and hence there is no
  subdivision involved in the constellations.  Null-dots represent nullary
  operations of the operads or polynomial monads of the Baez-Dolan
  construction~\ref{BD}, and nullary operations do arise.  Therefore the
  Baez-Dolan metatrees seem to be insufficient to reflect the Baez-Dolan
  construction and to describe opetopes.  Our zoom complexes may be
  what Baez and Dolan really envisaged with the notion of metatree.
  
%
%
%

  The second difference is of another nature: Baez and Dolan worked with
  planar trees, but introduced a notion of {\em combed tree}, in which the
  leaves are allowed to cross each other in any permutation.  The trees in
  Baez-Dolan metatrees are in fact combed.  These artefacts come from
  working with symmetric operads.  The effect on the definition of opetope
  is that each opetope comes equipped with an ordering of its faces.  We
  work instead with non-planar trees and polynomial monads, and the
  resulting opetopes (which agree with Leinster's, cf.~\ref{BD-formal}) are
  `un-ordered' like abstract geometric objects.  Planarity is revealed to
  be a special feature of dimension $3$, cf.~\ref{ex-1-2-3}.
  
  Let us remark that we think 
  the spheres are an important conceptual device for understanding
  opetopes in terms of sequences of trees.  Baez and Dolan stressed
  that a key feature of the slice construction is that operations are
  promoted to types, and reduction laws are promoted to operations.
  This two-level correspondence comes to the fore with the notion of
  zoom: the types are represented by the leaves, the operations are
  the dots, and the reduction laws are expressed by the spheres.  The
  zoom relation shifts dots to leaves and spheres to dots.
%
\end{blanko}

\section{Polynomial functors and polynomial monads}
\label{Sec:poly}

\begin{blanko}{Polynomial functors.}
  We recall some facts about polynomial functors. (Details for the notions 
  needed here can be found in \cite{Gambino-Kock}.  The manuscript
  \cite{Kock:NotesOnPolynomialFunctors} aims at eventually becoming a more
  comprehensive reference.)  A diagram of sets and set maps like this
  \begin{equation}\label{P}
    \polyFunct{I}{E}{B}{J}{s}{p}{t}{}
  \end{equation}
    gives rise to a {\em polynomial functor} $P:\Set/I \to 
    \Set/J$ defined by
  $$
  \Set/I \stackrel{s\upperstar }{\rTo} \Set/E 
  \stackrel{p\lowerstar }{\rTo} \Set/B 
  \stackrel{t_!}{\rTo} 
  \Set/J .
  $$
  Here lowerstar and lowershriek denote, respectively, the right 
  adjoint and the left adjoint of the pullback functor upperstar.
  In explicit terms, the functor is given by
    \begin{eqnarray*}
    \Set/I & \longrightarrow & \Set/J  \\
    {}[f:X\to I] & \longmapsto & \sum_{b\in B} \prod_{e\in E_b} X_{s(e)}
  \end{eqnarray*}
  where $E_b \df p^{-1}(b)$ and $X_i \df f^{-1}(i)$, and where the
  last set is considered to be over $J$ via $t_!$.
  
  We will always assume that $p:E \to B$ has finite fibres.  No finiteness
  conditions are imposed on the individual sets $I$, $J$, $E$, $B$, nor
  on the fibres of $s$ and $t$.
\end{blanko}
  
\begin{blanko}{Graphical interpretation.}\label{graphical}
  The following graphical interpretation links 
  polynomial functors to the tree structures of Section~\ref{Sec:opetopes}.
  (This interpretation is not a whim: there is a deeper relationship between polynomial
  functors and trees, analysed more closely in \cite{Kock:0807}.)
  The important aspects of an element $b\in B$ are: the fibre
  $E_b = p^{-1}(b)$ and the element $j\df t(b) \in J$.  We capture
  these data by picturing $b$ as a (non-planar) bouquet 
  (also called a corolla)
  \begin{center}\begin{texdraw}
    \linewd 0.5 \footnotesize
      \move (0 0) \lvec (0 20)
      \onedot
      \lvec (-10 40) 
      \move (0 20) \lvec (0 40)
      \move (0 20) \lvec (10 40)
      \htext ( 7 18) {$b$}
      \htext ( 0 -8) {$j$}
      \htext ( -12 47) {$e$}
      \htext ( 0 45) {$\ldots$}
      \move (0 53)
  \end{texdraw}\end{center}
  Hence each leaf is labelled by an element $e\in E_b$, and each
  element of $ E_b$ occurs exactly once.  In virtue of the map $s:E\to
  I$, each leaf $e\in E_b$ acquires furthermore an implicit decoration
  by an element in $I$, namely $s(e)$.

  An element in $E$ can be pictured as a bouquet of the same type, but 
  with one of the leaves marked (this mark chooses the element $e \in 
  E_b$, so this description is merely an expression of the natural 
  identification $E 
  = \coprod_{b\in B} E_b$). Then the map $p :E\to B$ consists in forgetting 
  this mark, and $s$ returns the $I$-decoration of the marked leaf.
\end{blanko}

\begin{blanko}{Evaluation of a polynomial functor.}\label{eval}
  Evaluating the polynomial functor $P$ on an object $f:X\to I$
  has the following graphical interpretation. The elements of $P(X)$ are
  bouquets as above, but where each leaf is furthermore decorated by
  elements in $X$ in a
  compatible way:
  \begin{center}\begin{texdraw}
    \linewd 0.5 \footnotesize
      \move (0 0) \lvec (0 20)
      \onedot
      \move (0 20) \lvec (-12 40)  \onedot
      \move (0 20) \lvec (-4.5 43)  \onedot
      \move (0 20) \lvec (4.5 43)  \onedot
      \move (0 20) \lvec (12 40)  \onedot
      \htext ( 7 17) {$b$}
      \htext ( 0 -8) {$j$}
      \htext ( -12 30) {$e$}
      \htext ( -4.3 32) {$\cdot$}
      \htext ( 0 33.4) {$\cdot$}
      \htext ( 4.5 32) {$\cdot$}
      \htext ( -14 48) {$x$}
      \htext ( -6 51) {$\cdot$}
      \htext ( 0 52.5) {$\cdot$}
      \htext ( 6 51) {$\cdot$}
      \htext ( 11 49) {$\cdot$}
      \move (0 61)
  \end{texdraw}\end{center}
  The compatibility condition for the decorations is that a leaf $e$  
  may have decoration $x$ only if $f(x)=s(e)$.
  The set of such $X$-decorated bouquets is
    naturally a set over $J$ via $t$ (return the decoration of the root edge). 
    More formally, $P(X)$ is the set over $B$ (and hence over $J$ via
    $t$) whose fibre over $b\in B$ 
  is the set of commutative triangles
  \begin{diagram}[w=2.2ex,h=4ex,tight]
  X    && \lTo    && E_b    \\
  &\rdTo<{f}    &      & \ldTo>{s}  &  \\
  &    & I .    & &
  \end{diagram}
\end{blanko}

\begin{blanko}{Composition of polynomial functors.}
  The composition of two polynomial functors is again polynomial; this
  is a consequence of distributivity and the Beck-Chevalley conditions
  \cite{Kock:NotesOnPolynomialFunctors}.  We are mostly interested in
  the case $J=I$ so that we can compose $P$ with itself.  The
  composite polynomial functor $P\circ P$ can be described in terms of
  grafting of bouquets: the base set for $P\circ P$, formally
  described as $p\lowerstar (B \times_I E)$, is the set of bouquets of
  bouquets (i.e.~two-level trees)
  \begin{center}\begin{texdraw}
    \linewd 0.5 \footnotesize
      \move (0 0) \lvec (0 20)
      \onedot
      \lvec (-15 40)  \onedot \lvec (-20 60)
      \move (-15 40) \lvec (-10 60)
      \move (0 20) \lvec (0 40) \onedot
      \move (0 20) \lvec (25 40) \onedot
      \lvec (25 60)
      \htext ( 7 18) {$b$}
      \htext ( -13 29) {$e$}
      \htext ( -22 40) {$c$}
      \htext ( 12 40) {$\dots$}
      \move (0 45)
  \end{texdraw}\end{center}
  The conditions on the individual bouquets are still in force: each dot
  is decorated by an element in $B$, and for a dot with decoration $b$
  the set of incoming
  edges is in specified bijection with the fibre $E_b$.  The compatibility
  condition for grafting is this:
  
  \bigskip
  
  \noindent
{\em Compatibility Condition}: for an edge $e$ coming out of
a dot decorated $c$, we have
$$
s(e) = t(c) .
$$

\end{blanko}

  
\begin{blanko}{Morphisms.}
  A cartesian natural transformation $u: P'\Rightarrow P$ between 
  polynomial functors corresponds to a commutative diagram
  \begin{equation}\label{morphism}
  \begin{diagram}[w=2.5ex,h=3ex,tight]
  &&E'\SEpbk&&\rTo^{p'}&&B'\\
  &\ldTo^{s'}&&&&&&\rdTo^{t'}\\
  I&&\dTo&&&&\dTo&&J \\
  &\luTo_s&&&&&&\ruTo_t\\
  &&E&&\rTo_p&&B\\
  \end{diagram}
  \end{equation}
  whose middle square is cartesian,
  cf.~\cite{Kock:NotesOnPolynomialFunctors}.  In other words, giving
  $u$ amounts to giving a $J$-map $u:B'\to B$ together with an $I$-bijection
  $E'_{b'} \isopil E_{u(b')}$ for each $b'\in B'$.
  
  Let $\Poly(I)$ denote the category whose objects are the polynomial
  endofunctors on $\Set/I$ as in (\ref{P}) and whose arrows are the
  cartesian natural transformations as in (\ref{morphism}).  This is a
  strict monoidal category under composition, and with the identity functor
  $\Id$ as unit object.  Note that a polynomial functor always preserves
  cartesian squares, and (under the assumption $E\to B$ finite)
  sequential colimits
  \cite{Kock:NotesOnPolynomialFunctors}.
\end{blanko}

\begin{blanko}{Polynomial monads.}
  By a {\em polynomial monad} we understand a polynomial endofunctor $P:\Set/I \to
  \Set/I$ with monoid structure in $\Poly(I)$.  
  In other words, there is specified a composition law
  $\mu:P\circ P \to P$ with unit $\eta: \Id \to P$, satisfying the usual
  associativity and unit conditions, and $\mu$ and $\eta$ are cartesian
  natural transformations.  Throughout we indicate monads by their functor
  part, confident that in each case it is clear what the
  natural-transformation part is, or explicitating it otherwise.

  The composition law is described graphically as an operation of
  contracting two-level trees (formal compositions of bouquets) to bouquets.
  
  We shall refer to $I$ as the set of {\em types} of $P$, denoted 
  $\operatorname{typ}(P)$, and $B$ as 
  the set of {\em operations}, denoted $\operatorname{op}(P)$.
  Since we have a unit, we can 
  furthermore think of $E$ as the set of {\em partial operations}, 
  i.e.~operations all of whose inputs except one are fed with a unit.
  The composition law can be described in terms of partial operations
  as a map
  $$B \times_I E \to B,$$
  consisting in substituting one operation into one input of 
  another operation, provided the types match: $t(b)=s(e)$.
\end{blanko}


\begin{blanko}{The free monad on a polynomial endofunctor.}
  (See also Gambino-Hyland~\cite{Gambino-Hyland}.)
  Given a polynomial endofunctor $P : \Set/I \to \Set/I$, a {\em 
  $P$-set} is a pair $(X,a)$ where $X$ is an object of $\Set/I$ and
  $a: P(X) \to X$ is an arrow in $\Set/I$ (not subject to any further 
  conditions).  A $P$-map from $(X,a)$ to $(Y,b)$ is an arrow $f: X\to Y$
  giving a commutative diagram
  \begin{diagram}[w=6ex,h=4.5ex,tight]
  P(X) & \rTo^{P(f)}  & P(Y)  \\
  \dTo<a  &    & \dTo>b  \\
  X  & \rTo_f  & Y .
  \end{diagram}
  Let $P\kat{-Set}/I$ denote the category of $P$-sets and $P$-maps.
  The forgetful functor $U : P\kat{-Set}/I \to \Set/I$ has a left adjoint
  $F$, the {\em free $P$-set functor}.  The monad $P\upperstar \df U \circ F : \Set/I
  \to \Set/I$ is the {\em free monad} on $P$.  This is a polynomial monad,
  and its set of operations is the set of $P$-trees, as we now explain.
\end{blanko}

\begin{blanko}{$P$-trees.}\label{P-trees}
  Let $P$ denote a polynomial endofunctor given by $I \leftarrow E \to
  B \to I$.  We define a {\em $P$-tree} to be a tree whose edges are
  decorated in $I$, whose nodes are decorated in $B$, and with the
  additional structure of a bijection for each node $n$ (with
  decoration $b$) between the set of input edges of $n$ and the fibre
  $E_b$, subject to the compatibility condition that such an edge
  $e\in E_b$ has decoration $s(e)$, and the output edge of $n$ has
  decoration $t(b)$.  Note that the $I$-decoration of the edges is
  completely determined by the node decoration together with the
  compatibility requirement, except for the case of a unit tree.
  
  Another description is useful: a $P$-tree is a tree with edge set 
  $A$, node set $N$, and node-with-marked-input-edge set $N'$, 
  together with a diagram
  \begin{diagram}[w=5ex,h=4ex,tight]
  A &\lTo & N' \SEpbk & \rTo & N & \rTo & A  \\
  \dTo<\alpha  &    & \dTo && \dTo>\beta  &    & \dTo>\alpha  \\
  I &\lTo  & E  & \rTo & B & \rTo & I .
  \end{diagram}
  Then the vertical maps $\alpha$ and $\beta$ express the decorations, and
  the commutativity and the cartesian condition on the middle square
  express the bijections and the compatibility condition.  The top row is a
  polynomial functor associated to a tree, and in short, a $P$-tree can be
  seen as a cartesian morphism from a tree to $P$ in a certain category of
  polynomial endofunctors \cite{Kock:0807}.
  
  An {\em isomorphism of $P$-trees} is an isomorphism of trees compatible 
  with the $P$-decorations.  It is clear that $P$-trees are rigid
  Denote by $\tr(P)$ the set of isomorphism classes of $P$-trees.
  This is the set of formal combinations of the operations of $P$, 
  i.e.~obtained by freely grafting elements of $B$ onto the
  leaves of elements of $B$, provided the decorations match (and formally
  adding a unit tree for each $i\in I$). 
  The set $\tr(P)$ has a natural map to $I$ by returning the root, and it can 
  be described as a least
  fixpoint for the polynomial endofunctor
  \begin{eqnarray*}
    \Set/I & \longrightarrow & \Set/I  \\
    X & \longmapsto & I + P(X);
  \end{eqnarray*}
  as such it is given explicitly as the colimit
  $$
  \tr(P) = \bigcup_{n\in \N} (I+P)^n(\emptyset) .
  $$
\end{blanko}

\begin{blanko}{Explicit description of the free monad on $P$.}
  A slightly more general fixpoint construction characterises the free 
  $P$-set monad $P\upperstar $:
  if $A$ is an object of 
  $\Set/I$, then $P\upperstar (A)$ is a least fixpoint for the 
  endofunctor $X\mapsto A+P(X)$.  In explicit terms,
  $$
  P\upperstar (A) = \bigcup_{n\in \N} (A+P)^n(\emptyset) .
  $$
  It is the set of (isomorphism classes of) $P$-trees with leaves decorated in $A$.  
  But this
  is exactly the characterisation of evaluation of a polynomial functor 
  (\ref{eval}) with operation set $\tr(P)$: let $\tr'(P)$
  denote the set of (isomorphism classes of) 
  $P$-trees with a marked leaf, then $P\upperstar : \Set/I \to
  \Set/I$ is the polynomial functor given by
  \polyFunct{I}{\tr'(P)}{\tr(P)}{I \ .}{}{}{}{} The
  maps are the obvious ones: return the marked leaf, forget the mark, and
  return the root edge, respectively.  The monad structure of $P\upperstar$ is
  described explicitly in terms of grafting of trees.  In a
  partial-composition description, the composition law is
  $$
  \tr(P)\times_I \tr'(P) \to \tr(P)
  $$
  consisting in grafting a tree onto the specified input leaf of
  another tree.  The unit is given by $I \to \tr(P)$ associating
  to $i\in I$ the unit tree with edge decorated by $i$.  (One can
  readily check that this monad is cartesian.)
\end{blanko}


\section{The Baez-Dolan construction for polynomial monads}
\label{Sec:BD}

Throughout this section, we fix a polynomial monad $P : \Set/I \to \Set/I$,
represented by \polyFunct{I}{E}{B}{I.}{}{}{}{} We shall associate to the
polynomial monad $P : \Set/I \to \Set/I$ another polynomial monad $P^+: \Set/B
\to \Set/B$.  The idea of this construction is due to Baez and
Dolan~\cite{Baez-Dolan:9702}, who realised it in the settings of symmetric
operads.  We first give a very explicit version for polynomial monads, and
show how to produce the opetopes from it by iteration, recovering the elementary
definition of opetopes given in \ref{def}.  It is the graphical interpretation
of polynomial functors that allows us to extract the combinatorics.
Afterwards we compare with Leinster's definition of opetopes \cite[\S
7.1]{Leinster:0305049}.  This is just a question of comparing our version of
the Baez-Dolan construction with Leinster's; the iterative construction 
of opetopes is exactly the same.

\subsection{Explicit construction}

\begin{blanko}{The Baez-Dolan construction for a polynomial monad.}\label{BD}
  Starting from our polynomial monad $P$, we describe explicitly a new
  polynomial monad $P^+$, the {\em Baez-Dolan construction} on $P$.
  The idea is to substitute into dots of trees instead
  of grafting at the leaves (so notice that this shift is like in a zoom
  relation).  Specifically,
  define $\tr\upperdot(P)$ to be the set of (isomorphism classes of) 
  $P$-trees with one marked dot.
  There is now a polynomial functor
  \polyFunct{B}{\tr\upperdot(P)}{\tr(P)}{B}{}{}{t}{P^+} 
  where $\tr\upperdot(P)
  \to \tr(P)$ is the forgetful map, $\tr\upperdot(P)\to B$ returns the
  bouquet around the marked dot, and $t:\tr(P)\to B$ comes from the
  monad structure on $P$: it amounts to contracting all inner edges
  (or setting a new dot in a unit tree).
  Graphically:
  
  \begin{equation}\label{slice1}\begin{texdraw}
    \move (-78 8) 
    \bsegment 
      \move (0 5) \lvec (0 17) \onedot
      \lvec (-10 30) \onedot  
      \lvec (-30 35)
      \move (-10 30) \lvec (-20 50) \onedot 
      \move (-10 30) \lvec (-5 45) \onedot \lvec (-10 65)
      \move (-5 45) \lvec (0 65)
      \move (0 17) \lvec (15 60)
      \move (0 17) \lvec (22 35) \onedot
      \htext (-15 24) {*}
    \esegment
    \htext (-120 40) {$\bigleftbrace{30}$}
    \htext (-40 40) {$\bigrightbrace{30}$}
	  
    \move (82 8)
    \bsegment 
      \move (0 5) \lvec (0 17) \onedot
      \lvec (-10 30) \onedot  
      \lvec (-30 35)
      \move (-10 30) \lvec (-20 50) \onedot 
      \move (-10 30) \lvec (-5 45) \onedot \lvec (-10 65)
      \move (-5 45) \lvec (0 65)
      \move (0 17) \lvec (15 60)
      \move (0 17) \lvec (22 35) \onedot
    \esegment
    \htext (40 40) {$\bigleftbrace{30}$}
    \htext (120 40) {$\bigrightbrace{30}$}

    \move (-140 -70)
    \bsegment 
      \move (0 5) \lvec (0 17) \onedot
      \lvec (-12 27)
      \move (0 17) \lvec (0 32)
      \move (0 17) \lvec (12 27)
    \esegment
    \htext (-160 -50) {$\bigleftbrace{15}$}
    \htext (-120 -50) {$\bigrightbrace{15}$}

    \move (140 -70)
    \bsegment 
      \move (0 5) \lvec (0 17) \onedot
      \lvec (-12 27)
      \move (0 17) \lvec (-5 30)
      \move (0 17) \lvec (5 30)
      \move (0 17) \lvec (12 27)
    \esegment
    \htext (120 -50) {$\bigleftbrace{15}$}
    \htext (160 -50) {$\bigrightbrace{15}$}
       
    \htext (0 -55) {$P^+$}
    \htext (120 -5) {\footnotesize $t$}
    \arrowheadtype t:F
    \arrowheadsize l:5 w:3

    \move (-15 40) \avec (15 40)
    \move (110 0) \avec (125 -23)
    \move (-110 0) \avec (-125 -23)
  \end{texdraw}\end{equation}
  (In this diagram as well as in the following diagrams of the same
  type, a symbol \raisebox{-4pt}{
  \begin{texdraw} \setunitscale 0.6
    \move (0 -17) 
    \bsegment 
     \move (0 5) \lvec (0 17) \onedot \lvec (-12 27) \move (0 17)
    \lvec (-5 30) \move (0 17) \lvec (5 30) \move (0 17) \lvec (12 27)
    \esegment 
    \htext (-20 0) {$\bigleftbrace{10}$} 
    \htext (20 0) {$\bigrightbrace{10}$} 
  \end{texdraw} }
  is meant to designate the set
  of {\em all} bouquets like this (with the appropriate decoration),
  but at the same time the specific figures representing each set are
  chosen in such a way that they match under the structure maps.)
  Note that since the forgetful map forgets a marked dot, the nullary 
  operations in $P^+$ are precisely the unit trees \inlineDotlessTree,
  one for each $i\in I$.
\end{blanko}
  
\begin{blanko}{Monad structure on $P^+$.}
  We first compute the value of $P^+$ on an object $C\to B$ of $\Set/B$.  Using
  the explicit graphical description of evaluation of a polynomial functor
  \ref{eval}, we see that the result is the set of $P$-trees with each node
  decorated by an element of $C$, compatibly with the arity map $C \to B$ (being
  a $P$-tree means in particular that each node already has a $B$-decoration;
  these decorations must match).
  
  We can now compute $P^+\circ P^+$: its set of operations is $P^+$ evaluated at
  $t: \tr(P)\to B$: that's the set of (isomorphism classes of) $P$-trees with nodes decorated by
  $P$-trees in such a way that the total bouquet of the decorating tree matches
  the local bouquet of the node it decorates.  Similarly, the set of `partial
  operations' for $P^+\circ P^+$ is the set of $P$-trees-with-a-marked-node,
  the marked node being decorated with a $P$-tree-with-a-marked-node, and the 
  remaining nodes being decorated by $P$-trees.
  
  Now the monad structure on $P^+$ is easy to describe:
   The composition law $P^+\circ
  P^+ \Rightarrow P^+$ consists in substituting each $P$-tree into the node it
  decorates.
  The substitution can be described in terms of a partial
  composition law
  $$
  \tr(P)\times_B \tr\upperdot(P)\to
  \tr(P)$$
  defined by substituting a $P$-tree into the marked dot of an 
  element in $\tr\upperdot(P)$, as indicated in this
  figure:
  \begin{equation}\label{BD-subst}\begin{texdraw}
  \scriptsize
    \move (0 0)
    \bsegment 
    \htext (0 -20) {\normalsize $F$}
    \linewd 0.3 
    \move (0 0) \lvec (0 20) \onedot
    \lvec (-15 30) \onedot
    \move (0 20) \lvec (-17 65) \htext (-19 70) {$x$}
    \move (0 20) \lvec (-5 70) \htext (-6 75) {$y$}
    \move (0 20) \lvec (14 35) \onedot
    \lvec (9 50) \onedot \lvec (9 70) \htext (9 75) {$z$}
    \move (14 35) \lvec (23 47) \onedot
    \lpatt (2 2)
    \move (5 35)  \freeEllipsis{29}{19}{35}
    \lpatt (2 2)
    \move (20 20) \clvec (60 -15)(90 -15)(110 15)
    \lpatt ()
    \arrowheadtype t:F
    \arrowheadsize l:5 w:3
    \avec (112.5 18.5)
    \esegment

    \move ( 130 0)
    \bsegment 
    \linewd 0.3 
    \move (0 0) \lvec (0 17) \onedot
    \lvec (-10 30) \onedot  
    \lvec (-30 35) \htext (-35 36) {$x$}
    \htext (-14 25) {$f$}
    \move (-10 30) \lvec (-20 55) \onedot \htext (-21 47) {$y$}
  \move (-10 30) \lvec (-5 50) \onedot \lvec (-10 70)
  \move (-5 50) \lvec (0 70)
  \move (0 17) \lvec (15 60) 
  \move (0 17) \lvec (22 35) \onedot 
  \htext (-3 43) {$z$}
  \lpatt (2 2) \move (-12 29) \lcir r:11
    \esegment
    
\htext (240 15) {\normalsize  resulting in }

\move (340 0) 
   \scriptsize
   \bsegment 
   \linewd 0.3 
   \move (0 0) \lvec (0 17) \onedot
   \lvec (-15 35) \onedot  
   \lvec (-35 70) \htext (-38 73) {$x$}
   \move (-15 35) \lvec (-30 40) \onedot
   \move (-15 35) \lvec (-25 85) \onedot \htext (-26 71) {$y$}
 \move (-15 35) 
 \lvec (-5 45) \onedot 
 \lvec (-5 60) \onedot \lvec (-5 80)  \onedot \lvec ( -10 95) 
 \move (-5 80) \lvec (0 95)

 \htext (-8 72) {$z$}
 \move (-5 45) \lvec (5 55) \onedot
 \move (0 17) \lvec (25 60) 
 \move (0 17) \lvec (22 35) \onedot 
 \lpatt (2 2)    \move (-12 47)  \freeEllipsis{26}{17}{32}

   \esegment

   \move (390 0)
 \end{texdraw}
 \end{equation}
  (The
  letters in the figure do not represent the decorations --- they are
  rather unique labels to express the involved bijections, and to
  facilitate comparison with Figure~(\ref{BD-subst-spheres}) below.)
  Of course the
  substitution makes sense only if the decorations match. 
  This means that $t(F)$, the `total bouquet' of the tree $F$, is the 
  same as the local bouquet of the node $f$.  Formally the substitution can
  be described as a pushout in a category of $P$-trees, cf.~\cite{Kock:0807}.
  
  The unit for the monad is given by the map $B \to \tr(P)$
  interpreting a bouquet as a tree with a single dot.  
  
  It is readily checked directly that the monad axioms hold.  (Alternatively
  this will follow from the proof of Theorem~\ref{ours=leinsters} where $P^+$ is
  shown isomorphic to something which is a monad by construction.)
\end{blanko}

\begin{blanko}{The BD construction in terms of nestings.}
  We have described the free-monad construction and the Baez-Dolan
  construction in terms of trees, but of course they can equally well
  be described in terms of nested spheres, as we shall now explain.
  The interplay between these two descriptions will lead directly to
  opetopes as defined in Section~\ref{Sec:opetopes}.
  Let us stress again that trees and nestings are just 
  different graphical expressions of the same combinatorial structure.
  However, some features of trees can be a little bit subtler to see
  in terms of nestings.
  
  The basic operations, the elements in $B$, are configurations of a sphere 
  with dots inside:
  \begin{center}
    \begin{texdraw}
      \scriptsize
      \bsegment
	\move (-2 20) \lcir r:20
	\move (7 24) \onedot
	\move (0 12) \onedot
	\move (-4 25) \onedot
	\htext (14 3) {$j$}
	\htext (0 20) {$e_1\dots$}
	\htext (-17 20) {$b$}
      \esegment
    \end{texdraw}
  \end{center}
  We call such a thing a {\em layer}.  The set of dots inside the
  sphere is in bijection with the set $E_b$, and via $s:E\to I$ these
  dots also carry an implicit decoration by elements in $I$, the input
  types.  The label $j$ on the outside of the sphere represents
  $t(b)$, the output.  We put the label $b$ on the inside of the
  sphere it decorates, since it mediates between the input devices
  (the dots) and the output device (the sphere), just as the dot of a
  bouquet mediates between the inputs (the leaves) and the output.

  Next, $\tr(P)$ is the set of (isomorphism classes of) arbitrary $P$-nestings,
  with layers decorated in $B$ and spheres and dots decorated in $I$ (subject to
  compatibility conditions), and $\tr'(P)$ is the set of (isomorphism classes
  of) arbitrary $P$-nestings (compatibly decorated) with a marked dot.  The
  substitution law for the free monad on $P$ is now described by substituting
  one $P$-nesting into a dot of another, provided the decorations match.  (This
  corresponds to grafting of trees.)
  
  For the Baez-Dolan construction (where we now suppose $P$ is a monad),
  $\tr\upperdot(P)$ is the set of (isomorphism classes of) $P$-nestings with a
  marked sphere, so here is the nesting version of Figure~(\ref{slice1}):

  \begin{center}\begin{texdraw}
    \footnotesize
    \move (-80 8) 
    \bsegment
      \move (0 30) \lcir r:30
      \move (-20 25) \lcir r:5
      \move (8 32) \lcir r:19
      \move (16 38) \lcir r:5
      \move (7 24) \lcir r:9
      \move (-14 14) \onedot
      \move (4 25) \onedot
      \move (10 23) \onedot
      \move (3 40) \onedot
      \htext (6 46.5) {\normalsize *}
    \esegment
    \htext (-120 40) {$\bigleftbrace{30}$}
    \htext (-40 40) {$\bigrightbrace{30}$}
	  
    \move (80 8)
    \bsegment
      \move (0 30) \lcir r:30
      \move (-20 25) \lcir r:5
      \move (8 32) \lcir r:19
      \move (16 38) \lcir r:5
      \move (7 24) \lcir r:9
      \move (-14 14) \onedot
      \move (4 25) \onedot
      \move (10 23) \onedot
      \move (3 40) \onedot
    \esegment
    \htext (40 40) {$\bigleftbrace{30}$}
    \htext (120 40) {$\bigrightbrace{30}$}

    \move (-140 -70)
    \bsegment
      \move (0 20) \lcir r:16
      \move (7 24) \onedot
      \move (0 12) \onedot
      \move (-4 25) \onedot
    \esegment
    \htext (-165 -50) {$\bigleftbrace{17}$}
    \htext (-115 -50) {$\bigrightbrace{17}$}

    \move (140 -70)
    \bsegment
      \move (0 20) \lcir r:16
      \move (-6 13) \onedot
      \move (1 18) \onedot
      \move (7 16) \onedot
      \move (1 27) \onedot
    \esegment
    \htext (115 -50) {$\bigleftbrace{17}$}
    \htext (165 -50) {$\bigrightbrace{17}$}
    
    \htext (0 -55) {$P^+$}
    \arrowheadtype t:F
    \arrowheadsize l:5 w:3

    \move (-15 40) \avec (15 40) 
    \move (-110 0) \avec (-125 -23) \htext (-122 -8){$s$}
    \move (110 0) \avec (125 -23) \htext (122 -8){$t$}
  \end{texdraw}\end{center}
  Note that the map $t$ consists in erasing all inner spheres,
  which is just the nesting equivalent of the tree operation
  of contracting all inner edges --- this is always possible for undecorated 
  nestings, but for this to make sense in the $P$-decorated 
  case we need the monad structure on $P$. The map $s$ consists in returning
  the {\em layer} determined by the marked sphere: this means the
  region delimited on the outside by the marked sphere itself and on the
  inside by its children, so the operation
  can also be described as taking the marked sphere and contracting each
  sphere inside it to a dot.  (This is the nesting equivalent
  of
  the tree operation of returning the `local bouquet' of a dot.)
  
  The substitution law is perhaps less obvious in this
  nesting interpretation.  Looking at Figure~(\ref{BD-subst})
  we see that for trees the substitution takes place at a specified
  dot, and consists in replacing its `local bouquet' by a more
  complicated tree,  so the operation is about refining the tree.
  Correspondingly for nestings, the operation is about
  refining the nesting by drawing some more spheres in the
  specified layer.  Here is the nesting version of
  Figure~(\ref{BD-subst}):
  \begin{equation}\label{BD-subst-spheres}\begin{texdraw}
  \scriptsize
    \move (0 0)
    \bsegment
    \move (-3 68) \lcir r:36
    \move (2.5 57.5) \freeEllipsis{12}{8}{66}
    \move (-15 90) \onedot
    \htext (-21 91) {$b$}

    \linewd 0.3 
    \move (-25 55) \onedot
    \htext ( -26 50) {$a$}
    \move  (5 58) \onedot 
  \htext (0 58) {$c$}
  \move (14 76) \lcir r:5
  \move (6 63) \freeEllipsis{26}{16}{63}
  \move (-28 72) \lcir r:5
  \esegment

    \move ( 130 0)
    \bsegment
    \move (-4 69) \fcir f:0.92 r:35
    \move (-4 69) \lcir r:35
    \move (2.5 57.5) \freeFillEllipsis{14}{8}{70}{1}
    \move (2.5 57.5) \freeEllipsis{14}{8}{70}
    \move (-15 81) \fcir f:1 r:8 
    \move (-15 81) \lcir r:8
    
    \linewd 0.3 
    \move (5 25) \onedot
    \move (0 50) \onedot \move (-21 55) \onedot
    \htext ( -22 50) {$a$}
    \move ( 25 35) \lcir r:5 
    \move  (5 65) \onedot 
    \htext (-18 81) {$b$}
	      \move (1 64) \lcir r:48
  \htext (35 58) {$f$}
  \htext (-1 59) {$c$}
    \esegment
	
    \move (0 0)
    \bsegment
  
    \lpatt (2 2)
      \move (24 92) \clvec (61 125)(90 125)(108 98)
      \lpatt ()
      \arrowheadtype t:F
      \arrowheadsize l:5 w:3
      \avec (112 93)

  \esegment
\htext (215 35) {\normalsize  giving }

\move (300 0) 
   \scriptsize
   \bsegment
   \move (-2 71) \fcir f:0.92 r:42
   \move (-2 71) \lcir r:42
   \move (2.5 57.5) \freeFillEllipsis{14}{8}{70}{1}
   \move (2.5 57.5) \freeEllipsis{14}{8}{70}
   \move (-20 90) \fcir f:1 r:8 
   \move (-20 90) \lcir r:8
   
   \linewd 0.3 
   \move (5 20) \onedot
   \move (0 50) \onedot \move (-25 55) \onedot
   \htext ( -26 50) {$a$}
   \move ( 28 30) \lcir r:5 
   \move  (5 65) \onedot 
   \htext (-23 91) {$b$}
	     \move (3 68) \lcir r:55
 \htext (43 58) {$f$}
 \htext (0 59) {$c$}
 \linewd 1
 \move (2.5 57.5) \freeEllipsis{18}{12}{70}
 \move (16 83) \lcir r:5
 \move (8 66) \freeEllipsis{32}{20}{63}
 \move (-30 72) \lcir r:5
   \esegment
  
   \move (390 40)
 \end{texdraw}
 \end{equation}
Again, the $B$-decorations have not been drawn; the letters serve only to 
specify the bijections, and to facilitate comparison with
Figure~(\ref{BD-subst}).

\end{blanko}

\begin{blanko}{The double Baez-Dolan construction (slice-twice construction).} 
  \label{slice-twice}%
  After applying the Baez-Dolan construction once (in
  its tree interpretation), we have a polynomial functor $B \leftarrow
  \tr\upperdot(P)\to \tr(P)\to B$ which is a monad for the
  operation of substituting one tree into a dot of another tree
  (subject to some book-keeping). Applying the construction a second time we get
  \begin{diagram}[w=2.2ex,h=3.7ex,tight]
&&\tr\upperdot(P^+)&&&\rTo&&&\tr(P^+)\\
&\ldTo^{s}&&&&&&&&\rdTo^{t}\\
\tr(P)&&&&&P^{++}&&&&&\tr(P)
\end{diagram}

  Let us spell out the details.  Unwinding the definitions, a $P^+$-tree is a
  tree $M$ whose dots are decorated by $P$-trees, and whose edges are decorated
  by elements in $B$, and with a specified bijection, for each node $n$ with
  decorating $P$-tree $T$, between the set of input edges of $n$ and the set of dots
  in $T$.  The decoration of such an input edge must be exactly the
  corresponding dot in $T$, interpreted as an element in $B$, and the
  output edge of a dot decorated by $T$ must be decorated by the total bouquet
  of $T$ (i.e.~the element of $B$ obtained by contracting all inner edges of $T$
  using the monad structure of $P$).
  The description of the elements in $\tr\upperdot(P^+)$ is similar,
  but with one node in $M$ marked.  The map $\tr\upperdot(P^+) \to \tr(P)$
  returns the $P$-tree decorating the marked node.
  
  The map $\tr(P^+) \to \tr(P)$ involves the monad law for $P^+$.  Namely, we
  contract each inner edge of $M$, by composing the two $P$-trees decorating the
  adjacent dots.  According to the composition law for $P^+$, this means 
  substituting the upper decorating $P$-tree into the designated dot of the 
  lower decorating $P$-tree.  (The designated dot is the one corresponding to the
  edge of $M$ we are contracting, and the substitution makes sense because of 
  the compatibility requirement of the decoration of $M$.)  In other words,
  this $P$-tree is obtained by successively substituting all the decorating 
  $P$-trees into each other according to the recipe specified by the tree $M$.

  Here is a drawing illustrating the notion of $P^+$-tree:
  \begin{center}\begin{texdraw}
    \linewd 0.5 \footnotesize
    \move (-20 -10)
    \bsegment
      \move (-5 0) \lvec (-5 20)
      \onedot
      \lvec (-15 40)  \onedot \lvec (-32 58) \onedot
      \move (-15 40) \lvec (-15 65)
      \move (-15 40) \lvec (0 58)
      \move (-5 20) \lvec (25 56) 
      \htext (-39 58) {$a$}
      \htext (-22 55) {$1$}
      \htext (-11 59) {$1$}
      \htext (-2 48) {$3$}
      \htext ( -6 32) {$3$}
      \htext (16 38) {$3$}
      \htext ( -22 40) {$b$}
      \htext (-11 19) {$c$}
      \htext (-1 10) {$5$}
      
      \htext ( -5 -10) {$M$}
\esegment
\htext (75 30) {with} 
      \move (155 40)
      \bsegment
      \bsegment
      \move (0 0) \lvec (0 15) \move (0 7) \onedot
            \htext (0 -7) {$1$}
	    \esegment
	    
      \move (20 0) 
      \bsegment
      \move (0 0)
      \lvec (0 15) \move (0 7) \onedot \lvec (-5 14) \move (0 7) 
      \lvec (5 14)
            \htext (0 -7) {$3$}
	    \esegment
	    \move (40 0)
	    \bsegment
      \move (0 0) \lvec (0 15) \move (0 7) \onedot \lvec (-7 11) \move (0 7) 
      \lvec (-4 14) \move (0 7) 
      \lvec (4 14) \move (0 7) 
      \lvec (7 11)
	    \htext (0 -7) {$5$}
	    \esegment
	    
	    \htext (80 3) {$\in \ B$,}
\esegment

\move ( 120 -20)

\bsegment
\move (-5 0)
  \bsegment
    \move (0 0)
    \lvec (0 30)
    \htext (0 -8) {$a$}
  \esegment
  \move (28 0)
  \bsegment
    \move (0 0)
    \lvec (0 30)
    \move (0 10) \onedot \lvec (-10 30) \move (-5 20) \onedot
    \move (0 10) \onedot \lvec (10 30) \move (5 20) \onedot
    \htext (0 -7) {$b$}
  \esegment
  \move (70 0)
  \bsegment
    \move (0 0)
    \lvec (0 10) \onedot
    \lvec (-10 20) \onedot \lvec (-18 30) \move (-10 20) \lvec (-10 32)
    \move (-10 20) \lvec (-2 30)
    \move (0 10)  \lvec (5 32) \move (0 10) \lvec (15 30)
    \htext (0 -8) {$c$}
  \esegment
	    \htext (120 12) {$\in \ \tr(P)$}
\esegment
  \end{texdraw}\end{center}
And here is the result of applying $t$ to it:
  \begin{center}\begin{texdraw}
    \linewd 0.5 \footnotesize
    \move (0 0)
    \bsegment
      \move (0 0) \lvec (0 20)
      \onedot
      \lvec (-15 40)  \onedot \lvec (-33 70) \move (-24 55) \onedot
      \move (-15 40) \lvec (-15 80)
      \move (-15 40) \lvec (4 82)
      \move (0 20) \lvec (25 75) 
      \move (0 20) \lvec (33 65) 
               \linewd 0.3 
 \lpatt (2 2)    
 \move (-6 55) \lcir r:6.5 \htext (-4 54) {$a$}

 \move (-13 50) \freeEllipsis{20}{17}{32} \htext ( -26 43) {$b$}

 \move (-5 45) \lcir r:33 

      \htext (-19 20) {$c$}
      
      \htext (-2 -10) {$t(M)$}
\esegment
  \end{texdraw}\end{center}
  Here the dashed spheres are drawn to indicate how the original $P$-trees $a$, $b$,
  and $c$ were substituted into each other:  the inner spheres represent the
  `scars' of the two substitutions, $a$ into a certain node of $b$, and $b$ into
  a certain node of $c$.  The outer sphere represents the tree $c$,
  corresponding to the `root dot' of $M$.  Altogether we see a constellation
  whose underlying nesting is precisely $M$, and
  whose underlying tree is a $P$-tree.
  
  This is general: the elements in $\tr(P^+)$ are obtained
  by successive substitutions of $P$-trees into nodes of a $P$-tree, and if for
  each such substitution we keep track of the surgery via the scar it left ---
  that's a sphere in the tree --- we obtain a $P$-constellation.  This is the
  content of the following theorem which also tells us that the $P^+$-tree can
  be recovered from the $P$-constellation.
\end{blanko}

\begin{blanko}{The $P$-constellation monad.}
  By a {\em $P$-constellation} we mean a constellation whose underlying tree
  is a $P$-tree. Let $\const(P)$ denote the set of isomorphism classes of 
  $P$-constellations (note that $P$-constellations are rigid objects).  
  Similarly, let $\const\uppercirc(P)$ denote the set of isomorphism classes of 
  $P$-constellations with a marked layer.

  Define a polynomial endofunctor by
  \begin{diagram}[w=2.5ex,h=3.7ex,tight]
&&\const\uppercirc(P)&&&\rTo&&&\const(P)\\
&\ldTo&&&&&&&&\rdTo\\
\tr(P)&&&&&&&&&&\tr(P)
\end{diagram}
Graphically,
  \begin{equation}\label{const-monad}\begin{texdraw}
    \move (-80 8) 
    \bsegment
    \setunitscale 0.7
    \move (0 -5) \lvec (0 25) \onedot
    \lcir r:10
    \lvec (-25 45) \onedot
    \lvec (-50 50) 
    \move (-25 45) \lvec (-20 95)
    \move (0 25) 
    \lvec (50 55)
    \move (20 37) \lcir r:6
    \move (0 25) \lvec (5 60) \onedot
    \lvec (10 95) 
    \move (-23 65) \lcir r:6
    \move (9 30) \freeEllipsis{27}{16}{30}
    \move (0 46) \lcir r:42
    \htext (17 18) {\normalsize *}
    \esegment
    \htext (-120 40) {$\bigleftbrace{30}$}
    \htext (-40 40) {$\bigrightbrace{30}$}
	  
    \move (80 8)
    \bsegment
    \setunitscale 0.7
    \move (0 -5) \lvec (0 25) \onedot
    \lcir r:10
    \lvec (-25 45) \onedot
    \lvec (-50 50) 
    \move (-25 45) \lvec (-20 95)
    \move (0 25) 
    \lvec (50 55)
    \move (20 37) \lcir r:6
    \move (0 25) \lvec (5 60) \onedot
    \lvec (10 95) 
    \move (-23 65) \lcir r:6
    \move (9 30) \freeEllipsis{27}{16}{30}
    \move (0 46) \lcir r:42
    \esegment
    \htext (40 40) {$\bigleftbrace{30}$}
    \htext (120 40) {$\bigrightbrace{30}$}

    \move (-143 -75)
    \bsegment
    \setunitscale 0.9
    \move (0 0) \lvec (0 15) \onedot
    \lvec (-17 30)
    \move (0 15) 
    \lvec (30 33)
    \move (15 24) \onedot
    \move (0 15) \lvec (3 47)
    \esegment
    \htext (-170 -55) {$\bigleftbrace{25}$}
    \htext (-110 -55) {$\bigrightbrace{25}$}

    \move (140 -75)
    \bsegment
    \setunitscale 0.8
    \move (0 -6) \lvec (0 12) \onedot
    \lvec (-15 24) \onedot
    \lvec (-30 27) 
    \move (-15 24) \lvec (-12 54)
    \move (0 12) 
    \lvec (30 30)
    \move (0 12) \lvec (3 33) \onedot
    \lvec (6 54) 
    \esegment
    \htext (110 -55) {$\bigleftbrace{25}$}
    \htext (170 -55) {$\bigrightbrace{25}$}
       

    \arrowheadtype t:F
    \arrowheadsize l:5 w:3

    \footnotesize
    
    \move (-15 40) \avec (15 40)
    \move (-110 0) \avec (-125 -23) \htext (-122 -8){$s$}
    \move (110 0) \avec (125 -23) \htext (122 -8){$t$}
  \end{texdraw}\end{equation}
  The structure maps are: $t$ returns the
  underlying tree of a constellation, and $s$ returns the tree
  contained in the marked layer.  The monad structure consists in
  substituting one constellation into the marked layer of another,
  provided of course their decorations match.
\end{blanko}

\begin{satz}\label{Thm:slicetwice}
  There is a natural bijection $\tr(P^+) = \const(P)$.  This bijection is 
  compatible with the structure maps described above, yielding an isomorphism
  of polynomial monads
  \begin{equation}\label{tr-const-diagram}
  \begin{diagram}[w=3.5ex,h=3.5ex,tight]
  &&\const\uppercirc(P)\SEpbk&&\rTo&&\const(P)\\
  &\ldTo&&&&&&\rdTo\\
  \tr(P)&&\dLig&&&&\dLig&&\tr(P) \\
  &\luTo&&&&&&\ruTo\\
  &&\tr\upperdot(P^+)&&\rTo&&\tr(P^+)\\
  \end{diagram}
  \end{equation}
\end{satz}

\begin{dem}
  {\em From $P$-constellation to $P^+$-tree.} Given a constellation
  $C$, we first get an abstract tree $M$ by taking the tree
  corresponding to the underlying nesting of $C$,
  cf.~\ref{corr}.  Let $L$ denote the set of layers, and $S$ the set
  of spheres and dots.  To each layer we associate its outside sphere
  (the output sphere), hence a map $L \to S$.  Let $\ov L$ denote the
  set of layers with a marked child, and consider the forgetful map to
  $L$; finally there is the obvious map $\ov L \to S$ returning the
  marked child.  These maps,
  $$
  S \leftarrow \ov L \to L \to S
  $$
  is the polynomial functor associated to the tree $M$ as in \ref{P-trees}.
  We must now decorate this
  tree by $P^+$, i.e., provide a diagram
  \begin{diagram}[w=5ex,h=4ex,tight]
  S & \lTo  & \ov L \SEpbk& \rTo & L & \rTo & S  \\
  \dTo<\alpha  & (3)   & \dTo<\gamma & (2)& \dTo>\beta  &  (1)  & \dTo>\alpha \\
  B  & \lTo  & \tr\upperdot(P) & \rTo & \tr(P) & \rTo & B .
  \end{diagram}
  
  To define $\alpha$: to each dot of $C$ we associate its local
  bouquet in the underlying $P$-tree of $C$.  To each sphere of $C$,
  intuitively we can just look which edges come into it and which edge
  goes out, and this defines the local bouquet of a sphere.  Note
  however that this description involves the monad structure of
  $P$, since in reality we are taking the $P$-tree $T$ contained in the
  sphere and then contracting this tree to a single bouquet
  $t(T)$.  The map $\beta$ is defined similarly: to each layer,
  return the $P$-tree seen in that layer.  This is the $P$-tree contained in
  the output sphere of the layer but with the subtrees in the children
  contracted (here again we use the monad structure of $P$).  With
  $\alpha$ and $\beta$ described this way, it is clear that square (1)
  commutes: both ways around the square amount to taking the bouquet
  around the output sphere of a given layer.
  
  To define $\gamma:\ov L \to \tr\upperdot(P)$, notice that the $P$-tree seen
  in a given layer has a node for each child sphere of the layer.  So given
  a layer with a marked child, return the $P$-tree seen in this layer (as
  in the definition of $\beta$), with the node marked that corresponds to
  the child.  Now (2) is commutative and cartesian by construction.
  
  Finally, both ways around the square (3) amount to returning the
  bouquet of the marked child, which is the same as the local bouquet
  of the node in the tree-with-marked-node corresponding to the 
  layer-with-marked-child.
    
  {\em From $P^+$-tree to $P$-constellation.}   
  A $P^+$-tree $M$ is viewed as a recipe for how to glue
  small $P$-trees together to a big $P$-tree, the small $P$-trees being those 
  that decorate the nodes of $M$.  We refer to $M$ as the {\em composition tree}.
  In the end the gluing loci 
  will sit as spheres in the resulting big $P$-tree.
  
  We start with the special case where the $P^+$-tree $M$ is the
  unit tree \inlineDotlessTree\nolinebreak, i.e., a single edge decorated by some
  bouquet $b\in B$.  We need a $P$-constellation whose nesting
  corresponds to a unit tree.  Hence this constellation has no spheres,
  and thus has just a single dot, so it amounts to giving a one-dot
  $P$-tree.  Obviously we just take $b$ itself, considered as a $P$-tree
  via the unit map for the monad.
  
  If the composition tree $M$ has just one dot $n$, this dot is decorated
  by a $P$-tree $T$ (of a certain type).  We need to provide a sphere 
  nesting with just one sphere, and we just take $T$ with a sphere 
  around it.
  
  If the composition tree $M$ has more than one dot, then it has inner
  edges, and each inner edge $a$, say from node $c$ down to node $r$
  represents a substitution: the $P$-tree $T_r$ decorating $r$ has a node
  for each input edge of $r$; by the compatibility condition, the node
  corresponding to edge $a$ is decorated $A=t(T_c)$, the output type
  of $T_c$.  Hence it makes sense to substitute $T_c$ into that node 
  of $T_r$,
  cf.~(\ref{BD-subst}).  We should perform the substitutions
  corresponding to all the inner edges of $M$.  By associativity of
  the substitution law, we can make the substitutions
  edge by edge in any order.
  
  Hence it is enough to explain what happens for a composition tree 
  with a single inner edge, i.e., a two-dot tree.
    Suppose the composition tree looks like this:
  \begin{equation}\label{recipe}\begin{texdraw}
    \linewd 0.5 \scriptsize
    \bsegment
        \htext (0 -18) {\normalsize $M$}

      \move (0 0) \lvec (0 18) \onedot
      \htext (-5 15) {$r$}
      \htext (-8 24) {$a$}
      \move (0 18) \lvec (-6 36) \onedot
		    \lvec (-25 42) 
      \move (-6 36) \lvec (-19 52) 
      \move (-6 36) \lvec (-10 57) 
      \move (-6 36) \lvec (0 55) 
      \move (-6 36) \lvec (8 53) 
      \move (0 18) \lvec (16 48) 
      \move (0 18) \lvec (24 43) 
      \move (0 18) \lvec (28 35) 
      \move (0 18) \lvec (31 28) 
      \htext (-12 33) {$c$}
    \esegment
  \end{texdraw}\end{equation}
  where node $c$ is decorated by the $P$-tree $T_c$ of output type $A
  \in B$, while node $r$ is decorated by the $P$-tree $T_r$ one of
  whose nodes $f$ is decorated by $A\in B$.  Now the substitution goes
  like this (cf.~(\ref{BD-subst})):
  
    \begin{equation}\label{subst-node}\begin{texdraw}
  \scriptsize
    \move (0 0)
    \bsegment 
    \htext (0 -20) {\normalsize $T_c$}
    \move (0 0) \lvec (0 20) \onedot
    \lvec (-15 30) \onedot
    \move (0 20) \lvec (-17 65) \htext (-19 70) {$x$}
    \move (0 20) \lvec (-5 70) \htext (-6 75) {$y$}
    \move (0 20) \lvec (14 35) \onedot
    \lvec (9 50) \onedot \lvec (9 70) \htext (9 75) {$z$}
    \move (14 35) \lvec (23 47) \onedot
    \lpatt (2 2)
    \move (5 35)  \freeEllipsis{29}{19}{35}
    \lpatt (2 2)
    \move (20 20) \clvec (60 -15)(90 -15)(110 15)
    \lpatt ()
    \arrowheadtype t:F
    \arrowheadsize l:5 w:3
    \avec (112.5 18.5)
    \esegment

    \move ( 130 0)
    \bsegment 
    \htext (0 -20) {\normalsize $T_r$}
    \move (0 0) \lvec (0 17) \onedot
    \lvec (-10 30) \onedot  
    \lvec (-30 35) \htext (-35 36) {$x$}
    \htext (-14 25) {$f$}
    \move (-10 30) \lvec (-20 55) \onedot \htext (-21 47) {$y$}
  \move (-10 30) \lvec (-5 50) \onedot \lvec (-10 70)
  \move (-5 50) \lvec (0 70)
  \move (0 17) \lvec (15 60) 
  \move (0 17) \lvec (22 35) \onedot 
  \htext (-3 43) {$z$}
  \lpatt (2 2) \move (-12 29) \lcir r:11
    \esegment
    
\htext (240 15) {\normalsize  resulting in }

\move (340 0) 
   \scriptsize
   \bsegment 
   \move (0 0) \lvec (0 17) \onedot
   \lvec (-15 35) \onedot  
   \lvec (-35 70) \htext (-38 73) {$x$}
   \move (-15 35) \lvec (-30 40) \onedot

   \move (-15 35) \lvec (-25 85) \onedot \htext (-26 71) {$y$}
 \move (-15 35) 
 \lvec (-5 45) \onedot 
 \lvec (-5 60) \onedot \lvec (-5 80)  \onedot \lvec ( -10 95) 
 \move (-5 80) \lvec (0 95)

 \htext (-8 72) {$z$}
 \move (-5 45) \lvec (5 55) \onedot
 \move (0 17) \lvec (25 60) 
 \move (0 17) \lvec (22 35) \onedot 
 \lpatt (2 2)    \move (-12 47)  \freeEllipsis{26}{17}{32}

   \esegment

   \move (390 0)
 \end{texdraw}
 \end{equation}
 This $P$-tree is the underlying $P$-tree of the constellation we are
 constructing.  There should be two spheres: one outer sphere
 (corresponding to the root edge of $M$) for which there is no choice,
 and one inner sphere corresponding to the inner edge in $M$.  This
 inner sphere has to be precisely the scar of the surgery.  (The remaining 
 edges of $M$ are leaves and correspond to dots in the constellation we are 
 constructing.)
 
 If the composition tree has more inner edges, each corresponding 
 substitution will produce a sphere in the final tree, and 
 clearly the nesting resulting from all the substitutions
 will correspond to the composition tree 
 as required.
 
  (A short remark concerning two degenerate cases: If $T_c$ is the unit tree
  \inlineDotlessTree decorated by $b\in B$, then its output type is the bouquet $b= $ 
  \inlineOnedotTree, sitting as dot
  $f$ in $T_r$.  The effect of the substitution in this case is simply to erase
  the dot $f$, leaving a null-sphere as scar.  If $T_c$ is a one-dot tree, then
  we are substituting a single dot into a another dot of the same type, and the
  resulting tree is unchanged, but a sphere is placed around this dot, as scar
  of the operation.  The fact that the underlying tree stays the same just says
  that one-dot trees are the units for the substitution law.)
  
  It is clear from the construction that we similarly get a bijection
  $\tr\upperdot(P^+)=\const\uppercirc(P)$ compatible with the `source' map
  and the forgetful map as in (\ref{tr-const-diagram}).  Commutativity of the 
  right-hand triangle in (\ref{tr-const-diagram}) is clear 
  from the explicit description of the `target' map given in \ref{slice-twice}.
\end{dem}

To appreciate this result, note that a $P^+$-tree is a complicated
structure: it is a whole collection of $P$-trees (the decorations)
satisfying a complicated set of compatibility conditions.  The theorem
shows that all these data can be encoded in a single $P$-constellation,
where there are no compatibility conditions to check!

\bigskip

The theorem has the following interesting corollary:
\begin{cor}
  For any polynomial monad $P$, any abstract tree admits a $P^+$-decoration.
\end{cor}
In contrast, it is not true that any tree admits a decoration by a monad not of
the form $P^+$.  For example, only linear trees can be decorated by the trivial
monad.

\begin{dem*}{Proof of the corollary.}
  By the theorem, a $P^+$-decoration of a tree is the same thing as
  a $P$-constellation.  But every abstract nesting can appear as 
  underlying nesting of a constellation.  In fact for any 
  $P$-tree, you can draw arbitrary nestings.
\end{dem*}

\subsection{The polynomial monads of opetopes}

We shall generate all the opetopes iteratively, starting from the
identity monad on $\Set$.

\begin{blanko}{The opetope monads and the opetopes.}\label{BD-def}
  Let $P^0$ denote the identity monad on $\Set$, 
  \polyFunct{1}{1}{1}{1 \ .}{}{}{}{P^0} 
  Let $P^k$ denote the $k$th iterated Baez-Dolan construction on
  $P^0$.  By definition, the set of {\em $k$-dimensional opetopes}
  $\fat Z^k$ is the set of types for $P^k$, or equivalently, for
  $k\geq 1$, the set of operations for $P^{k-1}$, or for $k\geq 2$,
  the set of (isomorphism classes of) $P^{k-2}$-trees.  Finally define
  $\ov{\fat Z}{}^{k+1}$ to be the set appearing in the polynomial representation 
  of $P^k$ like this:
  \polyFunct{\fat Z^k}{\ov{\fat Z}{}^{k+1}}{\fat Z^{k+1}}{\fat Z^k \ .}{s}{p}{t}{P^k} 
  We define the {\em target} of an opetope $Z\in \fat Z^{k+1}$ to be the 
  $k$-opetope $t(Z)$,
  and we define the {\em sources} of $Z\in \fat Z^{k+1}$ to be the $k$-opetopes
  $s(F)$ where $F$ runs through the fibre $p^{-1}(Z)$.
  
  (Sources and targets are perhaps easiest understood in terms of trees: an
  $(k+1)$-opetope $Z$ is a $P^{k-1}$-tree: this means its nodes are decorated by
  $k$-opetopes (the operations for $P^{k-1})$.  These are the sources of $Z$.
  The target of $Z$ is obtained by contracting each inner edge of the tree,
  correspondingly substituting the decorating $k$-opetopes into each other.
  We shall explain this in Section~\ref{Sec:calculus}.)
\end{blanko}

Before establishing the general result reconciling this definition of opetope
with the elementary combinatorial definition of \ref{def}, let us work
out this comparison in low dimensions.

\begin{blanko}{Basis for the construction.}
  According to the definition, $\fat Z^0$ and $\fat Z^1$ are both the singleton
  set, in agreement with \ref{def}.  We write $\fat Z^0 \df \{ \inlineDotlessTree \}$ and
  $\fat Z^1 \df \{\inlineOnedotTree\}$, to
  conform with the standard graphical interpretation (cf.~\ref{graphical})
  of $P^0$:
      \begin{center}\begin{texdraw}
    \footnotesize
    
    \move (0 25)
    \move (-55 8) 
    \bsegment 
    \move (0 -10) \lvec (0 10) \move (0 0) \onedot 
      \htext (-3 9) {\normalsize *}
    \htext (-12 0) {$\bigleftbrace{15}$}
    \htext (12 0) {$\bigrightbrace{15}$}
    \esegment
	  
    \move (55 8)
    \bsegment 
    \move (0 -10) \lvec (0 10) \move (0 0) \onedot 
    \htext (-11 0) {$\bigleftbrace{15}$}
    \htext (12 0) {$\bigrightbrace{15}$}
    \esegment

    \move (-105 -55) 
    \bsegment 
    \move (0 -10) \lvec (0 10) \move (0 0) 
    \htext (-11 0) {$\bigleftbrace{15}$}
    \htext (12 0) {$\bigrightbrace{15}$}
    \esegment
    
    \move (105 -55) 
    \bsegment 
    \move (0 -10) \lvec (0 10) \move (0 0) 
    \htext (-11 0) {$\bigleftbrace{15}$}
    \htext (12 0) {$\bigrightbrace{15}$}
    \esegment

    \arrowheadtype t:F
    \arrowheadsize l:5 w:3

    \move (-18 10) \avec (18 10)
    \move (-75 -12) \avec (-90 -35) 
    \move (75 -12) \avec (90 -35)
    
    \htext (0 -60) {\normalsize $P^0=\Id$}
  \end{texdraw}\end{center}
\end{blanko}

\begin{blanko}{First iteration of the Baez-Dolan construction.}\label{BD-on-trivial}
  Applying the Baez-Dolan construction to $P^0$ we get the polynomial monad
  $P^1: \Set \to \Set$, which is nothing but the free-monoid monad $X\mapsto 
  \sum_{n\in \N} X^n$.  Hence $\fat Z^2 = \N$, in agreement with \ref{def}.
  In graphical terms, $\fat Z^2$ is the set of (isomorphism classes of)
  $P^0$-trees, i.e.~linear trees, and the picture is:
    \begin{center}\begin{texdraw}
    \footnotesize
    \move (-65 8) 
    \bsegment
    \move (0 0) \lvec (0 50)
      \move (0 10) \onedot
      \move (0 20) \onedot
      \move (0 30) \onedot
      \move (0 40) \onedot
      \htext (-6 30) {\normalsize *}
    \esegment
    \htext (-90 30) {$\bigleftbrace{25}$}
    \htext (-40 30) {$\bigrightbrace{25}$}
	  
    \move (65 8)
    \bsegment
    \move (0 0) \lvec (0 50)
      \move (0 10) \onedot
      \move (0 20) \onedot
      \move (0 30) \onedot
      \move (0 40) \onedot
    \esegment
    \htext (40 30) {$\bigleftbrace{25}$}
    \htext (90 30) {$\bigrightbrace{25}$}

    \move (-120 -45) 
    \bsegment 
    \move (0 -10) \lvec (0 10) \move (0 0) \onedot 
    \esegment
    \htext (-130 -45) {$\bigleftbrace{15}$}
    \htext (-110 -45) {$\bigrightbrace{15}$}
    \move (120 -45) 
    \bsegment 
    \move (0 -10) \lvec (0 10) \move (0 0) \onedot 
    \esegment
    \htext (130 -45) {$\bigrightbrace{15}$}
    \htext (110 -45) {$\bigleftbrace{15}$}

    \arrowheadtype t:F
    \arrowheadsize l:5 w:3

    \move (-15 30) \avec (15 30) 
    \move (-100 0) \avec (-115 -23) 
    \move (100 0) \avec (115 -23)
    
        \htext (0 -50) {\normalsize $P^1$}

  \end{texdraw}\end{center}
  Note that $\fat Z^2$ is not yet the set of $P$-constellations for any $P$.
\end{blanko}

\begin{blanko}{Second iteration of the BD construction.}
  Performing the Baez-Dolan construction a second time defines $P^2$.
  By Theorem~\ref{Thm:slicetwice}, this is about setting spheres in
  the trees we have got, which are the linear trees.  So $P^2$ looks
  like this:
   \begin{center}\begin{texdraw}
     \footnotesize
     \move (-65 28) 
     \bsegment
       \linewd 0.5

     \move (0 -35)
     \lvec (0 35) 
     \move (0 -20) \onedot
     \move (0 8) \onedot
     \move (0 20) \onedot
     \htext (-8 9) {\normalsize *}

     \linewd 0.3

     \move (0 15) \lellip rx:10 ry:14
     \move (0 2) \lellip rx:16 ry:30
     \move (0 -20)
       \lcir r:5
     \move (0 -7)
       \lcir r:4
     \move (0 20)
       \lcir r:5

     \esegment
     \htext (-95 30) {$\bigleftbrace{32}$}
     \htext (-35 30) {$\bigrightbrace{32}$}
	   
     \move (65 28)
     \bsegment
       \linewd 0.5

     \move (0 -35)
     \lvec (0 35) 
     \move (0 -20) \onedot
     \move (0 8) \onedot
     \move (0 20) \onedot
       
     \linewd 0.3

     \move (0 15) \lellip rx:10 ry:14
     \move (0 2) \lellip rx:16 ry:30
     \move (0 -20)
       \lcir r:5
     \move (0 -7)
       \lcir r:4
     \move (0 20)
       \lcir r:5

     \esegment
     \htext (35 30) {$\bigleftbrace{32}$}
     \htext (95 30) {$\bigrightbrace{32}$}

     \move (-120 -50) 
     \bsegment 
     \move (0 -15) \lvec (0 15) 
     \move (0 -7) \onedot \move (0 7) \onedot 
     \esegment
     \htext (-135 -50) {$\bigleftbrace{20}$}
     \htext (-105 -50) {$\bigrightbrace{20}$}
     \move (120 -50) 
     \bsegment 
     \move (0 -17) \lvec (0 17) 
     \move (0 0) \onedot \move (0 10) \onedot \move (0 -10) \onedot 
     \esegment
     \htext (135 -50) {$\bigrightbrace{20}$}
     \htext (105 -50) {$\bigleftbrace{20}$}

     \arrowheadtype t:F
     \arrowheadsize l:5 w:3

     \move (-15 30) \avec (15 30) 
     \move (-100 0) \avec (-115 -23) 
     \move (100 0) \avec (115 -23)
     
             \htext (0 -60) {\normalsize $P^2$}

   \end{texdraw}\end{center}
   So $\fat Z^3 = \const(P^0)$ is the set of (isomorphism classes of)
   constellations whose underlying tree is linear.  This is also the set of 
   (isomorphism classes of) planar trees, in agreement with \ref{def}.
\end{blanko}

\begin{blanko}{Third iteration of the BD construction.}
  For the next iteration --- trees of trees of trees --- a new
  meta-device is needed, so we zoom: take
  the tree expression of the nesting and set spheres in it like
  in the previous step.  More precisely, by
  Theorem~\ref{Thm:slicetwice} the set $\fat Z^3$ (of constellations
  whose underlying tree is linear) is also the set of $P^1$-trees,
  i.e.~trees with a certain compatible decoration by linear trees,
  and we know that to specify such a tree is just to draw the tree
  corresponding to the nesting, with a specified bijection:
  all the decorations can then be read off this bijection.
  Applying now the Baez-Dolan construction a third time just amounts to
  freely drawing  spheres in these composition trees.
  Figure~(\ref{const-monad}) serves well as illustration of $P^3$,
  although it is not clear from the figure that the underlying tree is
  a $P_1$-tree --- but $P_1$-means planar tree.  In conclusion, the
  set of operations $\fat Z^4$ corresponds with the $4$-opetopes 
  defined in \ref{def} and explained in \ref{dim4}.
\end{blanko}

\begin{satz}\label{op=op}
  Let $\fat O^k$ denote the set of $k$-opetopes in the sense of
  Definition~\ref{def} (isomorphism classes of degree-$k$ zoom complexes with an
  initial condition).  We have for $k\geq 0$ natural bijections
  $$
  \fat O^k = \fat Z^k .
  $$
\end{satz}

\begin{dem}
  We already established the claim for opetopes of dimension $0$, $1$, $2$, and 
  $3$, and proceed from here by induction.
  By Definition~\ref{BD-def} and Theorem~\ref{Thm:slicetwice} we have
  $\fat Z^{k+3} \df \operatorname{typ}(P^{k+3}) =
  \operatorname{op}(P^{k+2}) = \tr(P^{k+1}) = \const(P^k)$, for $k\geq 0$.
  So the claim is 
  $$
  \fat O^{k+3} = \const(P^k) = \tr(P^{k+1}) ,
  $$
  and in the induction step we shall need the auxiliary statement that the
  spheres in the top constellation of the $(k+3)$-opetope correspond to the
  spheres in the
  $P^k$-constellation (and hence to the tree in the $P^{k+1}$-tree).

  For $k\geq 1$, suppose given a $P^k$-constellation.  That's a $P^k$-tree $M$
  with some spheres --- we forget the spheres for a short moment.  By induction,
  $M$ can be interpreted as a $(k+2)$-opetope $W$ (i.e.~a zoom complex of degree
  $k+2$), and by the auxiliary detail, the top constellation of $W$ has
  underlying nesting (composition tree) $M$.  Now put back the spheres on $M$ to
  form a zoom complex of degree $k+3$, i.e.~a $(k+3)$-opetope.  Conversely,
  given a $(k+3)$-opetope, let $M$ denote the underlying tree of the top
  constellation, and forget for a moment the spheres in $M$.  The other
  constellations in the zoom complex (i.e.~up to degree $k+2$) form a
  $(k+2)$-opetope $W$ with composition tree $M$.  By induction, $W$ can be
  interpreted as a $P^k$-tree, which by the auxiliary detail has underlying tree
  $M$.  That is, $M$ is a $P^k$-tree.  Putting back the spheres on $M$ makes it
  into a $P^k$-constellation.  In both directions of the argument, it is clear
  that spheres correspond to spheres as required in the auxiliary detail.
\end{dem}

\subsection{Comparison}

  There exist in the literature four variations of the notion of opetope, not
  only in formulation but also in content: the original definition of
  Baez-Dolan~\cite{Baez-Dolan:9702}, the multitopes of
  Hermida-Makkai-Power~\cite{Hermida-Makkai-Power:I}, the opetopes in terms of
  cartesian monads due to Leinster~\cite{Leinster:0305049}, and a modification
  of the Baez-Dolan notion due to Cheng~\cite{Cheng:0304277}.  The four notions
  have been compared by Cheng~\cite{Cheng:0304277}, \cite{Cheng:0304279}. 
  
  We shall establish rather easily that our notion coincides with Leinster's.
  Our description of Leinster's sequence of cartesian monads stresses that all
  these monads are polynomial, and exploits the graphical calculus for
  polynomial functors to provide the explicit combinatorial description that was
  previously lacking.
  
\begin{blanko}{The original Baez-Dolan construction.}
  Baez and Dolan~\cite{Baez-Dolan:9702} described the construction first for
  algebras for a symmetric operad,  then they applied it to symmetric operads
  by observing that symmetric operad are themselves algebras for some operad.
  This is why they had to use {\em symmetric} operads.
\end{blanko}
  
\begin{blanko}{Baez-Dolan construction and definition of opetopes, according
  to Leinster~\cite[7.1]{Leinster:0305049}.}\label{BD-Leinster} Let $\EE$ be a
  presheaf category, and let $T$ be a finitary cartesian monad on $\EE$.
  (Leinster's setup is slightly more general.)  Then there is a notion of
  $T$-operad: a $T$-operad is a monoid in the monoidal category $\EE/T1$ for a
  certain tensor product.  Leinster~\cite[Appendix D]{Leinster:0305049} shows
  that the forgetful functor from $T$-operads to $\EE/T1$ has a left adjoint,
  the free $T$-operad functor.  This adjunction generates a monad which by
  definition is $T^+$.  It is clear that $\EE/T1$ is again a presheaf category,
  and Leinster proves that $T^+$ is again a finitary cartesian monad, hence the
  construction can be iterated.
  
  Leinster now defines the opetopes by starting with the identity functor
  $T_0$ on $\Set$ letting $T_k$ denote the $k$th iterated Baez-Dolan
  construction, and defining the set of opetopes in dimension $k$ to be the set
  of types for $T_k$.
\end{blanko}

Our setup is a special case of Leinster's, where $\EE$ is a slice of $\Set$, and
$T$ is a polynomial monad.  Note that polynomial functors always preserve
pullbacks, and our assumption that the representing map $E \to B$ is finite
amounts to $T$ being finitary.

\begin{satz}\label{ours=leinsters}
  If $P$ is a polynomial monad, the explicit polynomial Baez-Dolan construction
  $P \mapsto P^+$ of \ref{BD} coincides with Leinster's version \ref{BD-Leinster}.
  In particular, the opetopes defined in \ref{def} and \ref{BD-def} coincide
  with Leinster's opetopes.
\end{satz}

For the proof, we first reformulate Leinster's construction and specialise it
to the polynomial case.

\begin{blanko}{Reformulation of Leinster's description.}
  The reformulation removes reference to operads
  and the tensor product of collections.  Let $P$ be a cartesian monad on a
  presheaf category $\EE$.  Then there is a natural equivalence of categories
  \begin{eqnarray}\label{slice-equiv}
  \Cart(\EE)/P & \isopil & \EE/P1 \\
  {}[Q\Rightarrow P] & \mapsto & [Q1\to P1] , \notag
  \end{eqnarray}
  where $\Cart(\EE)$ denotes the category of cartesian endofunctors and
  cartesian natural transformations.  This equivalence follows readily from the
  fact that a cartesian natural transformation is completely determined by its
  value on a terminal object.  The category of endofunctors over $P$ has an
  obvious monoidal structure given by composition, relying on the monad
  structure of $P$: the composite of $Q \to P$ with $R\to P$ is $R\circ Q \to
  P\circ P \to P$ and the unit is $\Id\to P$.  One slick way to define the {\em
  tensor product of collections} (cf.~Kelly~\cite{Kelly:operads}) is to transport
  this canonical strict monoidal structure on $\Cart(\EE)/P$ along the
  equivalence~(\ref{slice-equiv}); operads are just monoids in the monoidal
  category of collections $\EE/P1$.  It follows that the free-$P$-operad monad on
  $\EE/P1$ is equivalent to the free-$P$-monad monad on $\Cart(\EE)/P$.  This
  monad in turn is just a matter of applying the free-monad construction on
  $\Cart$: on an object $Q$ this gives $Q\upperstar $, and if $Q$ is over $P$
  then $Q\upperstar$ is over $P\upperstar $ which in turn is over $P$ in virtue
  of the monad structure on $P$.  In conclusion, Leinster's Baez-Dolan
  construction on $P$ consists is just the transportation along the equivalence
  (\ref{slice-equiv}) of the free-monad monad over $P$.
\end{blanko}

\begin{blanko}{Specialisation to the polynomial case.}
  \label{BD-formal}%
  Denote by $\Poly(I)$ the category whose objects are polynomial endofunctors on
  $\Set/I$ and whose arrows are the cartesian natural transformations.
  Suppose $P$ is a polynomial monad represented by
  $$
  I \leftarrow E \to B \to I .
  $$
  It is a basic fact \cite{Kock:NotesOnPolynomialFunctors} that any functor $Q$
  with a cartesian natural transformation to $P$ is polynomial again, so the
  equivalence (\ref{slice-equiv}) reads
  \begin{eqnarray}\label{poly-slice-eq}
    \Poly(I)/P & \isopil & \Set/B  \\
    {}[Q \Rightarrow P] & \longmapsto & [Q1\to P1=B] . \notag
  \end{eqnarray}
  The inverse equivalence takes an object $C \to B$ in $\Set/B$ to the object 
  $Q$ in
  $\Poly(I)/P$ given by the fibre square
  \begin{equation}\label{Q-from-C/B}
  \begin{diagram}[w=4.5ex,h=3.5ex,tight]
  &&E\times_B C & \rTo  & C  \\
  &\ldTo&\dTo  &    & \dTo  & \rdTo \\
  I & \lTo &E  & \rTo  & B & \rTo & I.
  \end{diagram}  
  \end{equation}
  Denote by $\kat{PolyMon}(I)$ the category of polynomial monads on 
  $\Set/I$, i.e.~the category of monoids in $\Poly(I)$.  
  The forgetful functor
  $\kat{PolyMon}(I)/P \to \Poly(I)/P$ has a left adjoint, the free
  $P$-monad functor, hence generating a monad $T_P : \Poly(I)/P \to
  \Poly(I)/P$, which we referred to above as the free-$P$-monad monad,
  and which is the BD construction on $P$ modulo 
  equivalence~(\ref{poly-slice-eq}).
  
\end{blanko}

\begin{dem*}{Proof of Theorem~\ref{ours=leinsters}.}
  In view of the preceding discussion, the claim of the theorem is that $T_P$
  and $P^+$ correspond to each other under the monoidal equivalence
  (\ref{poly-slice-eq}).  Here $P^+$ denotes the explicit Baez-Dolan 
  construction of \ref{BD}.

  We already computed the value of $P^+$ on an object $C\to B$ of $\Set/B$: the
  result is the set of $P$-trees with each node decorated by an element of $C$,
  compatibly with the arity map $C \to B$ (being a $P$-tree means in particular
  that each node already has a $B$-decoration; these decorations must match).
  We claim that this is the same thing as a $Q$-tree, where $Q$ corresponds to
  $C\to B$ under equivalence (\ref{poly-slice-eq}) as in
  diagram~(\ref{Q-from-C/B}).  Indeed, since the tree is already a $P$-tree, we
  already have $I$-decorations on edges, as well as bijections for each node
  between the input edges and the fibre $E_b$ over the decorating element $b\in
  B$.  But if $c\in C$ decorates this same node, then the cartesian square
  specifies a bijection between the fibre over $c$ and the fibre $E_b$ and hence
  also with the set of input edges.  So in conclusion, $P^+$ sends $C$ to the
  set of $Q$-trees.
  
  On the other hand, $T_P$ sends the corresponding polynomial
  functor $Q$ to the free monad on $Q$, with structure map to $P$ given by
  the monad structure on $P$.  Specifically, $T_P$ produces from $Q$ the
  polynomial monad given by
  \begin{diagram}[w=6ex,h=4.5ex,tight]
  \tr'(Q)\SEpbk & \rTo  & \tr(Q)   \\
  \dTo  &    & \dTo  \\
  \tr'(P) \SEpbk  & \rTo  & \tr(P)\\
  \dTo&&\dTo\\
  E & \rTo & B
  \end{diagram}
  so the two endofunctors agree on objects.  The same argument works for arrows,
  so the two endofunctors agree.
  
  To see that the monad structures agree, note that the set of operations for
  $P^+\circ P^+$ is the set of $P$-trees with nodes decorated by $P$-trees
  in such a way that the total bouquet of the decorating tree matches the
  local bouquet of the node it decorates.  The composition law $P^+\circ
  P^+ \Rightarrow P^+$ consists in substituting each tree into the node it
  decorates.  On the other hand, to describe the monad $T_P$ it is enough to
  look at the base sets, since each top set is determined as fibre product
  with $E$ over $B$.  In this optic, $T_P$ sends $B$ to $\tr(P)$, and
  $T_P\circ T_P$ sends $B$ to $\tr(P\upperstar)$, whose elements are
  (isomorphism classes of)
  $P$-trees with nodes decorated
  by $P$-trees, and edges decorated in $I$, subject to the usual
  compatibility conditions. Clearly the composition law $T_P\circ T_P
  \Rightarrow T_P$ corresponds precisely to the one we described for $P^+$.
  For both monads, the unit is described as associating to a bouquet the
  corresponding one-dot tree.
 
  In conclusion, the two constructions agree.
\end{dem*}

\section{Suspension and stable opetopes}
\label{Sec:stable}

We introduce the notion of suspension of opetopes, define stable opetopes,
and show that the accompanying monad is the least fixpoint for the Baez-Dolan 
construction (for pointed monads).

\begin{blanko}{Suspension.}\label{suspension}
  The {\em suspension} $S(X)$ of an $n$-opetope $X$ is the
  $(n+1)$-opetope  defined by setting
  \begin{eqnarray*}
    S(X)_0 & \df & \inlineTrivialOpetope  \\
    S(X)_{k+1} & \df & X_k  \quad \text{ for }  0\leq k \leq n.
  \end{eqnarray*}
  In other words, just prepend a new \inlineTrivialOpetope to the zoom 
  complex, raising the indices.

  The operations `source', `target', and `composition of opetopes' all 
  commute with suspension.  Indeed, these operations are defined on the top 
  constellations, and the repercussions down through the zoom complex
  can never reach the degree-$1$ term in the complex.
\end{blanko}

\begin{blanko}{Stable opetopes.}
  The suspension defines a map $S : \fat Z^n \to \fat Z^{n+1}$ for 
  each $n\geq 0$.
  Let $\fat Z^\infty$ denote the colimit of this sequence of maps,
  $$
  \fat Z^\infty = \bigcup_{n\geq 0} \fat Z^n .
  $$
  This is the set of all opetopes in all dimensions, where we 
  identify two opetopes if one is the suspension of the other. 
  The elements in $\fat Z^\infty$ are called {\em stable opetopes}. Note 
  that a stable opetope has a well-defined top constellation, and that
  therefore the notions of source, target, and composition make sense 
  for stable opetopes.

  Define $\ov{\fat Z}{}^\infty \df \cup_{n\geq 0} \ov{\fat 
  Z}{}^n$,
  the set of stable opetopes with a marked input facet.
  Now consider the {\em polynomial monad of stable opetopes}
  $$
  P^\infty:\Set/\fat Z^\infty \to \Set/\fat Z^\infty
  $$
  defined by the diagram \polyFunct{\fat Z^\infty}{\ov{\fat
  Z}{}^\infty}{\fat Z^\infty}{\fat Z^\infty}{s}{}{t}{} 
  As usual, $t$ returns the target, $s$ returns the source, and 
  $\ov{\fat Z}{}^\infty \to \fat Z^\infty$ is the forgetful map.  This polynomial
  functor is a least fixpoint for the pointed Baez-Dolan construction,
  as we shall now explain.
\end{blanko}

\begin{blanko}{The category of polynomial monads.}
    Let $\kat{PM}$ denote the category of all polynomial monads 
    \cite{Gambino-Kock}.  The 
  arrows in this category are diagrams
  \begin{equation}\label{alpha}
  \begin{diagram}[w=2.5ex,h=3.7ex,tight]
    &&E'\SEpbk&&\rTo&&B'\\
    &\ldTo&&&&&&\rdTo\\
    I'&&\dTo&&&&\dTo[midshaft]<\alpha&&I' \\
    \dTo&&E&&\rTo&&B&&\dTo\\
    &\ldTo&&&&&&\rdTo\\
    I&&&&&&&&I
  \end{diagram}
  \end{equation}
which respect the monad structure.  This is most easily expressed in 
the partial-composition viewpoint where it amounts to requiring that
these two squares commute:
\begin{diagram}[w=6ex,h=4.5ex,tight]
B' \times_{I'} E' & \rTo  & B' & \lTo & I'\\
\dTo  &    & \dTo &&\dTo  \\
B \times_I E & \rTo  & B & \lTo & I
\end{diagram}

  The suspension map $S: \fat Z^n \to \fat Z^{n+1}$ induces an arrow
  in $\kat{PM}$:
  $$
  S : P^n \to P^{n+1}
  $$
  In other words, there is a natural diagram
  \begin{diagram}[w=3ex,h=3.8ex,tight]
    &&\ov{\fat Z}{}^{n+1}\SEpbk&&\rTo&&\fat Z^{n+1}\\
    &\ldTo&&&&&&\rdTo\\
    \fat Z^{n}&&\dTo&&&&\dTo&&\fat Z^{n} \\
    \dTo&&\ov{\fat Z}{}^{n+2}&&\rTo&&\fat Z^{n+2}&&\dTo\\
    &\ldTo&&&&&&\rdTo\\
    \fat Z^{n+1}&&&&&&&&\fat Z^{n+1}
  \end{diagram}
  The middle square is cartesian because marking a sphere in the top
  constellation is independent of suspension.  It is a monad map since
  suspension commutes with partial composition.
\end{blanko}

\begin{prop}
  The Baez-Dolan construction is functorial: it defines a functor 
  $BD:\kat{PM} \to \kat{PM}$.
\end{prop}

\begin{dem}
  We have to explain what $BD$
  does on arrows (and then it will be clear that composition of arrows
  and identity arrows are respected).  The Baez-Dolan construction on
  $\alpha$ given in (\ref{alpha}) is:
  \begin{diagram}[w=2.5ex,h=3.7ex,tight]
    &&\tr\upperdot(P') &&\rTo&&\tr(P')\\
    &\ldTo&&&&&&\rdTo\\
    B'&&\dTo&&&&\dTo[midshaft]<{\alpha\upperstar }&&B' \\
    \dTo&&\tr\upperdot(P)&&\rTo&&\tr(P)&&\dTo\\
    &\ldTo&&&&&&\rdTo\\
    B&&&&&&&&B 
  \end{diagram}
  Here $\alpha\upperstar : \tr(P')\to\tr(P)$ is defined
  already on the level of the free-monad construction.  The right-hand 
  square commutes because $\alpha$ is a monad morphism.  The rest is 
  pure combinatorics, about setting marks in trees. Since $\alpha\upperstar $ 
  is defined `node-wise', there is also an evident map $\tr\upperdot(P')\to 
  \tr\upperdot(P)$ which makes the two other squares commute, and for 
  which the middle square is cartesian.  Finally one can check that $\alpha\upperstar $
  is a monad morphism:
  \begin{diagram}[w=7.5ex,h=4.5ex,tight]
  \tr(P') \times_{B'} \tr\upperdot(P') &\rTo & \tr(P')  & \lTo & B'\\
  \dTo  &    & \dTo &&\dTo  \\
  \tr(P) \times_B \tr\upperdot(P) & \rTo  & \tr(P) & \lTo & B
  \end{diagram}
  Again this is a purely combinatorial matter: the horizontal maps are defined in 
  terms of substituting trees into nodes of trees.  Since the two rows are 
  just two instances of this, but with different decorations, the diagram 
  commutes.
\end{dem}

\begin{blanko}{Pointed polynomial monads.}
  The Baez-Dolan functor has a rather boring least fixpoint: it is simply the 
initial polynomial monad $\emptyset \leftarrow \emptyset \to \emptyset 
\to \emptyset$.  We are more interested in the notion of pointed
polynomial monads and the pointed analogue of the Baez-Dolan functor.

By a {\em pointed polynomial monad} we understand a polynomial monad
equipped with a monad map from the trivial monad
\polyFunct{1}{1}{1}{1}{}{}{}{\Id} A morphism of pointed polynomial
monads is one that respects the map from $\Id$.  This defines a category
$\kat{PM}\lowerstar $.  If $i:\Id\to M$ is a pointed polynomial monad,
then $BD(M)$ is naturally pointed again, so the Baez-Dolan
construction defines also a functor $\kat{PM}\lowerstar \to
\kat{PM}\lowerstar $.  To see this, note that by functoriality we get
a map $BD(\Id) \stackrel{BD(i)}{\rTo} BD(M)$.  On the other hand we have
$\Id= P^0$, the polynomial monad of $0$-opetopes, and $BD(\Id) = P^1$, and
the suspension map provides $\Id \to BD(\Id)$.  (Note that $P^1 : \Set \to
\Set$ is the free-monoid monad.)

Now it follows readily from the standard Lambek iteration argument that
\end{blanko}

\begin{prop}
  The polynomial monad $P^\infty$ of stable opetopes is a least fixpoint 
  for the Baez-Dolan construction $BD:\kat{PM}\lowerstar \to \kat{PM}\lowerstar $.
\end{prop}

Indeed, $P^\infty$ can be characterised as the colimit of
$$
\Id  \rTo BD(\Id) \rTo BD^2(\Id) \rTo \dots
$$

\section{Calculus of opetopes --- example computations}

\label{Sec:calculus}

In this section we make explicit how to manipulate opetopes represented
as zoom complexes.  In particular we are concerned with calculating sources
and target of opetopes and the operation of gluing opetopes together.
A reader who has skipped Sections~\ref{Sec:poly} and \ref{Sec:BD} can take the 
following descriptions as definitions.

In this section, by {\em root dot} we mean the dot adjacent to the root
edge (if there are any dots).

\subsection{Faces}

We follow the polytope-inspired terminology for opetopes,
and call their input and output devices {\em facets} (i.e.~codimension-$1$ 
faces):

\begin{blanko}{Target.}
  The {\em target facet} of an $n$-opetope $X$ is the $(n-1)$-opetope
  obtained by omitting
  the top constellation $X_n$ and the last zoom in the zoom complex.
  The target is also called the 
  {\em output facet}.
\end{blanko}

\begin{blanko}{Sources.}
  Let $X$ be an $n$-opetope.  For each sphere $s$ in $X_n$, there is a
  {\em source facet} (or {\em input facet}), which is an
  $(n-1)$-opetope.  You can think of it as the part of the zoom
  complex you can see by looking only through the layer determined by
  $s$, i.e., the region in $X_n$ delimited on the outside by $s$
  itself and from the inside by the children of $s$.  

  So there are three steps in the computation of the source facet
  corresponding to $s$:
  \begin{itemize}
    \item[(i)] up in $X_n$, consider only the layer determined by $s$.
    In other words, restrict to the sphere $s$ and contract all
    spheres contained in $s$;
    
    \item[(ii)] perform certain corresponding operations on the
    spheres in $X_{n-1}$ and in all lower constellations, in order to 
    maintain the constellations in zoom relation;
  
    \item[(iii)] omit $X_n$.
  \end{itemize}
  In a moment we shall
  describe this in detail, but first it is convenient to introduce the
  notions of globs and drops:
\end{blanko}

\begin{blanko}{Globs.}
  An $n$-opetope whose top constellation $X_n$ has precisely one sphere
  is called a {\em glob}.  In this case, there is precisely one source
  facet, and this facet is isomorphic to the target facet.  For each
  $(n-1)$-opetope $F$ there is a unique $n$-glob whose target facet is $F$,
  obtained by drawing the tree corresponding to the nesting
  underlying $F_{n-1}$, and drawing a sphere around it all.  This is
  called the glob over $F$. In abstract terms, it is nothing but the 
  unit operation of type $F$, cf.~\ref{BD}.
  Hence the globs in dimension $n$ are in
  natural bijection with the $(n-1)$-opetopes, via the
  target map.  The term `glob' comes from the polytope-style of drawing
  opetopes: in dimension $2$ there is only one glob, which is pictured like 
  this:
\begin{equation}\label{glob}
\begin{texdraw}
  \move (0 25)
\move (0 0)
\clvec(15 30)(35 30)(50 0)
\move (0 0) \lvec(50 0)
\htext (25 10) {$\Downarrow$}
\end{texdraw}
\end{equation}
\end{blanko}

\begin{blanko}{Drops.}
  An opetope whose top constellation $X_n$ has no spheres is called a
  {\em drop}.  So a drop has no sources.  Since a constellation
  without spheres necessarily has a unique dot, $X_{n-1}$ has a unique
  sphere.  Hence the target of a drop is always a glob.  In
  particular the set of all $n$-drops is in bijection with
  the set of all $(n-2)$-opetopes, via the target map applied twice.
  Again the terminology comes from the polytope-style drawing of 
  opetopes, where in dimension $2$ one can draw the unique drop as
  \begin{center}
    \begin{texdraw}
    \htext (0 -20) {$\Downarrow$}
    \move (0.8 0) \clvec(-45 -45)(45 -45)(-0.8 0)
    \move (0.8 0) \fcir f:1 r:0.7
    \move (-0.8 0) \fcir f:1 r:0.7
    \move (0 -35)
    \end{texdraw}
  \end{center}
  Notice that also in dimension $3$ there is only one drop (since there is
  only one $1$-opetope): it is the $3$-opetope whose sole facet is 
  (\ref{glob}).
\end{blanko}

\begin{blanko}{Sphere operations.}
  The operations involved in computing sources can be described in
  terms of the following sphere operations on a constellation $X_i$:
  \begin{itemize}
    \item  Erase a sphere which is not the outer sphere.
  
    \item  Draw a new sphere around a dot or a sphere.
  
    \item  Contract a sphere to a dot.
  
    \item  Restrict to a sphere.
  \end{itemize}
  Each operation on $X_i$ implies certain other operations on
  $X_{i-1}$, ensuring that the resulting constellations are
  in zoom relation, and these operations in turn imply other
  operations on $X_{i-2}$, and so on.  (It is understood that the
  sequence of operations starts at the top constellation and 
  propagates
  downwards, so we will not have to worry about consequences on
  $X_{i+1}$ of an operation on $X_i$.)
\end{blanko}

\begin{blanko}{Erasing a sphere (not the outer sphere), or drawing a
  new sphere around a dot or a sphere.} These operations do not have
  any consequences in the constellation below.
\end{blanko}

\begin{blanko}{Contracting a sphere to a dot.}
  Let $s$ be a sphere in $X_i$, and let $T$ denote the tree it cuts.
  If there is at least one dot in $T$, then let $r$ denote the root
  dot of $T$.  Then we are contracting $s$ down to $r$.  In $X_{i-1}$
  we must erase the spheres corresponding to each non-root dot in $T$,
  and that's all.
  If there are no dots in $T$ ($T$ consists of just an edge), then
  we are contracting $s$ down to a new dot which we denote
  $s^\bullet$.  Since $T$ is just a single edge, the dot 
  $s^\bullet $ will have a unique child $c$ (either a dot or a leaf).
  In $X_{i-1}$ we have to draw a new sphere around the sphere or dot
  corresponding to $c$.
\end{blanko}

\begin{blanko}{Restricting to a sphere.}
  Let $s$ be a sphere in $X_i$.  Restricting to $s$ means erasing
  everything outside it.  The new root edge will be the root edge of
  the tree $T$ cut by $s$, and each leaf of $T$ will be labelled by
  the dot (or leaf) the edge was connecting to outside $s$.  For each
  dot $x$ that is descendant of $T$ but not in $T$ itself,
  contract the corresponding sphere $x^\circ$ in $X_{i-1}$.
  Finally, restrict to the sphere $r^\circ$ in $X_{i-1}$ corresponding 
  to the root dot $r$ of $T$.  (If $T$ contains no dot, i.e.~is just an 
  edge, then instead of a root dot it has a unique leaf $r$; in that case 
  we are restricting to
  the corresponding {\em dot} $r^\circ$ in $X_{i-1}$.)
\end{blanko}

\begin{eks}\label{exampleX}
  We will compute the sources of the following $5$-opetope:
\begin{center}\begin{texdraw}\scriptsize
    \bsegment 
    \linewd 0.5
      \move (0 0) \lvec (0 74)
      \move (0 37) \onedot \htext (5 37) {$1$}
    \linewd 0.3
      \move (0 37) \lcir r:10 \htext (14 37) {$4$}
      \move (0 37) \lcir r:20 \htext (24 37) {$3$}
      \move (0 37) \lcir r:30 \htext (34 37) {$2$}
    \esegment
    \htext (0 -20) {\normalsize $X_2$}
      \htext (50 -18) {$\zoom$}

     \htext (100 -20) {\normalsize $X_3$}
    \move (100 0)

    \bsegment 
    \linewd 0.5
      \move (0 0)\lvec (0 74)
    \linewd 0.3
      \move (0 37) \onedot \lcir r:10
      \move (0 21) \onedot
      \move (0 32) \lcir r:19
      \move (0 57) \onedot
      \move (0 37) \lcir r:30
      \htext (4 37) {$3$}
      \htext (4 21) {$2$}
      \htext (5 57) {$4$}
      \htext (12 30) {$7$}
      \htext (19 45) {$6$}
      \htext (25 60) {$5$}
      \htext (0 81) {$1$}
    \esegment
          \htext (160 -18) {$\zoom$}

    \move (220 0)

    \bsegment 
    \htext (0 -20) {\normalsize $X_4$}
    \linewd 0.5
      \move (0 -5)
      \lvec (0 20) \onedot
      \lvec (-50 70)
      \move (0 20)
      \lvec (15 40) \onedot
      \lvec (-20 95)
      \move (15 40) 
      \lvec (20 60)
      \onedot
      \lvec (45 75)
   \linewd 0.3
	\move (17.5 50) \lellip rx:11 ry:16
	\htext (-4 16) {$5$}
	\htext (20 40) {$6$}
	\htext (15 60) {$7$}
	\move (-1 18) \lcir r:9
	\move (10 40) \freeEllipsis{37}{26}{60}
	\move (0 46) \lcir r:44
	\move (-30 50) \lcir r:6
	\htext (-25 5) {$8$}
	\htext (-23 25) {$9$}
	\htext (-30 62) {$10$}
	\htext (-53 72) {$4$}
	\htext (-23 99) {$2$}
	\htext (48 77) {$3$}
	\htext (13 18) {$11$}
	\htext (33 50) {$12$}
  \esegment
  
      \htext (285 -18) {$\zoom$}

  \move (350 0)
  \bsegment 
  \htext (0 -20) {\normalsize $X_5$}
  \linewd 0.5 
  \move (0 -5) \lvec (0 20) \onedot
  \lvec (-15 30) \onedot
  \move (0 20) 
  \lvec (20 38) \onedot
  \lvec (28 52) \onedot
  \lvec (28 90)

  \move (20 38) \lvec (-20 60)
  \onedot \lvec (-30 85)
  \move (-20 60) \lvec (-10 96)
  
	\htext (-5 18) {$8$}
	\htext (-21 30) {$10$}
	\htext (25 37) {$9$}
	\htext (-27 58) {$12$}
	\htext (22 53) {$11$}
	\htext (28 95) {$5$}
	\htext (-32 89) {$6$}
	\htext (-9 101) {$7$}
	\linewd 0.3
		\move (-9 25) \freeEllipsis{12}{20}{54}
		\move (24 45) \lellip rx:12 ry:18
		\move ( -5 52) \lcir r:6
		\move (0 46) \lcir r:44

		\htext (-27 5) {$13$}
		\htext (14 14) {$14$}
		\htext (24 23) {$15$}
		\htext (-5 63) {$16$}
		
  \esegment

  \end{texdraw}
\end{center}
There are sources corresponding to the spheres $13$, $14$, $15$, and 
$16$; we will denote these source facets by $S13$, $S14$, $S15$, and 
$S16$.
\end{eks}


\begin{blanko}{Computation of source $S13$.}
  Step (i): contract
$14$, $15$, and $16$ in $X_5$:
\begin{center}\begin{texdraw}\scriptsize
  \bsegment
  \linewd 0.5 
    \move (0 0) \lvec (0 20) \onedot
    \move (0 20) \lvec (20 38) \onedot
    \lvec (28 88) \htext (28 93) {$5$}
    \move (20 38) \lvec (-20 60) \onedot 
    \lvec (-30 83) \htext (-32 87) {$6$}
    \move (-20 60) \lvec (-10 92) \htext (-9 97) {$7$}

      \htext (-7 21) {$14$}
      \htext (27 37) {$15$}
      \htext (-27 58) {$12$}
    \linewd 0.3
      \move ( -5 52) \onedot \htext (-5 60) {$16$}

      \move (0 46) \lcir r:38
      \htext (-27 11) {$13$}
      \htext (0 -20) {\text{layer `13'}}
  \esegment
\end{texdraw}\end{center}
Step (ii): perform the corresponding operations in
the lower constellations, according to the sphere operations rules.
This means deleting spheres $10$ and $11$, and drawing a new sphere
around sphere $12$ (corresponding to the contracted `empty' sphere
$16$).  Finally (iii), omit the top constellation. The end result is:
\begin{center}\begin{texdraw}\scriptsize
    \bsegment 
    \linewd 0.5
      \move (0 0) \lvec (0 74)
      \move (0 37) \onedot \htext (5 37) {$1$}
    \linewd 0.3
      \move (0 37) \lcir r:10 \htext (14 37) {$4$}
      \move (0 37) \lcir r:20 \htext (24 37) {$3$}
      \move (0 37) \lcir r:30 \htext (34 37) {$2$}
      \htext (0 -20) {\small $S13_2$}
    \esegment
      \htext (50 -18) {$\zoom$}

    \move (100 0)

    \bsegment 
    \linewd 0.5
      \move (0 0)\lvec (0 74)
    \linewd 0.3
      \move (0 37) \onedot \lcir r:10
      \move (0 21) \onedot
      \move (0 32) \lcir r:19
      \move (0 57) \onedot
      \move (0 37) \lcir r:30
      \htext (4 37) {$3$}
      \htext (4 21) {$2$}
      \htext (5 57) {$4$}
      \htext (12 30) {$7$}
      \htext (19 45) {$6$}
      \htext (25 60) {$5$}
      \htext (0 81) {$1$}
      \htext (0 -20) {\small $S13_3$}
    \esegment
          \htext (162.5 -18) {$\zoom$}

    \move (225 0)

    \bsegment 
    \linewd 0.5
      \move (0 -5)
      \lvec (0 20) \onedot
      \lvec (-50 70)
      \move (0 20)
      \lvec (15 40) \onedot
      \lvec (-20 95)
      \move (15 40) 
      \lvec (20 60)
      \onedot
      \lvec (45 75)
   \linewd 0.3
	\move (17.5 50) \lellip rx:11 ry:16
	\move (17.5 50) \lellip rx:16 ry:21
	\htext (-4 16) {$5$}
	\htext (20 40) {$6$}
	\htext (15 60) {$7$}
	\move (10 40) \freeEllipsis{37}{26}{60}
	\move (0 46) \lcir r:44
	\htext (-26 5) {$14$}
	\htext (-24 25) {$15$}
	\htext (-53 72) {$4$}
	\htext (-23 99) {$2$}
	\htext (48 77) {$3$}
	\htext (24 50) {$12$}
	\htext (-3 50) {$16$}
	\htext (0 -20) {\small $S13_4$}
  \esegment

  \end{texdraw}
\end{center}

\end{blanko}

\begin{blanko}{Computation of source $S14$.}
  Step (i): restrict to sphere $14$:
\begin{center}\begin{texdraw}\scriptsize
    \bsegment
    \linewd 0.5
      \move (0 0) \lvec (0 18) \onedot
      \lvec (-10 28) \onedot \htext (-10 35) {$10$}
      \move (0 18) \lvec (15 50) \htext (18 54) {$9$}
    \linewd 0.3
      \htext (4 14) {$8$}
      \move (0 26) \lcir r:20
      \htext (-16 6) {$14$}
  \esegment
\end{texdraw}\end{center}
Step (ii) amounts to contracting sphere $9$ in $X_4$, and hence
erasing sphere $6$ and $7$ down in $X_3$.  End result:
    
\begin{center}\begin{texdraw}\scriptsize
    \bsegment 
    \linewd 0.5
      \move (0 0) \lvec (0 74)
      \move (0 37) \onedot \htext (5 37) {$1$}
    \linewd 0.3
      \move (0 37) \lcir r:10 \htext (14 37) {$4$}
      \move (0 37) \lcir r:20 \htext (24 37) {$3$}
      \move (0 37) \lcir r:30 \htext (34 37) {$2$}
      \htext (0 -20) {\small $S14_2$}

    \esegment
        \htext (50 -18) {$\zoom$}

  \move (100 0)
 
  \bsegment 
  \linewd 0.5
    \move (0 0) \lvec (0 60)
  \linewd 0.3
    \move (0 41) \onedot
    \move (0 19) \onedot
    \move (0 30) \onedot
    \move (0 30) \lcir r:22
    \htext (5 30) {$3$}
    \htext (5 19) {$2$}
    \htext (5 41) {$4$}
    \htext (20 13) {$9$}
    \htext (0 65) {$1$}
    \htext (0 -20) {\small $S14_3$}

  \esegment
      \htext (150 -18) {$\zoom$}

  \move (200 0)
  
  \bsegment 
  \linewd 0.5
    \move (0 0) \lvec (0 18) \onedot
    \lvec (-25 43) \htext (-28 45) {$4$}
    \move (0 18) \lvec (0 58) \htext (0 63) {$2$}
    \move (0 18) \lvec (25 43) \htext (28 45) {$3$}
  \linewd 0.3
    \htext (-4 14) {$9$}
    \move (0 30) \lcir r:23
    \move (-10 28) \lcir r:6
    \htext (-10 38) {$10$}
    \htext (18 9) {$8$}
    \htext (0 -20) {\small $S14_4$}
  \esegment
\end{texdraw}\end{center}
\end{blanko}


\begin{blanko}{Computation of source $S15$.}
  Step (i): restrict to sphere $15$:
\begin{center}
  \begin{texdraw}
    \scriptsize
  \move (345 0)
    \bsegment
    \linewd 0.5
      \move (0 0)
      \lvec (0 18) \onedot
      \lvec (-10 50)      \htext (-12 54) {$12$}
\move (0 18)
      \lvec (8 32) \onedot
      \lvec (25 40) \htext (28 42) {$5$}
	\htext (-5 18) {$9$}
	\htext (6 38) {$11$}
	\linewd 0.3
	\move (0 26) \lcir r:20
	\htext (-16 6) {$15$}
  \esegment
\end{texdraw}\end{center}
    This implies (step (ii)) that in $X_4$ we have to restrict to sphere $9$ 
    and contract sphere $12$.  Down in $X_3$ this means erase sphere 
    $7$.  End result:
\begin{center}\begin{texdraw}\scriptsize
  \bsegment 
  \linewd 0.5
    \move (0 0) \lvec (0 74)
    \move (0 37) \onedot \htext (5 37) {$1$}
  \linewd 0.3
    \move (0 37) \lcir r:10 \htext (14 37) {$4$}
    \move (0 37) \lcir r:20 \htext (24 37) {$3$}
    \move (0 37) \lcir r:30 \htext (34 37) {$2$}
    \htext (0 -20) {\small $S15_2$}

  \esegment
      \htext (50 -18) {$\zoom$}

  \move (100 0)
   
  \bsegment 
  \linewd 0.5
    \move (0 0) \lvec (0 65)  
  \linewd 0.3
    \move (0 34) \onedot
    \move (0 24) \onedot
    \move (0 29)
    \lcir r:14
    \move (0 50) \onedot
    \move (0 34)
    \lcir r:25
    \htext (5 34) {$3$}
    \htext (5 24) {$2$}
    \htext (5 50) {$4$}
    \htext (18 35) {$12$}
    \htext (22 56) {$5$}
    \htext (0 70) {$1$}
    \htext (0 -20) {\small $S15_3$}

  \esegment
      \htext (150 -18) {$\zoom$}

    \move (200 0)
    
  \bsegment 
  \linewd 0.5
    \move (0 0)
    \lvec (0 23) \onedot
    \lvec (-25 48) \htext (-28 50) {$4$}
    \move (0 23)
    \lvec (8 40) \onedot
    \lvec (-10 62) \htext (-13 65) {$2$}
    \move (8 40)
    \lvec (25 50) \htext (28 51) {$3$}
    \htext (15 37) {$12$}
   \linewd 0.3
    \htext (-4 19) {$5$}
    \move (0 33) \lcir r:24
    \move (-2 22) \lcir r:9
    \htext (12 23) {$11$}
    \htext (18 12) {$9$}
    \htext (0 -20) {\small $S15_4$}

  \esegment

\end{texdraw}\end{center}
\end{blanko}

\begin{blanko}{Computation of source $S16$.}
  Step (i): restrict to sphere $16$:
\begin{center}\begin{texdraw}\scriptsize
  \bsegment
  \linewd 0.5
    \move (0 0) \lvec (0 30) \htext (0 35) {$12$}
    \move (0 16) \linewd 0.3 \lcir r:10 \htext (-14 6) {$16$}
  \esegment
\end{texdraw}\end{center}
Step (ii): the root of this subtree is the leaf $12$, so down in $X_4$ we have 
contract sphere $12$ and then restrict to the resulting dot $12$.
The contraction has the consequences in $X_3$ of erasing sphere $7$
(and we rename sphere $6$ to $12$). Restricting
to dot $12$ in $X_4$ means restricting to sphere $12$ in $X_3$.
Since dot $4$ is a descendant which is not inside sphere $12$, we have 
to contract sphere $4$ in $X_2$.  End result:
\begin{center}\begin{texdraw}\scriptsize
  \bsegment 
  \linewd 0.5
    \move (0 0) \lvec (0 50) 
    \move (0 25) \onedot \htext (5 25) {$4$}
  \linewd 0.3
    \move (0 25) \lcir r:10 \htext (14 25) {$3$}
    \move (0 25) \lcir r:20 \htext (24 25) {$2$}
    \htext (0 -20) {\small $S16_2$}

  \esegment
        \htext (45 -18) {$\zoom$}

  \move (90 0)
   
  \bsegment 
  \linewd 0.5
    \move (  0 5 ) \lvec ( 0 43)
  \linewd 0.3
    \move (0 29) \onedot
    \move (0 19) \onedot
    \htext (5 29) {$3$}
    \htext (5 19) {$2$}
    \move (0 24) \lcir r:14
    \htext (17 12) {$12$}
    \htext (0 48) {$4$}
    \htext (0 -20) {\small $S16_3$}
  \esegment
          \htext (130 -18) {$\zoom$}

  \move (170 0)

  \bsegment 
  \linewd 0.5
    \move (0 5)
    \lvec (0 20) \onedot
    \lvec (-9 38) \htext (-11 42) {$2$}
    \move (0 20)
    \lvec (9 38) \htext (11 42) {$3$}
    \htext (7 19) {$12$}
    \htext (0 -20) {\small $S16_4$}
  \esegment
\end{texdraw}\end{center}
\end{blanko}

\subsection{Composition tree and gluing}

\begin{blanko}{Composition tree.}
  The {\em composition tree} of an opetope is simply the tree
  corresponding to the nesting of the top constellation (with a
  specified correspondence).  It concisely expresses the incidence
  relations among the codimension-$1$ faces, and how these faces are
  attached to each other along codimension-$2$ faces.  We denote the 
  composition tree of $X$ by $\ct(X)$.

  In the composition tree $\ct(X)$, each dot $s$ corresponds to an input
  facet $S$ (codimension-$1$ face).  The last codimension-$1$ face of $X$,
  its target facet, is represented in the composition tree as the `total
  bouquet', i.e.~the bouquet obtained by contracting all inner edges
  (or setting a dot in the unit tree, if $X$ is a drop).

  The edges in $\ct(X)$ correspond to the codimension-$2$ faces of $X$:
  There is an incoming edge of dot $s$ for each input facet of $S$, and the
  output edge of $s$ represents the output facet of $S$.  In other words,
  an edge linking a dot $s$ to its parent dot $p$ represents the
  codimension-$2$ face along which $S$ is attached to $P$ (the face
  corresponding to $p$): this codimension-$2$ face is the target facet of
  $S$ and one specific source facet of $P$.  This source is easily
  determined: $p$ is a sphere in $X_n$ and $s$ is another sphere
  immediately contained in $p$.  When computing $P$ we contract the sphere
  $s$ to a dot, hence it becomes a sphere in $P_{n-1}$, and so represents a
  source facet of $P$.

  The leaves of $\ct(X)$ correspond to the dots in the top constellation,
  which in turn correspond to the spheres in $X_{n-1}$.  These are
  precisely the input facets of the target of $X$.  By the preceding
  discussion, each of these codimension-$2$ faces is also the source facet
  of exactly one source facet of $X$, namely the facet $S$ corresponding to
  the parent dot $s$ of the leaf.

  If there is a dot in $\ct(X)$ (i.e.~$X$ is not a drop), then the root dot
  determines a {\em bottom source}, characterised also as the source facet
  having the same target as the target of $X$ (corresponding to the output
  edge of $\ct(X)$).

  In summary we see that, except if $X$ is a drop, every codimension-$2$
  face of $X$ occurs exactly twice as a facet of a facet.  In fact, more
  generally, if $V$ is a codimension-$(k+2)$ face of an opetope $X$, and
  $F$ is a codimension-$k$ face of $X$ containing $V$, then the number of
  codimension-$(k+1)$ faces $E$ such that $V \subset E \subset F$ is either
  $1$, or $2$.  It is $1$ if and only if $F$ is a drop (in which
  case it is the drop on $E$ (which in turn is a glob on $V$)).  
\end{blanko}

\begin{blanko}{Example (continued from \ref{exampleX}).}
  For the opetope $X$ of the example above, the composition tree is
  \begin{center}\begin{texdraw}
    \linewd 0.5 \scriptsize
    \bsegment
      \move (0 0) \lvec (0 18) \onedot
      \lvec (-25 48) \htext (-30 50) {$12$}
      \htext (-6 15) {$13$}
      \move (0 18) \lvec (-6 36) \onedot
      \lvec (-15 58) \htext (-19 62) {$10$}
      \move (-6 36) \lvec (-1 62) \htext (0 66) {$8$}
      \move (0 18) \lvec (16 43) \onedot
      \lvec ( 15 62) \htext ( 15 67) {$11$}
      \move (16 43) \lvec (30 58) \htext (34 62) {$9$}
      \move (0 18) \lvec (30 30) \onedot \htext ( 37 30) {$16$}  
      \htext (1 36) {$14$}
      \htext (22 40) {$15$}
      \htext (0 -20) {\normalsize $\ct(X)$}
    \esegment
  \end{texdraw}\end{center}
  We see that $S13$ (corresponding to dot $13$) has four input facets
  (corresponding to the four input edges of dot $13$): the first one (leaf
  $12$) is left vacant, its three other input facets serve as gluing locus
  for the output facets of $S14$, $S15$, and $S16$.  In turn, $S14$ and
  $S15$ each has two input facets (which are not in use for gluing), while
  $S16$ has no input facets (i.e., $S16$ is a drop).  Note that the root
  edge represents the output facet of $S13$.
\end{blanko}

\begin{blanko}{Gluing and filling.}
  As explained in the proof of Theorem~\ref{Thm:slicetwice},
a decorated composition tree serves as a recipe for gluing together
$n$-dimensional opetopes $S_i$, producing one
big $n$-dimensional opetope $T$, and finally filling the whole
thing with an $n$-dimensional opetope $X$ in such a way that the
original opetopes $S_i$ become the input facets of $X$, and $T$ becomes the
output facet.

The first part consists in producing the `composite' opetope
$T$ from the $S_i$ according to the recipe specified by the
composition tree.  This can be done in steps: it is enough to explain
what happens when the composition tree has a single inner edge, i.e.,
a simple gluing.  The second part (\ref{filler}) consists in 
constructing the filling $(n+1)$-opetope $X$.
\end{blanko}

\begin{blanko}{Gluing.}\label{gluing}
  Given an $n$-opetope $R$ with a specified source $F$, and another
  $n$-opetope $S$ with target $F$, then their composite $T$ is again
  an $n$-opetope, whose target is the target of $R$, and whose set of
  sources is 
  $$
  \text{sources}(S) \cup \text{sources}(R) \shortsetminus
  \{F\}.
  $$
  The recipe composition tree looks something like this:
  \begin{equation}\begin{texdraw}
    \linewd 0.5 \scriptsize
    \bsegment
      \move (0 0) \lvec (0 18) \onedot
      \htext (-5 15) {$R$}
      \htext (-8 24) {$F$}
      \move (0 18) \lvec (-6 36) \onedot
		    \lvec (-25 42) 
      \move (-6 36) \lvec (-19 52) 
      \move (-6 36) \lvec (-10 57) 
      \move (-6 36) \lvec (0 55) 
      \move (-6 36) \lvec (8 53) 
      \move (0 18) \lvec (16 48) 
      \move (0 18) \lvec (24 43) 
      \move (0 18) \lvec (28 35) 
      \move (0 18) \lvec (31 28) 
      \htext (-12 33) {$S$}
    \esegment
  \end{texdraw}\end{equation}

  Every such situation arises as follows.  Write down an arbitrary
  $n$-opetope $R$ (but not a drop), pick one of its source facets, and
  write down this $(n-1)$-opetope $F$.  Next we need to provide an
  $n$-opetope $S$ having $F$ as its target.  By definition of the
  target map, $S$ is obtained from $F$ by drawing its composition tree
  and then drawing some arbitrary spheres in it.
\end{blanko}
 
\begin{blanko}{Example.}
  Let us illustrate the situation with an example.  Here
  is $S$:
  \begin{center}
  \begin{texdraw}
  \scriptsize

  \htext (-70 35) {\dots}
  \move (10 0)
  \bsegment 
  \htext (0 -20) {\normalsize $S_{n-1}$}

    \move (0 -5)
    \lvec (0 15) \onedot
    \lvec (-15 30) \onedot \lcir r:7
    \lvec (-6 50) \onedot
    \lvec (10 60) \onedot
    \lvec (20 85)
    \move (-6 50) \lvec (-20 80)
    \move (0 15) \lvec (45 50)
    \move (2 55) \linewd 0.3 \freeEllipsis{17}{10}{30}
    \htext (16 47) {$a$}
    \htext (-25 35) {$b$}
      \move (0 41) \lcir r:36
      \htext (32 15) {$c$}
      \htext ( 1 48) {$m$}

    \esegment
    
    \htext (75 -18){$\zoom$}

    \move (150 0)
    \bsegment 
    \htext (0 -20) {\normalsize $S_n$}
    \linewd 0.5
    \move (0 -5) \lvec (0 25) \onedot
    \linewd 0.3 \lcir r:10 \linewd 0.5
    \lvec (-25 45) \onedot
    \lvec (-50 50) \htext (-30 42) {$a$}
    \move (-25 45) \lvec (-20 95)
    \move (0 25) 
    \lvec (50 55)
    \htext (-5 23) {$c$}
    \move (20 37) \linewd 0.3 \lcir r:6 \linewd 0.5
    \move (0 25) \lvec (5 60) \onedot
    \lvec (10 95) \htext (0 60) {$b$}
    \move (-23 65) \linewd 0.3 \lcir r:6
	      \move (9 30) \freeEllipsis{27}{16}{30}
	      \move (0 46) \lcir r:42
	      \htext ( -55 50) {$m$}

    \esegment
    
    \move (300 0)
    \bsegment 
    \htext (0 -20) {\normalsize $\ct(S)$}
    \linewd 0.5
    \move (0 0) \lvec (0 20) \onedot
    \lvec (-15 30) \onedot
    \move (0 20) \lvec (-17 65) \htext (-19 70) {$a$}
    \move (0 20) \lvec (-5 70) \htext (-6 75) {$b$}
    \move (0 20) \lvec (14 35) \onedot
    \lvec (9 50) \onedot \lvec (9 70) \htext (9 75) {$c$}
    \move (14 35) \lvec (23 47) \onedot
    \esegment

    \end{texdraw}
  \end{center}
  And here comes $R$:

  \begin{center}
  \begin{texdraw}
  \scriptsize
  \htext (-70 40) {\dots}
  \move (10 0)

  \bsegment
  \htext (0 -20) {\normalsize $R_{n-1}$}

  \linewd 0.5
    \move (0 -5)
    \lvec (0 20) \onedot
    \lvec (-15 35) \onedot 
    \lvec (-6 55) \onedot
    \lvec (10 65) \onedot
    \lvec (20 105)
    \move (-6 55) \lvec (-20 100)
    \move (0 20) \lvec (45 55) \onedot 
    \lvec (45 95)
    \move (45 55) \lvec (65 75)
    \linewd 0.3 
    \move (2 60) \freeEllipsis{17}{10}{30}
    \htext (16 52) {$a$}
    \htext (7 36) {$y$}
      \move (0 46) \lcir r:36
      \htext (32 20) {$x$}
      \move (-6.5 27.5) \freeEllipsis{10}{17}{43}
      \move (4 49) \lellip rx:54 ry:48
      \htext ( 1 53) {$m$}

      \htext (-8 36) {$k$}
    \esegment
    
    \htext (85 -18){$\zoom$}

    \move (170 0)
    \bsegment 
    \htext (0 -20) {\normalsize $R_n$}
    \linewd 0.5
    \move (0 -5) \lvec (0 20) \onedot
    \lvec (0 45) \onedot \lvec (-25 55) \onedot
    \lvec ( -55 60) \move (-25 55) \lvec (-57 75)
    \htext ( -26 50) {$a$}
    \htext (5 45 ) {$x$}
    \move (0 20) \lvec (60 45)
    \linewd 0.3 
    \move ( 24 30) \lcir r:6 
    \linewd 0.5
    \move (0 45) \lvec (5 65) \onedot 
    \lvec (30 105) 
    \move (5 65) \lvec (-45 95)  
    \htext (-49 98) {$k$}
    \linewd 0.3 
    \move (-16 77) \lcir r:6 \htext (-25 72) {$b$}
	      \move (4 55) \freeEllipsis{18}{12}{80}
	      \htext (10 64) {$y$}
	      \move (-2 53) \lcir r:50
	      \move (-9 60) \freeEllipsis{35}{28}{-20}
  \htext (30 55) {$f$}
  \htext (-11 60) {$c$}
  \htext ( -60 60) {$m$}
    \esegment

    \move (330 0)
    \bsegment 
    \htext (0 -20) {\normalsize $\ct(R)$}
    \linewd 0.5
    \move (0 5) \lvec (0 22) \onedot
    \lvec (-10 35) \onedot
    \lvec (-30 40) \htext (-35 41) {$a$}
    \htext (-14 30) {$f$}
    \move (-10 35) \lvec (-20 55) \onedot \htext (-25 59) {$b$}
  \move (-10 35) \lvec (-5 50) \onedot \lvec (-10 70)
  \move (-5 50) \lvec (0 70)
  \move (0 22) \lvec (15 65) 
  \move (0 22) \lvec (22 40) \onedot 
  \htext (1 50) {$c$}
    \esegment

    \end{texdraw}
  \end{center}
%

\noindent
  Now $F$ is the target of $S$ and at the same time 
  the source of $R$ corresponding to 
  sphere $f$:
  \begin{center}
  \begin{texdraw}
  \scriptsize
  \htext (-80 30) {\dots}
  \move (0 0)

  \bsegment 
  \htext (0 -20) {\normalsize $F_{n-1}$}

  \linewd 0.5
    \move (0 -5)
    \lvec (0 15) \onedot
    \lvec (-15 30) \onedot \linewd 0.3 \lcir r:7 \linewd 0.5
    \lvec (-6 50) \onedot
    \lvec (10 60) \onedot
    \lvec (20 85)
    \move (-6 50) \lvec (-20 80)
    \move (0 15) \lvec (45 50)
    \linewd 0.3 
    \move (2 55) \freeEllipsis{17}{10}{30}
    \htext (16 47) {$a$}
    \htext (-25 35) {$b$}
      \move (0 41) \lcir r:36
      \htext (32 15) {$c$}
      \htext ( 1 48) {$m$}
    \esegment
    

    \move (150 0)
    \bsegment 
    \htext (0 -20) {\normalsize $\ct(F)$}
    \move (0 5) \lvec (0 25) \onedot
    \lvec (-15 35) \onedot
    \lvec (-35 45) \htext (-22 32) {$a$}
    \move (-15 35) \lvec (-20 55)
    \move (0 25) 
    \lvec (25 40)
    \htext (-5 23) {$c$}
    \move (0 25) \lvec (5 40) \onedot
    \lvec (10 55) \htext (0 40) {$b$}
    \esegment

   \end{texdraw}
   \end{center}

  We need to construct a new $n$-opetope $T$
  whose target is the same as the target of $R$.  This means that it
  differs from $R$ only in the top constellation, where the
  configuration of spheres is different.  The difference in sphere
  layout is expressed nicely in terms of the composition trees of $S$
  and $R$.  The recipe prescribes that we should glue $S$ onto the
  $F$-facet of $R$.  In terms of the composition trees of $S$ and $R$
  this means that we must substitute the whole tree $\ct(S)$
  into the node $f$ of $\ct(R)$.  Since the target of $S$ is
  $F$, this will again produce a valid decorated composition tree which
  will be $\ct(T)$.  In the current example, the situation is
  this:
   
  \begin{equation}\begin{texdraw}
  \scriptsize
    \move (0 0)
    \bsegment 
    \htext (0 -20) {\normalsize $\ct(S)$}
    \move (0 0) \lvec (0 20) \onedot
    \lvec (-15 30) \onedot
    \move (0 20) \lvec (-17 65) \htext (-19 70) {$a$}
    \move (0 20) \lvec (-5 70) \htext (-6 75) {$b$}
    \move (0 20) \lvec (14 35) \onedot
    \lvec (9 50) \onedot \lvec (9 70) \htext (9 75) {$c$}
    \move (14 35) \lvec (23 47) \onedot
    \lpatt (2 2)
    \move (5 35)  \freeEllipsis{29}{19}{35}
    \lpatt (2 2)
    \move (20 20) \clvec (60 -15)(90 -15)(110 15)
    \lpatt ()
    \arrowheadtype t:F
    \arrowheadsize l:5 w:3
    \avec (112.5 18.5)
    \esegment

    \move ( 130 0)
    \bsegment 
    \htext (0 -20) {\normalsize $\ct(R)$}
    \linewd 0.5
    \move (0 0) \lvec (0 17) \onedot
    \lvec (-10 30) \onedot  
    \lvec (-30 35) \htext (-35 36) {$a$}
    \htext (-14 25) {$f$}
    \move (-10 30) \lvec (-20 50) \onedot \htext (-25 54) {$b$}
  \move (-10 30) \lvec (-5 45) \onedot \lvec (-10 65)
  \move (-5 45) \lvec (0 65)
  \move (0 17) \lvec (15 60) 
  \move (0 17) \lvec (22 35) \onedot 
  \htext (1 45) {$c$}
  \lpatt (2 2) \move (-12 29) \lcir r:11
    \esegment
    
\htext (240 15) {\normalsize  resulting in }

\move (340 0) 
   \scriptsize
   \bsegment 
   \htext (0 -20) {\normalsize $\ct(T)$}
   \move (0 0) \lvec (0 17) \onedot
   \lvec (-15 35) \onedot  
   \lvec (-35 70) \htext (-38 73) {$a$}
   \move (-15 35) \lvec (-30 40) \onedot
   \move (-15 35) \lvec (-25 80) \onedot \htext (-27 87) {$b$}
 \move (-15 35) 
 \lvec (-5 45) \onedot 
 \lvec (-5 60) \onedot \lvec (-5 75)  \onedot \lvec ( -10 90) 
 \move (-5 75) \lvec (0 90)

 \htext (-10 74) {$c$}
 \move (-5 45) \lvec (5 55) \onedot
 \move (0 17) \lvec (25 60) 
 \move (0 17) \lvec (22 35) \onedot 
 \lpatt (2 2)    \move (-12 47)  \freeEllipsis{26}{17}{32}

   \esegment

   \move (390 0)
 \end{texdraw}
 \end{equation}
    
    The new dots that appear in the composition tree of $T$ specify that
    new spheres should be drawn in $R_n$ in order to obtain $T_n$.  These
    spheres are drawn in the layer between the sphere $f$ and the spheres
    contained in $f$.  The dot substitution performed on the composition
    trees is not enough information though:
    there is an ambiguity for the spheres
    corresponding to the childless dots in $\ct(T)$: where should those
    null-spheres be drawn?  But the missing bit is clearly encoded in $S_n$
    itself.  In fact, substituting $\ct(S)$ into the $f$ node of $\ct(R)$
    is just the composition-tree expression of copying over the non-outer
    spheres from $S_n$ to $R_n$: copy those four spheres, and paste them
    into the layer between the sphere $f$ and its children.  The children
    of $f$ (dots and spheres immediately contained in $f$) are in 1--1
    correspondence with the dots in $S_n$ (since $F$ is the target of $S$
    and the $f$ source of $R$).  Here is the result, with the four new
    spheres highlighted in fat black:
\begin{equation}\label{fat-black}
  \begin{texdraw}
  \scriptsize \linewd 0.3
  \htext (-100 50) {\dots}
  \move (10 0)

  \bsegment 
  \move (-7 73) \freeFillEllipsis{45}{42}{-20}{0.9}
  \move (-7 73) \freeEllipsis{45}{42}{-20}
  \move (2.5 57.5) \freeFillEllipsis{14}{8}{70}{1}
  \move (2.5 57.5) \freeEllipsis{14}{8}{70}
  \move (-25 95) \fcir f:1 r:9 
  \move (-25 95) \lcir r:9
  \htext (0 -20) {\normalsize $T_n$}
  \move (0 -5) \lvec (0 20) \onedot
  \lvec (0 50) \onedot \lvec (-25 55) \onedot
  \lvec ( -75 70) \htext (-80 70) {$m$}
\move (-25 55) \lvec (-67 115)
  \htext ( -26 50) {$a$}
  \move (0 20) \lvec (60 45)
  \move ( 24 30) \linewd 0.3 \lcir r:5 
\linewd 0.5
\move (0 50) \lvec (5 65) \onedot 
  \lvec (40 125) 
  \move (5 65) \lvec (-55 125)  
  \htext (-59 128) {$k$}
  \htext (-28 93) {$b$}
	    \move (-2 68) \lcir r:60
\htext (40 58) {$f$}
\htext (0 59) {$c$}
\linewd 1
\move (2.5 57.5) \freeEllipsis{18}{12}{70}
\move (16 83) \lcir r:5
\move (8 66) \freeEllipsis{32}{20}{63}
\move (-37 72) \lcir r:5
  \esegment

  \move (170 0)
  \bsegment 
  \htext (0 -20) {\normalsize $\ct(T)$}
  \linewd 0.5 
  \move (0 0) \lvec (0 17) \onedot
  \lvec (-15 35) \onedot  
  \lvec (-35 70) \htext (-38 73) {$a$}
  \move (-15 35) \lvec (-30 40) \onedot
  \move (-15 35) \lvec (-25 80) \onedot \htext (-27 87) {$b$}
\move (-15 35) 
\lvec (-5 45) \onedot 
\lvec (-5 60) \onedot \lvec (-5 75)  \onedot \lvec ( -10 90) 
\move (-5 75) \lvec (0 90)

\htext (-10 74) {$c$}
\move (-5 45) \lvec (5 55) \onedot
\move (0 17) \lvec (25 60) 
\move (0 17) \lvec (22 35) \onedot 

  \esegment

\end{texdraw}
\end{equation}
\end{blanko}

\begin{blanko}{The filler.}\label{filler}
  The filling $(n+1)$-opetope $X$ should have $T$ as target, so $X_k =
  T_k$ for $k\leq n$.  The underlying tree of $X_{n+1}$ must be the
  composition tree of $T$; it remains to draw some spheres in this
  tree.  These spheres are determined by the original recipe
  composition tree (Figure~(\ref{recipe})): there are precisely two
  spheres to be drawn, corresponding to the two dots $S$ and $R$ in 
  the composition tree: one
  sphere is the outer sphere (corresponding to the root dot $R$), the
  other sphere is the `scar' of the gluing operation (corresponding to
  $S$) --- this sphere was already drawn dashed in 
  Figure~(\ref{subst-node}).
  
  So here is the final $X$ of our running example: 
  \begin{center}
    \begin{texdraw}
    \scriptsize \linewd 0.3
    \htext (-90 50) {\dots}
    \move (20 0)

    \bsegment 
    \htext (0 -20) {\normalsize $X_n = T_n$}
    \linewd 0.5
    \move (0 -5) \lvec (0 20) \onedot
    \lvec (0 50) \onedot \lvec (-25 55) \onedot
    \lvec ( -75 70) \htext (-80 70) {$m$}
  \move (-25 55) \lvec (-67 115)
    \htext ( -26 50) {$a$}
    \move (0 20) \lvec (60 45)
    \move ( 24 30) \linewd 0.3 \lcir r:5 
    \linewd 0.5
    \move (0 50) \lvec (5 65) \onedot 
    \lvec (40 125) 
    \move (5 65) \lvec (-55 125)  
    \htext (-59 128) {$k$}
    \move (-25 95) \linewd 0.3 \lcir r:9
    \htext (-28 93) {$b$}
    \move (2.5 57.5) \freeEllipsis{14}{8}{70}
	      \move (-2 68) \lcir r:60
	      \move (-7 73) \freeEllipsis{45}{42}{-20}
  \htext (40 58) {$f$}
  \htext (0 59) {$c$}
  \linewd 1
  \move (2.5 57.5) \freeEllipsis{18}{12}{70}
  \move (16 83) \lcir r:5
  \move (8 66) \freeEllipsis{32}{20}{63}
  \move (-37 72) \lcir r:5
    \esegment
    
          \htext (95 -18) {$\zoom$}

    \move (170 0)
    \bsegment 
    \htext (0 -20) {\normalsize $X_{n+1}$}
    \linewd 0.5
    \move (0 -5) \lvec (0 17) \onedot
    \lvec (-15 35) \onedot  
    \lvec (-45 80) \htext (-48 83) {$a$}
    \move (-15 35) \lvec (-30 40) \onedot
    \move (-15 35) \lvec (-25 75) \onedot \htext (-20 75) {$b$}
  \move (-15 35) 
  \lvec (-5 45) \onedot 
  \lvec (-5 60) \onedot \lvec (-5 75)  \onedot \lvec ( -13 100) 
  \move (-5 75) \lvec (3 100)

  \htext (-10 74) {$c$}
  \move (-5 45) \lvec (5 55) \onedot
  \move (0 17) \lvec (38 85) 
  \move (0 17) \lvec (22 35) \onedot 
  \move (-12 47)  \freeEllipsis{26}{17}{32}
  \move (0 50) \lcir r:45

    \esegment

    \move (300 10)
    \bsegment
      \move (0 0) \lvec (0 18) \onedot
      \move (0 18) \lvec (-6 36) \onedot
		    \lvec (-25 42) 
      \move (-6 36) \lvec (-19 52) 
      \move (-6 36) \lvec (-10 57) 
      \move (-6 36) \lvec (0 55) 
      \move (-6 36) \lvec (8 53) 
      \move (0 18) \lvec (16 48) 
      \move (0 18) \lvec (24 43) 
      \move (0 18) \lvec (28 35) 
      \move (0 18) \lvec (31 28) 
      
      \htext (0 -30) {\normalsize $\ct(T)$}
    \esegment

  \end{texdraw}
  \end{center}
  It is clear from the construction that it has $S$ and $R$ as
  sources and $T$ as target.
  
\end{blanko}



\renewcommand{\thesection}{\Alph{section}}

\section*{Appendix: Machine implementation}
\setcounter{section}{1}
\setcounter{lemma}{0}

%

Our description of opetopes naturally lends itself towards machine 
implementation.  The involved data grow only linearly with 
the dimension of the opetopes, and being fundamentally a tree 
structure, it is straightforward to encode in XML, as we
shall now explain.

\begin{blanko}{Trees-only representation.}
  For the sake of machine implementation, we have adopted a variation
  of the trees-only representation of opetopes given in 
  \ref{opetope-metatree}: instead of having the white dots (i.e.~the 
  null-spheres) explicitly, we let each null-dot refer to the unique 
  child of the corresponding null-sphere in the previous constellation
   (be it
  a dot or a leaf).
  Now, more than one null-sphere may sit on the same edge, in which
  case it is not enough for the corresponding null-dots to refer to
  that edge.  But the fact that these spheres sit on the same edge
  means there is induced an ordering among them, and this ordering can
  be expressed on the level of null-dots by letting them refer to each
  other in a chain, with only the last null-dot referring to something
  in the previous constellation (corresponding to the null-sphere
  farthest away from the root).  This system in turn requires some
  careful book-keeping in connection with sphere operations, since the
  reference of null-dot $x$ to a null-dot $y$ becomes invalid if $y$ is
  contracted.  Keeping track of these references is not difficult, but
  tedious and unenlightening.
\end{blanko}

\begin{blanko}{File format.}
  XML (Extensible Mark-up Language,
  cf.~\texttt{http://www.w3.org/XML/}) is a lot like HTML, except that
  you define your own tags to express a grammar.  This is done in a
  {\em Document Type Definition (DTD)}.  The opetope DTD looks like
  this:

\begin{lstlisting}
<!ELEMENT opetope (constellation+)>
<!ELEMENT constellation (dot|leaf)>
<!ELEMENT dot (dot|leaf)*>
<!ELEMENT leaf EMPTY>

<!ATTLIST opetope name CDATA #REQUIRED>
<!ATTLIST constellation name CDATA #REQUIRED>
<!ATTLIST dot name CDATA #REQUIRED ref CDATA #IMPLIED>
<!ATTLIST leaf name CDATA #REQUIRED>
\end{lstlisting}
  The first block declares the tags for \texttt{opetope},
  \texttt{constellation}, \texttt{dot}, and \texttt{leaf}, specifying
  which sort of children they can have. In the second block it is 
  specified that each tag must have a \texttt{name} attribute, and 
  that the
  \texttt{dot} tag is also allowed an optional \texttt{ref} attribute,
  used only for null-dots.
  
  Here is an XML representation 
  of the zoom complex in Example~\ref{opEx} interpreted as a 
  $5$-opetope:

\begin{lstlisting}
<?xml version="1.0" encoding="UTF-8"?>
<!DOCTYPE opetope SYSTEM "opetope.dtd">
<opetope name="Z">
    <constellation name="Z4">
	<dot name="b">
	    <dot name="a">
		<leaf name="1"/>
	    </dot>
	    <dot name="c">
		<leaf name="2"/>
		<leaf name="3"/>
	    </dot>
	</dot>
    </constellation>
    <constellation name="Z5">
	<dot name="p">
	    <dot name="x" ref="b"/>
	    <dot name="y">
		<leaf name="a"/>
		<leaf name="b"/>
		<leaf name="c"/>
	    </dot>
	</dot>
    </constellation>
    <constellation name="ct(Z)">
	<dot name="s">
	    <dot name="w" ref="a"/>
	    <leaf name="p"/>
	    <leaf name="x"/>
	    <leaf name="y"/>
	</dot>
    </constellation>
</opetope>
\end{lstlisting}
(The indentation is only for the benefit of the human reader; the
XML parser ignores whitespace between the tags.)  Notice how the 
null-dots $x$ and $w$ are provided with a reference to dots in the 
previous constellations, indicating where the corresponding spheres 
belong.
\end{blanko}

\begin{blanko}{Scripts.}
  The algorithms for sphere operations have been implemented in the
  scripting language Tcl, using the tDOM extension
  (cf.~\texttt{http://www.tdom.org/}) for parsing and manipulating XML.
  There are among other things procedures for computing sources, targets,
  and compositions, and writing the results back to new XML files.  These
  scripts can be run from the unix prompt, provided Tcl and the tDOM
  extension are available on the system.  The script 
  \texttt{computeAllFacets}
  takes as argument the name of an opetope XML file, and computes all its
  \mbox{codimension-$1$} faces, writing the resulting opetopes to separate
  XML files.  The script \texttt{glueOnto} takes three arguments: the
  bottom opetope (name of XML file), the name of the gluing locus, and the
  top opetope (as XML file).  The result is written to a new XML file.
  
  Precise instruction for installation and usage can be found in the readme
  file and manual pages accompanying the scripts.  XML files for all the
  examples of this paper are also included, together with the XML
  representation of a $10$-opetope with $15$ input facets.
\end{blanko}

\begin{blanko}{Automatic generation of graphical representation.}
  DOT\footnote{See {\sc E.~Gansner, E.~Koutsofios, 
  {\rm and } S.~North}, {\em Drawing graphs with {DOT}},
  \texttt{http://www.research.att.com/sw/tools/graphviz/dotguide.pdf}.} is a language for specifying abstract graphs in terms
  of node-edge incidences, and generate a graphical representation of
  the graph, for example in PDF format.
We provide a short Tcl script
  \texttt{opetope2pdf} which produces a dot file from an opetope XML
  file, and, if the dot interpreter is present on the system, also
  generates a pdf file.  This can be helpful to get an overview of a
  complicated opetope and its faces, but unfortunately the output is
  not quite as nice as the drawings in this paper (hand-coded \LaTeX);
  specifically, there is no support for drawing the spheres.
\end{blanko}


\noindent Here is what the output looks like when the script is run on
the XML file listed above:

\vspace{-10pt}

\begin{center}
\includegraphics[width=3in]{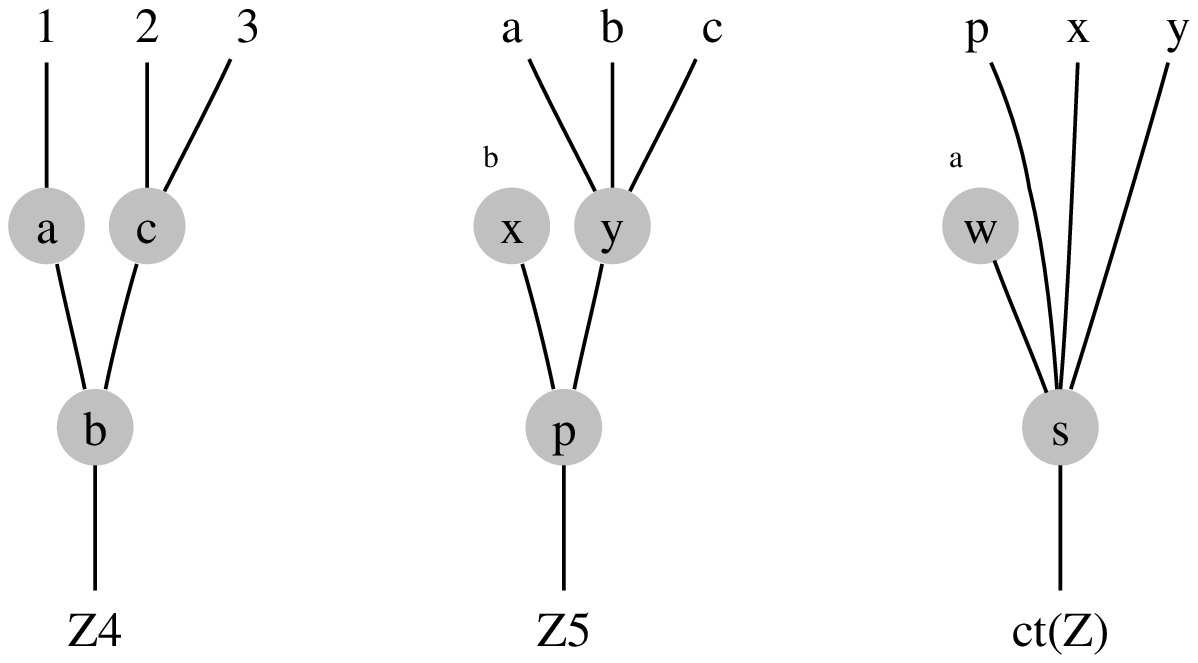}
\end{center}

\vfill

%
%

\medskip \small
\noindent 
{\sc Departament de Matem\`atiques -- Universitat Aut\`onoma de Barcelona
-- 08193 Bellaterra (Barcelona) -- Spain} 

\medskip \noindent
{\sc D\'epartement de math\'ematiques -- Universit\'e du Qu\'ebec \`a 
Montr\'eal -- Case postale 8888, succursale centre-ville --
Montr\'eal (Qu\'ebec), H3C 3P8 -- Canada}

\medskip \noindent
{\sc Department of Mathematics,
Division  of ICS --
Macquarie  University --
NSW  2109
-- Australia}

\medskip \noindent
{\sc Istituto per le Applicazioni del Calcolo and Institute of Complex Systems 
-- \linebreak National Research Council of Italy -- Via dei Taurini 19, 00185 Rome --
Italy}

\label{lastpage}

\end{document}